\documentclass[12pt]{article}
\usepackage{full page, amsfonts, tikz, comment,
amssymb, amscd, amsmath, graphicx, %showkeys,
wrapfig, stmaryrd, hyperref}
\usetikzlibrary{shapes}
\allowdisplaybreaks

 \title{{\bf Modular invariance of (logarithmic) intertwining operators}}
 \author{Yi-Zhi Huang}
    \date{}
    
\begin{document}
    \bibliographystyle{alpha}
    \maketitle
\newtheorem{thm}{Theorem}[section]
\newtheorem{defn}[thm]{Definition}
\newtheorem{prop}[thm]{Proposition}
\newtheorem{cor}[thm]{Corollary}
\newtheorem{lemma}[thm]{Lemma}
\newtheorem{rema}[thm]{Remark}
\newtheorem{app}[thm]{Application}
\newtheorem{prob}[thm]{Problem}
\newtheorem{conv}[thm]{Convention}
\newtheorem{conj}[thm]{Conjecture}
\newtheorem{cond}[thm]{Condition}
    \newtheorem{exam}[thm]{Example}
\newtheorem{assum}[thm]{Assumption}
     \newtheorem{nota}[thm]{Notation}
\newcommand{\halmos}{\rule{1ex}{1.4ex}}
\newcommand{\pfbox}{\hspace*{\fill}\mbox{$\halmos$}}

\newcommand{\bkdiamd}{\;\raisebox{.5pt}{\tikz{\node[draw,scale=0.35,diamond,fill=black](){};}}\;}

\newcommand{\nn}{\nonumber \\}

 \newcommand{\res}{\mbox{\rm Res}}
 \newcommand{\ord}{\mbox{\rm ord}}
\renewcommand{\hom}{\mbox{\rm Hom}}
\newcommand{\edo}{\mbox{\rm End}\ }
 \newcommand{\pf}{{\it Proof.}\hspace{2ex}}
 \newcommand{\epf}{\hspace*{\fill}\mbox{$\halmos$}}
 \newcommand{\epfv}{\hspace*{\fill}\mbox{$\halmos$}\vspace{1em}}
 \newcommand{\epfe}{\hspace{2em}\halmos}
\newcommand{\nord}{\mbox{\scriptsize ${\circ\atop\circ}$}}
\newcommand{\wt}{\mbox{\rm wt}\ }
\newcommand{\swt}{\mbox{\rm {\scriptsize wt}}\ }
\newcommand{\lwt}{\mbox{\rm wt}^{L}\;}
\newcommand{\rwt}{\mbox{\rm wt}^{R}\;}
\newcommand{\slwt}{\mbox{\rm {\scriptsize wt}}^{L}\,}
\newcommand{\srwt}{\mbox{\rm {\scriptsize wt}}^{R}\,}
\newcommand{\clr}{\mbox{\rm clr}\ }
\newcommand{\tr}{\mbox{\rm Tr}}
\renewcommand{\H}{\mathbb{H}}
\newcommand{\C}{\mathbb{C}}
\newcommand{\Z}{\mathbb{Z}}
\newcommand{\R}{\mathbb{R}}
\newcommand{\Q}{\mathbb{Q}}
\newcommand{\N}{\mathbb{N}}
\newcommand{\CN}{\mathcal{N}}
\newcommand{\F}{\mathcal{F}}
\newcommand{\I}{\mathcal{I}}
\newcommand{\V}{\mathcal{V}}
\newcommand{\one}{\mathbf{1}}
\newcommand{\BY}{\mathbb{Y}}
\newcommand{\ds}{\displaystyle}

        \newcommand{\ba}{\begin{array}}
        \newcommand{\ea}{\end{array}}
        \newcommand{\be}{\begin{equation}}
        \newcommand{\ee}{\end{equation}}
        \newcommand{\bea}{\begin{eqnarray}}
        \newcommand{\eea}{\end{eqnarray}}
         \newcommand{\lbar}{\bigg\vert}
        \newcommand{\p}{\partial}
        \newcommand{\dps}{\displaystyle}
        \newcommand{\bra}{\langle}
        \newcommand{\ket}{\rangle}

        \newcommand{\ob}{{\rm ob}\,}
        \renewcommand{\hom}{{\rm Hom}}

\newcommand{\A}{\mathcal{A}}
\newcommand{\Y}{\mathcal{Y}}
\newcommand{\End}{\mathrm{End}}
 \newcommand{\rad}{\mbox{\rm rad}}
\renewcommand{\l}{\llfloor}
\renewcommand{\r}{\rrfloor}

\newcommand{\dlt}[3]{#1 ^{-1}\delta \bigg( \frac{#2 #3 }{#1 }\bigg) }

\newcommand{\dlti}[3]{#1 \delta \bigg( \frac{#2 #3 }{#1 ^{-1}}\bigg) }

 \makeatletter
\newlength{\@pxlwd} \newlength{\@rulewd} \newlength{\@pxlht}
\catcode`.=\active \catcode`B=\active \catcode`:=\active
\catcode`|=\active
\def\sprite#1(#2,#3)[#4,#5]{
   \edef\@sprbox{\expandafter\@cdr\string#1\@nil @box}
   \expandafter\newsavebox\csname\@sprbox\endcsname
   \edef#1{\expandafter\usebox\csname\@sprbox\endcsname}
   \expandafter\setbox\csname\@sprbox\endcsname =\hbox\bgroup
   \vbox\bgroup
  \catcode`.=\active\catcode`B=\active\catcode`:=\active\catcode`|=\active
      \@pxlwd=#4 \divide\@pxlwd by #3 \@rulewd=\@pxlwd
      \@pxlht=#5 \divide\@pxlht by #2
      \def .{\hskip \@pxlwd \ignorespaces}
      \def B{\@ifnextchar B{\advance\@rulewd by \@pxlwd}{\vrule
         height \@pxlht width \@rulewd depth 0 pt \@rulewd=\@pxlwd}}
      \def :{\hbox\bgroup\vrule height \@pxlht width 0pt depth
0pt\ignorespaces}
      \def |{\vrule height \@pxlht width 0pt depth 0pt\egroup
         \prevdepth= -1000 pt}
   }
\def\endsprite{\egroup\egroup}
\catcode`.=12 \catcode`B=11 \catcode`:=12 \catcode`|=12\relax
\makeatother

\def\hboxtr{\FormOfHboxtr} % Only necessary if
%\kern... is wanted
\sprite{\FormOfHboxtr}(25,25)[0.5 em, 1.2 ex] % Resolution ca. 200x340 dpi.

:BBBBBBBBBBBBBBBBBBBBBBBBB | :BB......................B |
:B.B.....................B | :B..B....................B |
:B...B...................B | :B....B..................B |
:B.....B.................B | :B......B................B |
:B.......B...............B | :B........B..............B |
:B.........B.............B | :B..........B............B |
:B...........B...........B | :B............B..........B |
:B.............B.........B | :B..............B........B |
:B...............B.......B | :B................B......B |
:B.................B.....B | :B..................B....B |
:B...................B...B | :B....................B..B |
:B.....................B.B | :B......................BB |
:BBBBBBBBBBBBBBBBBBBBBBBBB |

\endsprite
\def\shboxtr{\FormOfShboxtr} % Only necessary if
%\kern... is wanted
\sprite{\FormOfShboxtr}(25,25)[0.3 em, 0.72 ex] % Resolution ca. 200x340 dp

:BBBBBBBBBBBBBBBBBBBBBBBBB | :BB......................B |
:B.B.....................B | :B..B....................B |
:B...B...................B | :B....B..................B |
:B.....B.................B | :B......B................B |
:B.......B...............B | :B........B..............B |
:B.........B.............B | :B..........B............B |
:B...........B...........B | :B............B..........B |
:B.............B.........B | :B..............B........B |
:B...............B.......B | :B................B......B |
:B.................B.....B | :B..................B....B |
:B...................B...B | :B....................B..B |
:B.....................B.B | :B......................BB |
:BBBBBBBBBBBBBBBBBBBBBBBBB |

\endsprite

\makeatletter
\DeclareFontFamily{U}{tipa}{}
\DeclareFontShape{U}{tipa}{m}{n}{<->tipa10}{}
\newcommand{\arc@char}{{\usefont{U}{tipa}{m}{n}\symbol{62}}}%

\newcommand{\arc}[1]{\mathpalette\arc@arc{#1}}

\newcommand{\arc@arc}[2]{%
  \sbox0{$\m@th#1#2$}%
  \vbox{
    \hbox{\resizebox{\wd0}{\height}{\arc@char}}
    \nointerlineskip
    \box0
  }%
}
\makeatother

\date{}
\maketitle

\begin{abstract}
Let $V$ be a 
$C_2$-cofinite vertex operator algebra without nonzero elements of negative weights. 
We prove the conjecture that the 
spaces spanned by analytic extensions of pseudo-$q$-traces
($q=e^{2\pi i\tau}$) shifted 
by $-\frac{c}{24}$ of products of geometrically-modified (logarithmic) 
intertwining operators among grading-restricted generalized 
$V$-modules are invariant under modular transformations.  
The convergence and analytic extension result needed to formulate 
this conjecture and some 
consequences on such shifted pseudo-$q$-traces were proved 
by Fiordalisi in \cite{F1} and \cite{F} using the method 
developed in \cite{H-mod-inv-int}. 
The method that we use to prove this conjecture 
is based on the theory of the associative 
algebras $A^{N}(V)$ for $N\in \N$, their graded modules 
and their bimodules introduced and studied by the 
author in \cite{H-aa-va} and \cite{H-aa-int-op}.
This modular invariance result gives a construction of
$C_2$-cofinite genus-one logarithmic conformal field theories
from the corresponding genus-zero logarithmic 
conformal field theories. 
\end{abstract}

%\tableofcontents \vspace{1em}
%\noindent{\large \bf References}\hfill %290

%\newpage

\renewcommand{\theequation}{\thesection.\arabic{equation}}
\renewcommand{\thethm}{\thesection.\arabic{thm}}
\setcounter{equation}{0} \setcounter{thm}{0} 

\section{Introduction}

\renewcommand{\theequation}{\thesection.\arabic{equation}}
\renewcommand{\thethm}{\thesection.\arabic{thm}}
\setcounter{equation}{0} \setcounter{thm}{0}

The modular invariance of (logarithmic) intertwining operators
was conjectured by the author almost twenty years ago\footnote{The conjecture 
appeared first in some grant proposals by the author and 
later in \cite{F1}.}. 
It is a conjecture on logarithmic conformal field theories.

In this paper, we prove this conjecture. 
In the language of conformal field theory, this modular invariance result 
says that in the $C_{2}$-cofinite case, genus-one logarithmic 
conformal field theories can be constructed by sewing the corresponding 
genus-zero logarithmic conformal field theories. 
We expect that this modular invariance result will play an important role
in the study of problems and conjectures on 
$C_2$-cofinite logarithmic conformal field theories.

Modular invariance plays a crucial role in 
the construction and study of
conformal field theories. 
In 1988, Moore and Seiberg \cite{MS1} \cite{MS2} 
conjectured that 
for rational conformal 
field theories, the space spanned by
$q$-traces ($q=e^{2\pi i\tau}$) shifted by $-\frac{c}{24}$ of 
products of $n$ chiral vertex operators (intertwining operators)
is invariant under modular transformations for each $n\in \N$. 
(For simplicity, we shall use “shifted (pseudo-)$q$-traces” 
below to mean “(pseudo-)$q$-traces
($q=e^{2\pi i\tau}$) shifted by $-\frac{c}{24}$.”)
Based on this conjecture in the most important case of $n=2$
and the conjecture that 
intertwining operators have operator product expansion, they 
derived the Verlinde formula \cite{V}
and discovered a mathematical structure now called 
modular tensor category (see \cite{T} for a precise definition and 
its connection with three-dimensional topological 
quantum field theories). 

In 1990, Zhu \cite{Z} 
proved a special case of the modular invariance 
conjecture of Moore and Seiberg under suitable 
conditions formulated precisely also in \cite{Z}. 
Let $V$ be a vertex operator algebra
satisfying the conditions that 
(i) $V$ has no nonzero elements 
of negative weights, 
(ii) $V$ is $C_{2}$-cofinite (that is, $\dim V/C_{2}(V)<\infty$,
where $C_{2}(V)=\langle \res_{x}x^{-2}Y_{V}(u, x)v\mid u, v\in V\rangle$), 
and (iii) every lower-bounded generalized $V$-module 
is completely reducible. Zhu proved in \cite{Z}
that for each $n\in \N$ and each subset $\{v_{1}, \dots, 
v_{n}\}$ of $V$, 
the space spanned by the analytic extensions of shifted
$q$-traces of 
products of suitably modified vertex operators associated to $v_{1}, \dots, v_{n}$
for $V$-modules
is invariant under modular 
transformations. In particular, the space spanned by the 
vacuum characters of irreducible $V$-modules is invariant 
under modular transformations. In 2000, using the 
the method developed by Zhu in \cite{Z}, Miyamoto in \cite{M1}
generalized this modular 
invariance result of Zhu to the space spanned
by the analytic extensions of shifted
$q$-traces  of 
products of suitably modified vertex operators associated to $v_{1}, \dots v_{n-1}$ 
for $V$-modules and 
one suitably modified
intertwining operator of a special type 
associated to an element $w\in W$ for each $n\in \Z_{+}$,
each subset $\{v_{1}, \dots, v_{n-1}\}$ of $V$, each $V$-module $W$ and 
each $w\in W$.
Unfortunately,  the method used in \cite{Z} and \cite{M1}
cannot be used to prove the modular invariance 
conjecture of Moore and Seiberg in the most important case
$n=2$ and the important cases $n>2$,

In 2002, Miyamoto \cite{M} 
proved a nonsemisimple generalization 
of Zhu's theorem. It was observed in \cite{M} 
that in the nonsemisimple case, the space of shifted $q$-traces of 
suitably modified vertex operators for 
grading-restricted generalized $V$-modules is in general 
not modular invariant and one needs shifted 
pseudo-$q$-traces of suitable operators on 
grading-restricted generalized $V$-modules introduced and studied in
\cite{M}.
Let $V$ be a vertex operator algebra satisfying 
Conditions (i) and (ii) above and also satisfying 
an additional condition that (iv) there are no finite-dimensional
irreducible $V$-modules. 
Miyamoto proved that for each $n\in \N$ and each subset 
$\{v_{1}, \dots, v_{n}\}$ of $V$, the space spanned by the 
analytic extensions of shifted
pseudo-$q$-traces  of 
products of suitably modified vertex operators associated to $v_{1}, \dots,v_{n}$
for grading-restricted 
generalized  $V$-modules
is invariant under modular 
transformations. The fact that Condition (iv) above is needed in \cite{M}
was pointed out explicitly in Remark 3.3.5 in \cite{ArN}. 
There are examples of vertex operator algebras satisfying 
Conditions (i) and (ii) above but not Condition (iv) (for example,
$\mathcal{W}_{2, 3}$; see Remark 3.3.5 
in \cite{ArN} and \cite{AM}). 

In 2003, the author proved in \cite{H-mod-inv-int} the 
modular invariance 
conjecture of Moore and Seiberg under the same conditions 
as in \cite{Z}. The precise statement of the modular invariance 
theorem in \cite{H-mod-inv-int} is that for a vertex operator algebra
$V$ satisfying the same conditions (i)\footnote{Note that the statement of the modular
 invariance theorem  in \cite{H-mod-inv-int} also 
has an additional condition $V_{(0)}=\C\one$. This condition is 
added in \cite{H-mod-inv-int} because Theorem 7.2 in 
\cite{H-mod-inv-int} needs some results of \cite{ABD}, which
are in turn proved using a result of Buhl \cite{B} giving a spanning set of 
a weak $V$-module. It is this result in \cite{B} 
that needs the condition $V_{(0)}=\C\one$. 
But in Lemma 2.4 in \cite{M}, Miyamoto obtained such a spanning set of a weak module 
without using the condition $V_{(0)}=\C\one$. So the modular invariance theorem 
in \cite{H-mod-inv-int} in fact does not need this condition.}, (ii) and 
(iii) above as in \cite{Z} 
and for each $n\in 
\N$, each set of $n$ grading-restricted  generalized $V$-modules
$W_{1}, \dots, W_{n}$  and each set $\{w_{1}, \dots, w_{n}\}$
for $w_{1}\in W_{1}, \dots, w_{n}\in W_{n}$, 
the space spanned by the analytic extensions of shifted
$q$-traces of 
products of geometrically-modified 
intertwining operators associated to $w_{1}, \dots, w_{n}$
is invariant under modular transformations. 
As is mentioned above, the method used by Zhu and Miyamoto
cannot be used to prove this conjecture of Moore and Seiberg
in the main important cases $n\ge 2$. 
A completely different method was developed in 
\cite{H-mod-inv-int}
to prove this conjecture in these cases, including the most important case
$n=2$. 
Using this modular invariance 
and the associativity of intertwining operators
proved in \cite{H-diff-eqn}, the author proved 
the Verlinde conjecture and
Verlinde formula in \cite{H-verlinde}
and the rigidity and modularity 
of the braided tensor category of $V$-modules in \cite{H-rigidity}. 

Around 2003, after the work of Miyamoto \cite{M} and 
the proof in \cite{H-mod-inv-int}
of the modular invariance conjecture
of Moore and Seiberg, the following modular invariance 
was conjectured by the author:
For a $C_{2}$-cofinite 
vertex operator algebra $V$ without nonzero elements of negative weights
and for each $n\in \N$,
each set of $n$ grading-restricted generalized $V$-modules 
$W_{1}, \dots, W_{n}$ and each set $\{w_{1}, \dots, w_{n}\}$
for $w_{1}\in W_{1}, \dots, w_{n}\in W_{n}$, 
the space of analytic extensions of shifted
pseudo-$q$-traces of 
products of geometrically-modified 
(logarithmic) intertwining operators associated to 
$w_{1}, \dots, w_{n}$ 
is invariant under modular 
transformations. We put the word ``logarithmic'' in parenthesis before ``intertwining 
operators" since both intertwining operators without logarithm
and  logarithmic intertwining operators are needed. In the main body
of this paper, we shall omit ``(logarithm)" so that intertwining 
operators in general might have the logarithms of the variables involved.

In fact, in the original version of this conjecture, the convergence 
and analytic extension of shifted pseudo-$q$-traces  of 
products of geometrically-modified 
(logarithmic) intertwining operators is also part of the conjecture.
In 2015, using the method developed in 
\cite{H-mod-inv-int}, Fiordalisi \cite{F1} \cite{F} proved 
that such shifted pseudo-$q$-traces are 
convergent absolutely in a suitable region and can be analytically
extended to multivalued analytic functions on a maximal region. 
He also proved that these multivalued analytic functions 
satisfy the genus-one associativity, genus-one commutativity 
and other properties using the associativity  and commutativity 
of (logarithmic) intertwining operators and other properties for a $C_{2}$-cofinite
vertex operator algebra proved in \cite{H-cofiniteness}. 
In addition, he proved in \cite{F1} \cite{F}
that the solution space of the differential equations
used to prove the convergence and analytic extension property of
shifted pseudo-$q$-traces discussed above is invariant under 
the modular transformations. 
As in the proof in \cite{H-mod-inv-int} of the modular invariance conjecture of 
Moore and Seiberg, by the genus-one associativity proved in \cite{F1} and \cite{F}, 
the modular invariance conjecture in the case $n\ge 2$ 
can be reduced to the modular invariance conjecture in the case $n=1$.
Here the author would like to emphasize the importance of the 
convergence and analytic extension results proved in 
\cite{H-mod-inv-int}, \cite{F1} and \cite{F}. Without the convergence 
and analytic extension 
results in \cite{H-mod-inv-int}, \cite{F1} and \cite{F}, we could not
even formulate the modular invariance conjecture of Moore and Seiberg
and the modular invariance conjecture for (logarithmic) intertwining
operators above. See \cite{H-conv-cft} for a survey on 
the convergence and analytic extension results and conjectures
in the approach to conformal field theory using the representation theory of vertex 
operator algebras. 

To prove the modular invariance conjecture for (logarithmic) intertwining
operators above, 
we need to show that the modular transformations of 
the multivalued analytic 
functions obtained from shifted
pseudo-$q$-traces of geometrically-modified
(logarithmic) intertwining 
operators are sums of the multivalued analytic 
functions obtained from such shifted pseudo-$q$-traces.
In \cite{Z}, \cite{M1} and \cite{H-mod-inv-int},  
Zhu algebra $A(V)=A_{0}(V)$ associated to $V$, 
its modules and its bimodules are needed 
to prove the modular invariance. 
In \cite{M}, the 
generalizations $A_{n}(V)$ for $n\in \N$ of Zhu algebra by Dong, Li and Mason 
\cite{DLM} and their modules are needed to prove 
the modular invariance. On the other hand, 
the additional condition (Condition (iv) above) in \cite{M} 
that there are no finite-dimensional
irreducible $V$-modules
is needed exactly because the associative algebras $A_{n}(V)$ for $n\in \N$
cannot be used to handle the case that there 
exist finite-dimensional
irreducible $V$-modules. More importantly, to study general (logarithmic)
intertwining operators using the theory of associative algebras, 
the associative algebras $A_{n}(V)$ for $n\in \N$
are not enough. Therefore, even if Condition (iv) above is satisfied, 
we will not be able to prove this modular invariance 
using only $A_{n}(V)$ for $n\in \N$ and their bimodules 
introduced and studied in \cite{HY}. 

In \cite{H-aa-va}, the author 
introduced new associative algebras $A^{\infty}(V)$ and $A^{N}(V)$ 
for $N\in \N$ associated 
to a vertex operator 
algebra $V$. These associative algebras 
contain $A_{n}(V)$ for $n\in \N$ as (very small) subalgebras. 
In fact, $A_{n}(V)$ for $n\in \N$ are all algebras of zero modes
but acting on different homogeneous subspaces of lower-bounded generalized 
$V$-modules.
On the other hand, the associative algebras $A^{\infty}(V)$ and 
$A^{N}(V)$ for $N\in \N$ introduced by
the author in \cite{H-aa-va} are algebras of all modes, including all nonzero modes.
In \cite{H-aa-int-op},
the author introduced bimodules $A^{\infty}(W)$ and 
$A^{N}(W)$ for $N\in \N$ for these new associative 
algebras associated to a lower-bounded generalized $V$-module $W$
and proved that the spaces of (logarithmic) intertwining 
operators are linearly isomorphic to the corresponding spaces 
of module maps between suitable modules for 
these associative algebras. In Section 3 of
the present paper, we also introduced associative algebras 
$\widetilde{A}^{\infty}(V)$ and $\widetilde{A}^{N}(V)$ for $N\in \N$
isomorphic to $A^{\infty}(V)$ and 
$A^{N}(V)$ for $N\in \N$, respectively, and the corresponding bimodules
$\widetilde{A}^{\infty}(W)$ and $\widetilde{A}^{N}(W)$ for $N\in \N$ 
associated to a lower-bounded generalized $V$-module $W$. We
then transport the results on the associative algebras $A^{\infty}(V)$ and 
$A^{N}(V)$ for $N\in \N$, their
graded modules and bimodules obtained in \cite{H-aa-va} and 
\cite{H-aa-int-op} to the corresponding results on
$\widetilde{A}^{\infty}(V)$ and 
$\widetilde{A}^{N}(V)$ for $N\in \N$, their graded modules and bimodules.

In this paper, using the results of Fiordalisi in \cite{F1} and \cite{F}
and the results on these new associative algebras, their modules
and bimodules in \cite{H-aa-va}, \cite{H-aa-int-op} and Section 3 of the present paper, 
we prove the modular invariance conjecture for
(logarithmic) intertwining operators discussed above. 
For the precise statement, see Theorem \ref{mod-inv}. 
Note that all the modular invariance theorems mentioned above 
are special cases of Theorem \ref{mod-inv}. In particular, 
in the special case studied in \cite{M} that the intertwining 
operators involved are vertex operators for grading-restricted 
generalized $V$-modules, we obtain a proof of the modular invariance result 
in \cite{M} without Condition (iv) above requiring that there are no 
finite-dimensional irreducible $V$-modules. Also, as special cases 
of the proofs of (\ref{main-lemma-3}) and (\ref{grading-slf}), 
we recover the proof by McRae \cite{Mc} of Propositions 4.4 in \cite{M} and 
obtain a proof of Proposition 4.5 in \cite{M} (see 
Remarks \ref{miyamoto-prop-4.4}
and \ref{miyamoto-prop-4.5}). 

Here we give a sketch of the proof of this 
modular invariance conjecture: As is mentioned above,
we need only prove the modular invariance conjecture 
in the case $n=1$.  To prove the
modular invariance conjecture in this case, 
we introduce a notion of genus-one $1$-point 
conformal block  labeled by a grading-restricted generalized 
$V$-module $W$. Then the conjecture follows
from the following two results:
(1) For grading-restricted generalized $V$-modules $W$ and 
$\widetilde{W}$, the modular 
transformation of the analytic extension of the shifted 
pseudo-$q$-trace of a geometrically-modified intertwining 
operator of type $\binom{\widetilde{W}}{W\widetilde{W}}$
is a genus-one $1$-point 
conformal block labeled by $W$ 
(see Proposition \ref{mod-transf-conf-bk}). (2) Every 
genus-one $1$-point 
conformal block  labeled by $W$
is the sum of the analytic extensions of the shifted 
pseudo-$q$-traces of geometrically-modified intertwining 
operators of type 
$\binom{\widetilde{W}_{i}}{W\widetilde{W}_{i}}$
for finitely many grading-restricted generalized $V$-modules 
$\widetilde{W}_{i}$ (see Theorem \ref{conf-bk=>pseudo-q-tr}). 
The proof of (1) can be obtained
easily from the properties of the shifted pseudo-q-traces of 
(logarithmic) intertwining operators
in \cite{F1} and \cite{F}.
To prove (2), we first prove 
that a genus-one $1$-point 
conformal block labeled by $W$
gives a symmetric linear function
on $\widetilde{A}^{N}(W)$ satisfying some additional properties
for each $N\in \N$. This proof 
is technically the most difficult part of this paper. 
Then using the results on symmetric 
linear functions and pseudo-traces 
proved by Miyamoto \cite{M}, Arike \cite{Ar}
and Fiordalisi \cite{F1} \cite{F}, 
we prove  that for $N$ sufficiently large,
these symmetric linear functions 
on $A^{N}(W)$ are in fact finite sums of pseudo-traces 
of suitable linear operators on $A^{N}(V)$-modules. 
Next we use a main result in \cite{H-aa-int-op} 
to show that for $N$ sufficiently large, pseudo-traces 
of such linear operators on $A^{N}(V)$-modules
are in fact obtained from shifted pseudo-$q$-traces of geometrically-modified 
(logarithmic) intertwining operators. 
Finally we show that when $N$ is sufficiently large,
the genus-one $1$-point 
conformal block labeled by $W$ that we start with 
is equal to the sum of the analytic extensions
of the shifted pseudo-$q$-traces of these geometrically-modified
(logarithmic) intertwining operators. 

The present paper is organized as follows:
In Section 2, we recall some basic definitions and results
needed in this paper. In Subsection 2.1, we recall 
the definition of pseudo-traces and the results about 
pseudo-traces and symmetric linear functions obtained 
by  Miyaomoto \cite{M}, Arike \cite{Ar} and Fiordalisi \cite{F1}. 
In Subsection 2.2, we recall the results of 
Fiordalisi  in \cite{F1} and \cite{F} on 
the convergence and analytic extensions of shifted pseudo-$q$-traces of 
products of intertwining 
operators and their properties. In Section 3, for a vertex operator
algebra $V$, we study further the $A^{\infty}(V)$-bimodule $A^{\infty}(W)$ and 
the $A^{N}(V)$-bimodules $A^{N}(W)$ for $N\in \N$
constructed from a lower-bounded generalized $V$-module $W$
in \cite{H-aa-int-op}. Then
we introduce and study new associative algebras $\widetilde{A}^{\infty}(V)$ 
and $\widetilde{A}^{N}(V)$
$N\in \N$, the $\widetilde{A}^{\infty}(V)$-bimodule 
$\widetilde{A}^{\infty}(W)$ and 
$\widetilde{A}^{N}(V)$-bimodules 
$\widetilde{A}^{N}(W)$ for $N\in \N$ in this section.
In Section 4, we give two constructions of 
symmetric linear functions on $\widetilde{A}^{N}(V)$-bimodules 
$\widetilde{A}^{N}(W)$ for a grading-restricted generalized 
$V$-module $W$. The first construction is given by shifted 
pseudo-$q$-traces of intertwining operators. 
The second construction is given by suitable maps 
satisfying properties involving 
the Weierstrass $\wp$- and $\zeta$-functions. 
In Section 5, we prove our modular invariance 
theorem (Theorem \ref{mod-inv})
by proving Proposition \ref{mod-transf-conf-bk}
and Theorem \ref{conf-bk=>pseudo-q-tr} (the two results (1)
 and (2) discussed above). We also have two appendices. 
 In Appendix A, we recall some basic facts on  
 the Weierstrass $\wp$-function
$\wp_{2}(z; \tau)$
the Weierstrass $\zeta$-function $\wp_{1}(z; \tau)$
and the Eisenstein series $G_{2}(\tau)$. In Appendix B, 
we prove and collect a number of identities involving 
the binomial coefficients. The proof of the modular invariance 
theorem in this paper depends heavily on these identities.

\paragraph{Acknowledgment} Francesco Fiordalisi 
worked jointly with me in 2015 on an early method to the 
modular invariance conjecture of (logarithmic) 
intertwining operators.
I am very grateful to him for many discussions 
when he did his Ph.D. thesis research in \cite{F1} and \cite{F}. I am also very 
grateful to Robert McRae for many discussions, especially  
on his  proof in \cite{Mc} of Proposition 4.4 in \cite{M} (see 
Remark \ref{miyamoto-prop-4.4}), his observation that 
the condition $V_{(0)}=\C\one$ is not needed in the proof 
that a grading-restricted generalized $V$-module is $C_{n}$-cofinite 
(see the proof of Theorem \ref{finite-d-A-N-W}) 
and the possible applications of 
the modular invariance theorem proved in this paper. 
Some adjustments in the statements and proofs of 
several results and the corrections of many typos
in the present paper were done when I was visiting the School of Mathematical 
Sciences in Shanghai Jiaotong University. I would like to thank Jinwei 
Yang for his hospitality during the visit. I am also grateful to Hao Zhang 
for noticing a mistake in the formulation of Theorem \ref{ar-sym}
of Miyamoto and Arike in the previous version of 
this paper.  (This wrong formulation and the discussions
in several places using this formulation 
have been corrected in this version).

\renewcommand{\theequation}{\thesection.\arabic{equation}}
\renewcommand{\thethm}{\thesection.\arabic{thm}}
\setcounter{equation}{0} \setcounter{thm}{0}

\section{Pseudo-traces, symmetric linear functions and 
pseu\-do-$q$-traces}

We recall some basic definitions and results in this section. 
In Subsection 2.1, we recall the basic definitions and results on 
pseudo-traces and symmetric linear functions 
introduced and obtained by Miyamoto \cite{M}, Arike \cite{Ar} and 
Fiordalisi \cite{F1}.
In Subsection 2.2, we recall the results obtained by Fiordalisi \cite{F1} \cite{F}
on genus-one correlation functions
constructed using shifted pseudo-$q$-traces of products of (logarithmic) 
intertwining operators. 

\subsection{Pseudo-traces and symmetric linear functions}

In this subsection, we first  recall the definition of pseudo-traces 
for a finitely generated right projective module
$M$ for a finite-dimensional associative algebra $A$ introduced by 
Miyamoto \cite{M} and further studied
by Arike \cite{Ar}. Then we 
recall a result on symmetric linear functions 
on $A$ obtained by Miyamoto \cite{M} and Arike \cite{Ar}
and a result on symmetric linear
functions on finite-dimensional $A$-bimodules 
obtained by Fiordalisi \cite{F1}. See \cite{Ar} and Fiordalisi \cite{F1} 
for details.

Let $A$ be a finite-dimensional associative algebra.
Recall that a right $A$-module $M$ is said to be projective if every short exact sequence 
$$0\to K\to P\to M\to 0$$
for right $A$-modules $K$ and $P$ splits. 

Let $M$ be a finitely generated right $A$-module. A projective basis for $M$ is
a pair of sets $\{m_i\}_{i=1}^n\subset M$,
$\{\alpha_i\}_{i=1}^n\subset \hom_A(M, A)$
such that for all $m\in M$, 
$$m = \sum_{i=1}^n m_i\alpha_{i}(m).$$ 
A finitely generated right $A$-module $M$ has a projective basis if and only if it is projective.

Let $M$ be a right $A$-module. Then
$\hom_A(M,A)$ has a left $A$-module structure given by
$$(a\alpha)(m) = a\alpha(m)$$
for $a\in A, \alpha\in\hom_A(M,A)$ and $m\in M$.
We also have a contraction  map
\begin{align*}
  \pi_{M}: M\otimes_A \hom_A(M,A) &\rightarrow A \\
  m\otimes_{A} \alpha & \mapsto \alpha(m)
\end{align*}

For right $A$-modules $M_{1}$ and $M_{2}$,
let
\begin{align*}
\tau_{M_{1}, M_{2}}: M_{2}\otimes_A \hom_A(M_{1},A) &\to \hom_A(M_{1}, M_{2})\nn
m_{2}\otimes_{A} \alpha&\mapsto \tau_{M_{1}, M_{2}} (m_{2}\otimes_{A} \alpha)
\end{align*}
be the natural linear map defined by $(\tau_{M_{1}, M_{2}} (m_2\otimes \alpha))(m_{1})=m_{2}\alpha(m_{1})$ for
$m_{1}\in M_{1}$, $m_{2}\in M_{2}$ and $\alpha\in \hom_A(M_{1}, M_{2})$.
In the case that $M_{1}=M_{2}=M$, let $\tau_M = \tau_{M, M}$.

For  a finitely generated right $A$-module $M$, the map $\tau_M: M\otimes_A \hom_A(M,A)
 \rightarrow \End_{A}M$ is an isomorphism so that $\tau_M^{-1}: \End_{A}M
\to M\otimes_A \hom_A(M,A)$ exists. The
Hattori-Stallings trace of an endomorphism
$\alpha\in \End_{A}M$ is the element
$$
\tr_M \alpha= \pi_{M}(\tau_M^{-1}(\alpha)) + [A,A]
$$
of $A/[A,A]$. 
For finitely generated  projective right $A$-modules $M_{1}$, $M_{2}$, 
$f\in \hom_A(M_{1}, M_{2})$ and $g\in \hom_A(M_{2}, M_{1})$, we have
$$ \tr_{M_{2}} f\circ g=\tr_{M_{1}} g\circ f.$$

A linear function $\phi: A\to \C$ is said to be symmetric
if  $\phi(ab) = \phi(ba)$ for all $a,b\in A$. We denote the the
space of symmetric linear functions on $A$
by $SLF(A)$. Then $SLF(A)$ is linearly isomorphic to $(A/[A,A])^*$.
Any symmetric linear function on $A$ defines a
symmetric bilinear form
\begin{align*}
\langle\ ,~ \rangle &: A\times A \rightarrow \C
\end{align*}
by $\langle a, b\rangle = \phi (ab)$. A symmetric linear function $\phi$  
on $A$ is said to be 
nondegenerate if the corresponding bilinear form is nondegenerate.
The radical of a symmetric linear function on $A$ is defined to 
be the two sided ideal 
$$\rad~\phi = 
\big\{a\in A\ |\ \langle a, b\rangle =0 \quad \text{for all } b\in A\big\}$$
of $A$. A symmetric function $\phi$ is nondegenerate if
and only if $\rad\ \phi = \{0\}$.

Let $M$ be a finitely generated projective right $A$-module and
$\{m_i\}_{i=1}^n, \{\alpha_i\}_{i=1}^n$ a projective basis.
The pseudo-trace function $\phi_M$ on $\End_{A}M$ associated to 
a linear symmetric function $\phi$ on $A$
is the map $\phi_M = \phi\circ \tr_M: \End_{A}M\to \C$. 
We can express the pseudo-trace
of  $\alpha\in \End_{A}M$ in terms of the projective basis as
$$\phi_{M}(\alpha)=\phi\left(\sum_{i=1}^n \alpha_i(\alpha(m_i))\right).$$

For right projective $A$-modules $M_1$ and $M_2$,
  $\alpha \in \hom_A(M_1,M_2)$ and $\beta \in \hom_A(M_2,M_1)$, we have
  $$
    \phi_{M_1} (\beta \circ \alpha) =
    \phi_{M_2} (\alpha \circ \beta).
  $$

  A symmetric algebra (or Frobenius algebra) is an
  associative algebra equipped with a nondegenerate symmetric linear function.  
A basic algebra is an associative algebra $A$ such that $A/J(A)$ is isomorphic
  to $\C^n$ for some $n\in \N$, where $J(A)$ is the Jacobson radical
  of $A$. For an associative algebra $A$, we have a complete set of 
  orthogonal primitive idempodents. Each idempodent in this set gives
  a right $A$-module of $A$. But these right $A$-modules 
  might be equivalent. Then we can find a complete set 
$\{e_{ij}\mid i=1, \dots, n, \; j=1, \dots, n_{i}\}$ of
  orthogonal primitive idempodents such that 
  the right $A$-modules $e_{ik}A\simeq e_{il}A$ for $i$, $k, l=1, \dots, n_{i}$ 
  and $e_{ik}A\not\simeq e_{jl}A$ for $i\ne j$, $k=1, \dots, n_{i}$,
  $l=1, \dots, n_{j}$, where $\simeq$ means equvalence of 
  right $A$-modules. 
  Since this set is complete, we have 
  $$1_{A}=\sum_{i=1}^{n}
\sum_{j=1}^{n_{i}}e_{ij}.$$
  Let $e=e_{11}+\cdots +e_{n1}$. The algebra  $eAe$
  is called a basic algebra of $A$.
  
\begin{thm}[Miyamoto \cite{M}, Arike  \cite{Ar}] \label{ar-sym}
  Let $A$ be a finite-dimensional associative algebra and $\phi\in SLF(A)$.
  Let 
$A=A_{1}\oplus \cdots \oplus A_{n}$ be a decomposition of $A$
as a direct sum of indecomposable $A$-bimodules. 
Let $\phi_{i}=\phi|_{A_{i}}$
for $i=1, \dots, n$. For each $i$, let $P_{i}=\bar{e}_{i}
(A/{\rm Rad}(\phi_{i}))\bar{e}_{i}$ and 
$M_{i}=(A/{\rm Rad}(\phi_{i}))\bar{e}_{i}$ be a basic algebra 
of $A/{\rm Rad}(\phi_{i})$ and the corresponding $A$-$P_{i}$-bimodule, 
where $\bar{e}_{i}=\bar{e}_{i, 11}+\cdots +\bar{e}_{i, n_{i}1}$
and 
$$\{\bar{e}_{i, jk}\mid j=1, \dots, n_{i}, k=1, \dots, n_{ij}\}$$
is a complete set of orthogonal primitive idempodents 
of $A/{\rm Rad}(\phi_{i})$ such that 
$$1_{A/{\rm Rad}(\phi_{i})}=\sum_{j=1}^{n_{i}}
\sum_{k=1}^{n_{ij}}\bar{e}_{i, jk},$$
$\bar{e}_{i,jl}(A/{\rm Rad}(\phi_{i}))\simeq 
\bar{e}_{i,jm}(A/{\rm Rad}(\phi_{i}))$
and $\bar{e}_{i,jl}(A/{\rm Rad}(\phi_{i}))\simeq 
\bar{e}_{i,km}(A/{\rm Rad}(\phi_{i}))$ for $l\ne m$.
Then for $i=1, \dots, n$, 
$P_i$ are basic symmetric algebras with symmetric linear functions
given by  $\phi_i$ (still denoted by $\phi_{i}$), 
and $M_i$ are $A$-$P_i$-bimodules ,
   finitely generated and projective as right $P_i$-modules,
  such that
  $$
    \phi(a) = \sum_{i=1}^n (\phi_i)_{M_i}(a)
  $$
  where in each term in the right-hand side, 
$a\in A$ is viewed 
an element of $\End_{P_i}M_i$ given by the left action $a$ on $M_{i}$.
  Furthermore, if $\nu$ is an element of $\rad\ \phi$, that is,
  $$
    \phi(\nu a) = 0
  $$
  for all $a\in A$, then 
  $\nu M_i = 0$.
\end{thm}

 Let $A$ be an associative algebra and $M$ an $A$-bimodule. A linear
 function $\phi:M\rightarrow \C$ is said to be symmetric if for all
  $m\in M$ and $a\in A$,
  $$\phi(am) = \phi(ma).$$
Let $P$ be another associative algebra and 
$U$ an $A$-$P$-bimodule. 
Then the endomorphism ring $\End_{P} U$ is an $A$-$A$-bimodule with the actions given by $(a\tau)(u) = a(\tau(u))$ and $(\tau a)(u) = \tau(au)$
for $a\in A$, $\tau\in \End_{P} U$ and $u\in U$.

Now assume that $P$ is finite dimensional and $U$ is finitely 
generated and projective as a right
$P$-module. Then for  $\phi\in SLF(P)$, the pseudo-trace $\phi_U$ on $\End_P U$ 
is a symmetric linear function on the $A$-bimodule $\End_P U$.

Let $M$ be an $A$-bimodule and let $\hom_{A,P}(M\otimes_A U, U)$ be the
set of all $A$-$P$-bimodule maps from $M\otimes_A U$ to $U$. Let
$f\in \hom_{A,P}(M\otimes_A U, U)$. Then for any $m\in M$, the map
\begin{align*}
  U&\rightarrow U\\
  u&\mapsto f(m\otimes u)
\end{align*}
is an element of  $\End_P U$. Define $T_f : M \rightarrow \End_{P} U$ by 
$T_{f}(m)=(u \mapsto f(m\otimes u))$.
Then the map $T_f$ is an $A$-bimodule homomorphism.
Let $\phi^f_U:M \rightarrow \C$ be the linear function on $M$ defined by
$\phi^f_U(m)=\phi_U(T_f(m))$.
Then the linear function $\phi^f_U$ is symmetric.

We have obtained a map $SLF(P)\otimes\hom_{A,P}(M\otimes_AU, U) \rightarrow SLF(M)$. We now want to give an ``inverse'' map in a suitable sense. 

We consider the trivial square-zero extension $\bar{A}=A\oplus M$
of $A$ by $M$ with the product given by 
$$(a_{1}, m_{1})(a_{2}, m_{2})
=(a_{1}a_{2}, a_{1}m_{2}+m_{1}a_{2})$$
for $(a_{1}, m_{1}), (a_{2}, m_{2})\in \bar{A}$. 
Let $\phi\in SLF(M)$. We extend $\phi$ 
to a linear function $\bar{\phi}$ on $\bar{A}$ by 
$\bar{\phi}(a, m)=\phi(m)$. Then 
\begin{align*}
\bar{\phi}((a_{1}, m_{1})(a_{2}, m_{2}))
&=\bar{\phi}(a_{1}a_{2}, a_{1}m_{2}+m_{1}a_{2})\nn
&=\phi(a_{1}m_{2}+m_{1}a_{2})\nn
&=\phi(m_{2}a_{1}+a_{2}m_{1})\nn
&=\bar{\phi}(a_{2}a_{1}, m_{2}a_{1}+a_{2}m_{1})\nn
&=\bar{\phi}((a_{2}, m_{2})(a_{1}, m_{1})),
\end{align*}
which says that $\bar{\phi}$ is in fact symmetric. 

Let $1_{A}=d_{1}+\cdots +d_{n}$, where $1_{A}$ 
is the identity of $A$ and $d_{1}, \dots, d_{n}$
are orthogonal primitive central idempotents. Then 
$A=A_{1}\oplus \cdots \oplus A_{n}$ be a decomposition of $A$
as a direct sum of indecomposable $A$-bimodules, where $A_{i}=Ae_{i}$ for 
$i=1, \dots, n$. Then
we also have a decomposition 
$\bar{A}=A+M=\bar{A}_{1}\oplus \cdots
\oplus \bar{A}_{n}$ as a direct sum of indecomposable $A$-bimodules, where 
$\bar{A}_{i}=\bar{A}e_{i}$ for 
$i=1, \dots, n$. 
Let $\bar{\phi}_{i}=\bar{\phi}|_{\bar{A}_{i}}$
for $i=1, \dots, n$. 

\begin{thm}[Fiordalisi \cite{F1}] \label{slf=>p-tr}
  Let $A$ be a finite-dimensional associative algebra and 
$M$ a finite-dimensional $A$-bimodule, and let $\phi\in SLF(M)$. 
Let $1_{A}=d_{1}+\cdots +d_{n}$, where $1_{A}$ 
is the identity of $A$ and $d_{1}, \dots, d_{n}$
are orthogonal primitive central idempotents such that 
$A=A_{1}\oplus \cdots \oplus A_{n}$ is a decomposition of $A$
as a direct sum of $A$-bimodules, where $A_{i}=Ad_{i}$ for 
$i=1, \dots, n$,
 and $\bar{A}=A+M=\bar{A}_{1}\oplus \cdots
\oplus \bar{A}_{n}$ is a decomposition of $\bar{A}$, where 
$\bar{A}_{i}=\bar{A}d_{i}$ for 
$i=1, \dots, n$. For each $i$, let 
$\bar{P}_{i}=\bar{e}_{i}(\bar{A}/{\rm Rad}(\bar{\phi}_{i}))\bar{e}_{i}$
and $U_{i}=(\bar{A}/{\rm Rad}(\bar{\phi}_{i}))\bar{e}_{i}$,
where $\bar{e}_{i}=\bar{e}_{i, 11}+\dots +\bar{e}_{i, n_{i}1}$
and 
$$\{\bar{e}_{i, jk}\mid j=1, \dots, n_{i}, k=1, \dots, n_{ij}\}$$
is a complete set of orthogonal primitive idempodents 
of $A/{\rm Rad}(\phi_{i})$ such that 
$$1_{\bar{A}/{\rm Rad}(\bar{\phi}_{i})}=1_{A/{\rm Rad}(\phi_{i})}
=\sum_{j=1}^{n_{i}}
\sum_{k=1}^{n_{ij}}\bar{e}_{i, jk},$$
$\bar{e}_{i,jl}(\bar{A}/{\rm Rad}(\bar{\phi}_{i}))\simeq
\bar{e}_{i,jp}(\bar{A}/{\rm Rad}(\bar{\phi}_{i}))$
and $\bar{e}_{i,jl}(\bar{A}/{\rm Rad}(\bar{\phi}_{i}))
\simeq \bar{e}_{i,kp}(\bar{A}/{\rm Rad}(\bar{\phi}_{i}))$ for $l\ne p$.
Then for $i=1,\dots, n$, $\bar{P}_i$ are basic symmetric 
algebras with symmetric linear functions given by $\bar{\phi}_{i}$, 
$U_i$ are $\bar{A}$-$\bar{P}_i$-bimodules, 
finitely generated and projective as right $\bar{P}_i$-modules.
For each $i$, let $f_{i}\in \hom(M\otimes U_{i}, U_{i})$ be defined by
$f_{i}(m\otimes u_{i})=(0, m)u_{i}$ for $m\in M$ and $u_{i}\in U_{i}$. 
Then $f_i\in \hom_{A,P}(M\otimes_A U_i, U_i)$ and 
for any $m\in M$,
$$\phi(m) = \sum_{i=1}^n (\phi_{i})_{U_i}^{f_i}(m).$$
Moreover, if $\nu$ is an element in $A$ such that $\phi(\nu m) = 0$
for all $m\in M$, then the modules $U_i$ can be chosen in such a way that
$\nu U_i = 0$ for $i=1,\ldots, n$.
\end{thm}
\pf The proof of this result is essentially in \cite{F1}.
But since we need to give 
$\bar{P}_{i}$, $\bar{\phi}_{i}$, $U_{i}$ and $f_{i}$ explicitly, we give a 
complete proof here. 

We  apply Theorem \ref{ar-sym}
to the algebra $\bar{A}$. Then
$\bar{P}_{i}$ for $i=1, \dots, n$ are 
basic symmetric algebras with symmetric linear functions given
by $\bar{\phi}_{i}$ (still denoted by $\bar{\phi}_{i}$),
and $U_{i}$ for $i=1, \dots, n$ are $\bar{A}$-$\bar{P}_{i}$-bimodules,
finitely generated and projective as right $\bar{P}_{i}$-modules,
such that 
$$\bar{\phi}(\bar{a})=\sum_{n=1}^{n}(\phi_{i})_{U_{i}}(\bar{a})$$
for $\bar{a}\in \bar{A}$, where in each term in the right-hand side,
$\bar{a}\in \bar{A}$ is viewed as an element of 
$\text{End}_{\bar{P}_{i}}U_{i}$
given by the left action of $\bar{a}$ on $U_{i}$.
By definition, we have
$\bar{\phi}(a, m)=\phi(m)$ for $a\in A$ and $m\in M$. 
Then we have 
$$\phi(m)=\bar{\phi}(0, m)=
\sum_{n=1}^{n}(\phi_{i})_{U_{i}}(0, m)$$
for $m\in M$. Since $A$ is a subalgebra of $\bar{A}$, $U_{i}$ 
is also a left $A$-module. By the definition of $f_{i}$, we have
\begin{align*}
f_{i}(ma\otimes u)&=(0, ma)u_{i}
=((0, m)(a, 0))u_{i}
=(0, m)(au_{i})
=f_{i}(m\otimes au_{i}),\\
f_{i}(am\otimes u)&=(0, am)u_{i}
=((a, 0)(0, m))u_{i}
=a((0, m)u_{i})
=af_{i}(m\otimes u),\\
f_{i}(m\otimes u)p_{i}&=((0, m)u_{i})p_{i}
=(0, m)(u_{i}p_{i})
=f_{i}(m\otimes u_{i}p_{i})
\end{align*}
for $a\in A$, $m\in M$, $u_{i}\in  U_{i}$ and $p_{i}\in \bar{P}_{i}$.
So $f_{i}\in \hom_{A, \bar{P}_{i}}(M\otimes_{A} U_{i}, U_{i})$.
Then 
$$(\phi_{i})_{U_{i}}(0, m)=(\phi_{i})_{U_{i}}^{f_{i}}(m).$$
Thus we obtain
$$\phi(m)=
\sum_{n=1}^{n}(\phi_{i})_{U_{i}}^{f_{i}}(m).$$
\epfv

\subsection{Genus-one correlation functions from shifted pseudo-$q$-traces}

In this subsection, we recall the results on genus-one correlation functions 
constructed from shifted pseudo-$q$-traces of products of intertwining 
operators obtained by Fiordalisi
 in \cite{F1} and \cite{F}. 

Let $V$ be a vertex operator algebra. In this section, though some of the results hold
for more general vertex operator algebras, we assume that 
$V$ has no nonzero elements of
negative weights  (that is, $V_{(n)}=0$ for $n<0$) and satisfies the 
$C_{2}$-cofiniteness conditions (that is, $\dim V/C_{2}(V)<\infty$, 
where $C_{2}(V)=\langle \res_{x}x^{-2}Y_{V}(u, x)v\mid u, v\in V\rangle$).

Let $P$ be a finite-dimensional associative algebra equipped 
with a symmetric linear function $\phi$. 
A grading-restricted (or lower-bounded)
generalized (or ordinary) $V$-module $W$ equipped with a right $P$-module structure 
such that the the vertex operators on $W$ commute with the right actions of elements of $P$
is called a grading-restricted (or lower-bounded) generalized (or ordinary) $V$-$P$-bimodule. 

In this section, we shall consider mostly grading-restricted generalized $V$-modules and 
grading-restricted generalized $V$-$P$-bimodules. 
We shall also consider intertwining operators without logarithm 
and logarithmic intertwining operators. For simplicity, starting from now on, we shall 
call all these simply intertwining operators, no matter 
whether they contain or do not 
contain the logarithm of the variable.  

Let $W_{1}$, $W_{2}$ and $W_{3}$ be grading-restricted generalized $V$-$P$-bimodules.
For an intertwining operator $\Y$ of type $\binom{W_{3}}{W_{1}W_{2}}$
and $w_{1}\in W_{1}$, we say that $\Y(w_{1}, x)$ is compatible 
with $P$ or $\Y(w_{1}, x)$ is $P$-compatible if 
$$(\Y(w_{1}, x)w_{2})p=\Y(w_{1}, x)(w_{2}p)$$
for $w_{2}\in W_{2}$ and $p\in P$. 
An intertwining operator of type $\binom{W_{3}}{W_{1}W_{2}}$ compatible with $P$ or 
a $P$-intertwining operator of type $\binom{W_{3}}{W_{1}W_{2}}$
is an intertwining operator $\Y$ of type $\binom{W_{3}}{W_{1}W_{2}}$ such that 
$\Y(w_{1}, x)$ is compatible with $P$ for every $w_{1}\in W_{1}$. More generally, let 
$W_{1}, \dots, W_{n}$, $\widetilde{W}_{1}, \dots, \widetilde{W}_{n-1}$
be grading-restricted generalized $V$-modules and let 
$\widetilde{W}_{0}$ and $\widetilde{W}_{n}$ be grading-restricted generalized $V$-$P$-bimodules. 
Let $\Y_{1}, \dots, \Y_{n}$ be intertwining operators of type 
$\binom{\widetilde{W}_{0}}{W_{1}\widetilde{W}_{1}}, \dots, 
\binom{\widetilde{W}_{n-1}}{W_{n}\widetilde{W}_{n}}$, respectively. 
For $w_{1}\in W_{1}, \dots, w_{n}\in W_{n}$, we say that the product
$\Y_{1}(w_{1}, x_{1})\cdots \Y_{n}(w_{n}, x_{n})$ is compatible with 
$P$ or is $P$-compatible if 
$$(\Y_{1}(w_{1}, x_{1})\cdots \Y_{n}(w_{n}, x_{n})\tilde{w}_{n})p
=\Y_{1}(w_{1}, x_{1})\cdots \Y_{n}(w_{n}, x_{n})(\tilde{w}_{n}p)$$
for $\tilde{w}_{n}\in \widetilde{W}_{n}$
and $p\in P$.
We say that the product of $\Y_{1}, \dots, \Y_{n}$ is compatible with $P$ if 
$\Y_{1}(w_{1}, x_{1})\cdots \Y_{n}(w_{n}, x_{n})$
is compatible with $P$ for all $w_{1}\in W_{1}, \dots, w_{n}\in W_{n}$.
If for $w_{1}\in W_{1}, \dots, w_{n}\in W_{n}$,
$\Y_{1}(w_{1}, x_{1})\cdots \Y_{n}(w_{n}, x_{n})$ is compatible with $P$,
then the coefficients of 
$\Y_{1}(w_{1}, x_{1})\cdots \Y_{n}(w_{n}, x_{n})$ 
are elements of $\hom_{P}(\widetilde{W}_{n}, \widetilde{W}_{0})$. 

From \cite{H-cofiniteness}, we know that the conditions to apply the results in 
\cite{HLZ6} and \cite{HLZ7} are
satisfied. In particular, the associativity and commutativity of intertwining operators hold.
In \cite{F1} and \cite{F}, Fiordalisi studied these properties in the case of 
grading-restricted generalized $V$-$P$-bimodules and intertwining 
operators compatible with $P$. But  the main results in \cite{F1} and \cite{F} 
also hold in such more general settings
with the same proofs except that the associativity and commutativity of intertwining 
operators compatible with $P$ should be replaced by the versions of the associativity and 
commutativity below. For these associativity and commutativity of intertwining operators,
we give complete proofs. 

For simplicity, we use
$$\langle w_{4}', \Y_{1}(w_{1}, z_{1})\Y_{2}(w_{2}, z_{2})w_{3}\rangle$$
to denote
$$\langle w_{4}', \Y_{1}(w_{1}, x_{1})\Y_{2}(w_{2}, x_{2})w_{3}\rangle
\lbar_{x_{1}^{n}
=e^{n \log z_{1}}, \; x_{2}^{n}
=e^{n \log z_{2}}, \; \log x_{1}=\log z_{1}, \; \log x_{2}=\log z_{2}},$$
where for $z \in \C^{\times}$, $\log z = \log |z|+i \arg z$ for 
$0\le \arg z < 2\pi$. Similarly, we have the notation
$$\langle w_{4}', \Y_{3}(\Y_{4}(w_{1}, z_{1}-z_{2})w_{2}, z_{2})w_{3}\rangle.$$
We shall use these and similar notations throughout the present paper.

\begin{prop}\label{assoc-P}
Let $W_{1}, W_{2}, W_{3}, W_{4}, W_{5}, W_{6}$ be grading-restricted generalized $V$-modules and 
$\Y_{1}$, $\Y_{2}$, $\Y_{3}$ and $\Y_{4}$ intertwining operators
of types $\binom{W_{4}}{W_{1}W_{5}}$, $\binom{W_{5}}{W_{2}W_{3}}$,
$\binom{W_{4}}{W_{6}W_{3}}$
and $\binom{W_{6}}{W_{2}W_{3}}$, respectively, such that 
$$\langle w_{4}', \Y_{1}(w_{1}, z_{1})\Y_{2}(w_{2}, z_{2})w_{3}\rangle
=\langle w_{4}', \Y_{3}(\Y_{4}(w_{1}, z_{1}-z_{2})w_{2}, z_{2})w_{3}\rangle$$
for $w_{1}\in W_{1}$, $w_{2}\in W_{2}$, $w_{3}\in W_{3}$ and $w_{4}'\in W_{4}'$
in the region $|z_{1}|>|z_{2}|>|z_{1}-z_{2}|>0$. 
Assume that $W_{3}$ and $W_{4}$ are grading-restricted generalized $V$-$P$-bimodules.
Then $\Y_{1}(w_{1}, x_{1})\Y_{2}(w_{2}, x_{2})$ is compatible with $P$ 
if and only if the coefficients of $\Y_{3}(\Y_{4}(w_{1}, x_{0})w_{2}, x_{2})$
as a formal series in powers of $x_{0}$ and nonnegative powers of $\log x_{0}$
 is compatible with $P$. In particular, in the case that
the product of $\Y_{1}$ and $Y_{2}$ 
is compatible with $P$ and
$W_{6}$ is spanned by  the coefficients
of $\Y_{4}(w_{1}, x)w_{2}$ for $w_{1}\in W_{1}$ and $w_{2}\in W_{2}$,
the product of $\Y_{1}$ and $\Y_{2}$ 
is compatible with $P$ if and only if $\Y_{3}$ is compatible with $P$. 
\end{prop}
\pf If $\Y_{1}(w_{1}, x_{1})\Y_{2}(w_{2}, x_{2})$ is compatible with $P$, then in the region 
$|z_{1}|>|z_{2}|>|z_{1}-z_{2}|>0$, we have 
\begin{align}\label{assoc-P-1}
&\langle w_{4}', (\Y_{3}(\Y_{4}(w_{1}, z_{1}-z_{2})w_{2}, z_{2})w_{3})p\rangle\nn
&\quad =\langle p w_{4}', \Y_{3}(\Y_{4}(w_{1}, z_{1}-z_{2})w_{2}, z_{2})w_{3}\rangle\nn
&\quad =\langle p w_{4}', \Y_{1}(w_{1}, z_{1})\Y_{2}(w_{2}, z_{2})w_{3}\rangle\nn
&\quad =\langle w_{4}', (\Y_{1}(w_{1}, z_{1})\Y_{2}(w_{2}, z_{2})w_{3})p\rangle\nn
&\quad =\langle w_{4}', \Y_{1}(w_{1}, z_{1})\Y_{2}(w_{2}, z_{2})(w_{3}p)\rangle\nn
&\quad =\langle w_{4}', \Y_{3}(\Y_{4}(w_{1}, z_{1}-z_{2})w_{2}, z_{2})(w_{3}p)\rangle
\end{align}
for $w_{1}\in W_{1}$, $w_{2}\in W_{2}$, $w_{3}\in W_{3}$, $w_{4}'\in W_{4}'$
and $p\in P$. 
Since both sides of (\ref{assoc-P-1}) are convergent in the region $|z_{2}|>|z_{1}-z_{2}|>0$,
the left-hand and right-hand sides of (\ref{assoc-P-1}) are equal in this larger region. 
In this region, we can take the coefficients of 
$\Y_{4}(w_{1}, z_{1}-z_{2})w_{2}$ in both sides of (\ref{assoc-P-1}). 
Thus the coefficients of $\Y_{3}(\Y_{4}(w_{1}, x_{0})w_{2}, x_{2})$
as a formal series in powers of $x_{0}$ and nonnegative powers of $\log x_{0}$
 is compatible with $P$.
If the product of $\Y_{1}$ and $\Y_{2}$ is compatible with $P$ and 
$W_{6}$  is spanned by  the coefficients
of $\Y_{4}(w_{1}, x_{0})w_{2}$ for $w_{1}\in W_{1}$ and $w_{2}\in W_{2}$,
then we obtain 
$$\langle w_{4}', (\Y_{3}(w_{6}, z_{2})w_{3})p\rangle
=\langle w_{4}', \Y_{3}(w_{6}, z_{2})(w_{3}p)\rangle$$
for $w_{3}\in W_{3}$, $w_{6}\in W_{6}$, $w_{4}'\in W_{4}'$ and $p\in P$ in the region $z_{2}\ne 0$. 
This shows that $\Y_{3}$ is compatible with $P$. 

Conversely, if the coefficients of $\Y_{3}(\Y_{4}(w_{1}, x_{0})w_{2}, x_{2})$ 
as a series in powers of $x_{0}$ and nonnegative powers of $\log x_{0}$
are compatible with $P$, then in the region 
$|z_{1}|>|z_{2}|>|z_{1}-z_{2}|>0$, we have 
\begin{align}\label{assoc-P-2}
&\langle w_{4}', (\Y_{1}(w_{1}, z_{1})\Y_{2}(w_{2}, z_{2})w_{3})p\rangle\nn
&\quad =\langle p w_{4}', \Y_{1}(w_{1}, z_{1})\Y_{2}(w_{2}, z_{2})w_{3}\rangle\nn
&\quad =\langle p w_{4}', \Y_{3}(\Y_{4}(w_{1}, z_{1}-z_{2})w_{2}, z_{2})w_{3}\rangle\nn
&\quad =\langle w_{4}', (\Y_{3}(\Y_{4}(w_{1}, z_{1}-z_{2})w_{2}, z_{2})w_{3})p\rangle\nn
&\quad =\langle w_{4}', \Y_{1}(w_{1}, z_{1})\Y_{2}(w_{2}, z_{2})(w_{3}p)\rangle\nn
&\quad =\langle w_{4}', \Y_{1}(w_{1}, z_{1})\Y_{2}(w_{2}, z_{2})(w_{3}p)\rangle
\end{align}
for $w_{1}\in W_{1}$, $w_{2}\in W_{2}$, $w_{3}\in W_{3}$, $w_{4}'\in W_{4}'$
and $p\in P$. 
Since both sides of (\ref{assoc-P-2}) are convergent in the region $|z_{1}|>|z_{2}|>0$,
the left-hand and right-hand sides of (\ref{assoc-P-1}) are equal in this larger region. 
This shows that $\Y_{1}(w_{1}, x_{1})\Y_{2}(w_{2}, x_{2})$ is compatible with $P$. 
If $\Y_{3}$ is compatible with $P$, then the coefficients of $\Y_{3}(\Y_{4}(w_{1}, x_{0})w_{2}, x_{2})$ 
as a series in powers of $x_{0}$ and nonnegative powers of $\log x_{0}$
are compatible with $P$. Thus $\Y_{1}(w_{1}, x_{1})\Y_{2}(w_{2}, x_{2})$ 
is compatible with $P$ for all $w_{1}\in W_{1}$ and $w_{2}\in W_{2}$, that is, 
the product of $\Y_{1}$ and $\Y_{2}$ is compatible with $P$. 
\epfv

For two analytic functions $f$ and $g$ on two regions, we shall use $f \sim g$ 
to mean that $f$ and $g$ are analytic extensions of each other. 

\begin{prop}\label{commu-P}
Let $W_{1}, W_{2}, W_{3}, W_{4}, W_{5}, W_{6}$ be grading-restricted generalized $V$-modules and 
$\Y_{1}$, $\Y_{2}$, $\Y_{5}$ and $\Y_{6}$ intertwining operators
of types $\binom{W_{4}}{W_{1}W_{5}}$, $\binom{W_{5}}{W_{2}W_{3}}$,
$\binom{W_{4}}{W_{2}W_{6}}$
and $\binom{W_{6}}{W_{1}W_{3}}$, respectively, such that 
$$\langle w_{4}', \Y_{1}(w_{1}, z_{1})\Y_{2}(w_{2}, z_{2})w_{3}\rangle\sim 
\langle w_{4}', \Y_{5}(w_{2}, z_{2})\Y_{6}(w_{1}, z_{1})w_{3}\rangle$$
for $w_{1}\in W_{1}$, $w_{2}\in W_{2}$, $w_{3}\in W_{3}$ and $w_{4}'\in W_{4}'$.
Assume that $W_{3}$ and $W_{4}$ are grading-restricted generalized $V$-$P$-bimodules.
Then $\Y_{1}(w_{1}, z_{1})\Y_{2}(w_{2}, z_{2})$ is compatible with $P$ 
if and only if $\Y_{5}(w_{2}, z_{2})\Y_{6}(w_{1}, z_{1})$ is compatible with $P$.
In particular, the product of $\Y_{1}$ and $\Y_{2}$ is compatible with $P$ if and only if 
the product of $\Y_{5}$ and $\Y_{6}$ is compatible with $P$. 
\end{prop}
\pf
We need only prove the ``if'' part; the ``only if" part is obtained by symmetry. 

If $\Y_{5}(w_{2}, z_{2})\Y_{6}(w_{1}, z_{1})$ is compatible with $P$, then
we have 
\begin{align}\label{commu-P-1}
&\langle w_{4}', (\Y_{1}(w_{1}, z_{1})\Y_{2}(w_{2}, z_{2})w_{3})p\rangle\nn
&\quad =\langle p w_{4}', \Y_{1}(w_{1}, z_{1})\Y_{2}(w_{2}, z_{2})w_{3}\rangle\nn
&\quad \sim \langle p w_{4}', \Y_{5}(w_{2}, z_{2})\Y_{6}(w_{1}, z_{1})w_{3}\rangle\nn
&\quad = \langle w_{4}', (\Y_{5}(w_{2}, z_{2})\Y_{6}(w_{1}, z_{1})w_{3})p\rangle\nn
&\quad =\langle w_{4}', \Y_{5}(w_{1}, z_{1})\Y_{6}(w_{2}, z_{2})(w_{3}p)\rangle\nn
&\quad \sim\langle w_{4}', \Y_{1}(w_{1}, z_{1})\Y_{2}(w_{2}, z_{2})(w_{3}p)\rangle.
\end{align}
Since both sides of (\ref{commu-P-1}) are analytic functions in the same region 
$|z_{1}|>|z_{2}|>0$, $0\le \arg z_{1}, \arg z_{2}<2\pi$, these two sides must be equal, proving that $\Y_{1}(w_{1}, z_{1})\Y_{2}(w_{2}, z_{2})$ is compatible with $P$. 
If the product of $\Y_{5}$ and $\Y_{6}$ is compatible with $P$, then 
$\Y_{1}(w_{1}, z_{1})\Y_{2}(w_{2}, z_{2})$ is compatible with $P$ for all 
$w_{1}\in W_{1}$ and $w_{2}\in W_{2}$. Thus
the product of $\Y_{1}$ and $\Y_{2}$ is compatible with $P$. 
\epfv

Let $W$ be a grading-restricted generalized $V$-$P$-bimodules which is projective 
as a right $P$-module. Then for each $n\in \C$, the homogeneous subspace 
$W_{[n]}$ of $W$ of conformal weight $n$
is a finite-dimensional right projective $P$-module. 
Let $\phi$ be a symmetric linear function on $P$. 
We have a pseudo-trace function $\phi_{W_{[n]}}$ on $\text{End}_{P} W_{[n]}$ (see Subsection 2.1).
We also know that $W$ is of finite length. Let $K$ be the length of $W$. Then 
we know that $L_{W}(0)_{N}^{K+1}$ is $0$. 
For $f\in \hom_{P}(W, \overline{W})$,
define
$$\tr_{W}^{\phi}f q^{L_{W}(0)}
=\sum_{n\in \C}\sum_{k=0}^{K}\phi_{W_{[n]}}\left(\pi_{n}f
\frac{L_{W}(0)_{N}^{k}}{k!}\right)\lbar_{W_{[n]}}(\log q)^{k}q^{n},$$
where for $n\in \C$, $\pi_{n}$ is the projection from $W$ to $W_{[n]}$.

For a lower-bounded generalized $V$-module $W$,
as in \cite{H-mod-inv-int}, let 
$$\mathcal{U}_{W}(x)=(2\pi ix)^{L_{W}(0)}e^{-L^{+}_{W}(A)}\in (\mbox{\rm End}\;W)\{x\}[\log x],$$
where $(2\pi i)^{L_{W}(0)}=e^{(\log 2\pi +i \frac{\pi}{2})L_{W}(0)}$, $x^{L_{W}(0)}
=x^{L_{W}(0)_{S}}e^{(\log x)L_{W}(0)_{N}}$, 
$L^{+}_{W}(A)=\sum_{j\in \mathbb{Z}_{+}}
A_{j}L_{W}(j)$ and $A_{j}$ for $j\in \N$ are given by 
$$\frac{1}{2\pi i}\log(1+2\pi i y)=\left(\exp\left(\sum_{j\in \mathbb{Z}_{+}}
A_{j}y^{j+1}\frac{\partial}{\partial y}\right)\right)y.$$
For $\Y=Y_{W}$,  (1.5) in \cite{H-mod-inv-int} gives
\begin{equation}\label{U-1-v-op}
\mathcal{U}_{W}(1)Y_{W}\left(v, \frac{1}{2\pi i}\log (1+x)\right)
=Y_{W}((1+x)^{L_{V}(0)}\mathcal{U}_{V}(1)v, x)\mathcal{U}_{W}(1)
\end{equation}
for $v\in V$. Let $\Y$ be an intertwining operator 
of type $\binom{W_{3}}{WW_{2}}$, where $W_{2}$ and $W_{3}$ 
are also lower-bounded generalized $V$-modules.
For $z\in \C$, let $q_{z}=e^{2\pi i z}$. Then as in \cite{H-mod-inv-int},
we call $\mathcal{Y}(\mathcal{U}_{W_{1}}(q_{z})w_{1}, q_{z})$
a geometrically-modified intertwining operator. 

Let 
$W_{1}, \dots, W_{n}$, $\widetilde{W}_{1}, \dots, \widetilde{W}_{n-1}$
be grading-restricted generalized $V$-modules and let 
$\widetilde{W}_{0}=\widetilde{W}_{n}$ be a grading-restricted generalized $V$-$P$-bimodule
which is projective as a right $P$-module. 
Let $\Y_{1}, \dots, \Y_{n}$ be intertwining operators of types
$\binom{\widetilde{W}_{0}}{W_{1}\widetilde{W}_{1}}, \dots, 
\binom{\widetilde{W}_{n-1}}{W_{n}\widetilde{W}_{n}}$, respectively. 
We assume that
\begin{equation}\label{prod-int-comp-P}
\mathcal{Y}_{1}(\mathcal{U}_{W_{1}}(q_{z_{1}})w_{1}, x_{1})
\cdots
\mathcal{Y}_{n}(\mathcal{U}_{W_{n}}(q_{z_{n}})w_{n}, x_{n})
\end{equation}
is compatible with $P$. 
If  the product of 
$\Y_{1}, \dots, \Y_{n}$ is compatible with $P$, then this assumption 
is true for any $w_{1}\in W_{1}, \dots, w_{n}\in W_{n}$. 
Since  the product 
$$\mathcal{Y}_{1}(\mathcal{U}_{W_{1}}(q_{z_{1}})w_{1}, q_{z_{1}})
\cdots
\mathcal{Y}_{n}(\mathcal{U}_{W_{n}}(q_{z_{n}})w_{n}, q_{z_{n}})$$
of geometrically-modified
intertwining operators 
is absolutely convergent to an element of $\hom(W, \overline{W})$, 
it is in fact absolutely convergent to an element of $\hom_{P}(W, \overline{W})$.
Then we have the pseudo-$q_{\tau}$-trace shifted by $-\frac{c}{24}$
or simply the shifted pseudo-$q_{\tau}$-trace
\begin{equation}\label{pseudo-tr-prod-int}
\tr_{\tilde{W}_{n}}^{\phi}\mathcal{Y}_{1}(\mathcal{U}_{W_{1}}(q_{z_{1}})w_{1}, q_{z_{1}})
\cdots
\mathcal{Y}_{n}(\mathcal{U}_{W_{n}}(q_{z_{n}})w_{n}, q_{z_{n}})
q_{\tau}^{L(0)-\frac{c}{24}}
\end{equation}
of products of $n$ geometrically-modified intertwining operators. 

We now state several results of Fiordalisi in \cite{F1} and \cite{F} 
generalizing the corresponding results in \cite{H-mod-inv-int} in the semisimple case.
As we discussed above, our statements of these results are slightly more 
general than those in \cite{F1} and \cite{F} but the proofs there in fact 
already gave these results. 

\begin{thm}[Convergence and analytic extension \cite{F1} \cite{F}]
\label{F-genus1-conv}
For $w_{1}\in W_{1}, \dots$, $w_{n}\\ \in W_{n}$ such that 
(\ref{prod-int-comp-P}) is compatible with $P$,
the series (\ref{pseudo-tr-prod-int}) 
is absolutely convergent in the region 
$1>|q_{z_{1}}|>\cdots >|q_{z_{n}}|>|q_{\tau}|>0$
and can be analytically extended to a multivalued analytic function in the region 
$\Im(\tau)>0$, $z_{i}\ne z_{j}+k\tau+l$ for $i\ne j$, $k, l\in \Z$. 
\end{thm}

For $w_{1}\in W_{1}, \dots, w_{n}\in W_{n}$, we denote 
the multivalued analytic function in Theorem \ref{F-genus1-conv} by
\begin{equation}\label{analy-ext-pseudo-tr}
\overline{F}^{\phi}_{\mathcal{Y}_{1}, \dots, \mathcal{Y}_{n}}(w_{1}, \dots, w_{n};
z_{1}, \dots, z_{n}; \tau).
\end{equation}
Note that the multivalued analytic function (\ref{analy-ext-pseudo-tr})
has a particular branch (usually called a preferred branch by the author)
in the region $|q_{z_{1}}|>\cdots >|q_{z_{n}}|>|q_{\tau}|>0$
given by (\ref{pseudo-tr-prod-int}). 
Such a function of $z_{1}, \dots, z_{n}$ and $\tau$ 
is called a genus-one $n$-point 
correlation function from a shifted pseudo-$q_{\tau}$-trace.

\begin{thm}[Genus-one commutativity \cite{F1} \cite{F}]\label{F-genus1-comm}
For $w_{1}\in W_{1}, \dots, w_{n}\in W_{n}$ such that 
(\ref{prod-int-comp-P}) is compatible with $P$ and for  $1\le k\le n-1$,
there exist grading-restricted generalized 
$V$-modules $\widehat{W}_{k}$ and intertwining operators
$\widehat{\Y}_{k}$ and $\widehat{\Y}_{k+1}$
of types $\binom{\widehat{W}_{k}}{W_{k}\widetilde{W}_{k+1}}$
and $\binom{\widetilde{W}_{k-1}}{W_{k+1}\widehat{W}_{k}}$, respectively,
such that 
\begin{align*}
&\mathcal{Y}_{1}(\mathcal{U}_{W_{1}}(q_{z_{1}})w_{1}, x_{1})
\cdots \mathcal{Y}_{k-1}(\mathcal{U}_{W_{k-1}}(q_{z_{k-1}})w_{k-1}, x_{k-1})
\widehat{\Y}_{k+1}(\mathcal{U}_{W_{k+1}}(q_{z_{k+1}})w_{k+1}, x_{k+1})\cdot\nn
&\quad\quad\cdot
\widehat{\Y}_{k}(\mathcal{U}_{W_{k}}(q_{z_{k}})w_{k}, x_{k})
\mathcal{Y}_{k+2}(\mathcal{U}_{W_{k+2}}(q_{z_{k+2}})w_{k+2}, x_{k+2})\cdots
\mathcal{Y}_{n}(\mathcal{U}_{W_{n}}(q_{z_{n}})w_{n}, x_{n})
\end{align*}
is compatible with $P$ and
\begin{align*}
&\overline{F}^{\phi}_{\mathcal{Y}_{1}, \dots, \mathcal{Y}_{n}}(w_{1}, 
\dots, w_{n};
z_{1}, \dots, z_{n}; \tau)\nn
&\quad =\overline{F}^{\phi}_{\mathcal{Y}_{1}, \dots, \mathcal{Y}_{k-1},
\widehat{\Y}_{k+1}, \widehat{\Y}_{k}, \mathcal{Y}_{k+2}
\dots, \mathcal{Y}_{n}}(w_{1}, 
\dots, w_{k-1}, w_{k+1}, w_{k}, w_{k+2}, \dots, w_{n};\nn
&\quad\quad\quad\quad\quad\quad\quad\quad\quad\quad\quad\quad\quad
\quad\quad\quad
z_{1}, \dots, z_{k-1}, z_{k+1}, z_{k}, z_{k+2}, 
\dots, z_{n}; \tau).
\end{align*}
More generally, for any $\sigma\in S_{n}$, 
there exist grading-restricted generalized  $V$-modules $\widehat{W}_{i}$ for $i=1, \dots, n-1$ 
and intertwining operators $\widehat{\Y}_{i}$ of 
types $\binom{\widehat{W}_{i-1}}{W_{\sigma(i)}\widehat{W}_{i}}$ for $i=1, \dots, n$
(where $\widehat{W}_{0}=\widehat{W}_{n}=\widetilde{W}_{0}=\widetilde{W}_{n}$), respectively, 
such that 
$$\widehat{\mathcal{Y}}_{1}(\mathcal{U}_{W_{\sigma(1)}}
(q_{z_{\sigma(1)}})w_{\sigma(1)}, x_{\sigma(1)})
\cdots
\widehat{\mathcal{Y}}_{n}(\mathcal{U}_{W_{\sigma(n)}}(q_{z_{\sigma(n)}})w_{\sigma(n)}, x_{\sigma(n)})$$
 is compatible with $P$ and 
\begin{align*}
\overline{F}^{\phi}_{\mathcal{Y}_{1}, \dots, \mathcal{Y}_{n}}(w_{1}, 
\dots, w_{n};
z_{1}, \dots, z_{n}; \tau)\ =\overline{F}^{\phi}_{\widehat{\Y}_{1}, \dots, 
\widehat{\Y}_{n}}(w_{\sigma(1)}, 
\dots, w_{\sigma(n)};
z_{\sigma(1)}, \dots, z_{\sigma(n)}; \tau).
\end{align*}
\end{thm}

\begin{thm}[Genus-one associativity \cite{F1} \cite{F}]\label{F-genus1-assoc}
For $w_{1}\in W_{1}, \dots, w_{n}\in W_{n}$ such that 
(\ref{prod-int-comp-P}) is compatible with $P$ and for $1\le k\le n-1$, 
there exist a grading-restricted generalized $V$-module $\widehat{W}_{k}$ and 
intertwining operators $\widehat{\Y}_{k}$ and 
$\widehat{\Y}_{k+1}$ of 
types $\binom{\widehat{W}_{k}}{W_{k} W_{k+1}}$ and 
$\binom{\widetilde{W}_{k-1}}{\widehat{W}_{k}\widetilde{W}_{k+1}}$, respectively, 
such that the coefficients as a series in powers of $x_{0}$ and nonnegative powers of 
$\log x_{0}$ of
\begin{align*}
&\mathcal{Y}_{1}(\mathcal{U}_{W_{1}}(q_{z_{1}})w_{1}, x_{1})
\cdots \mathcal{Y}_{k-1}(\mathcal{U}_{W_{k-1}}(q_{z_{k-1}})w_{k-1}, x_{k-1})
\widehat{\Y}_{k+1}(\mathcal{U}_{W_{k+1}}(q_{z_{k+1}})
\widehat{\Y}_{k}(w_{k}, x_{0})w_{k+1}, x_{k+1})\cdot\nn
&\quad\quad\cdot
\mathcal{Y}_{k+2}(\mathcal{U}_{W_{k+2}}(q_{z_{k+2}})w_{k+2}, x_{k+2})\cdots
\mathcal{Y}_{n}(\mathcal{U}_{W_{n}}(q_{z_{n}})w_{n}, x_{n})
\end{align*}
is compatible with $P$ and
\begin{align*}%\label{g-1-iter}
&\overline{F}^{\phi}_{\mathcal{Y}_{1}, \dots, \mathcal{Y}_{k-1},
\widehat{\Y}_{k+1}, \mathcal{Y}_{k+2}, \dots,
\mathcal{Y}_{n}}(w_{1}, 
\dots, w_{k-1}, \widehat{\Y}(w_{k}, z_{k}-z_{k+1})
w_{k+1}, \nn
&\quad\quad\quad\quad\quad\quad\quad\quad\quad\quad\quad\quad\quad\quad
w_{k+2}, \dots, w_{n};
z_{1}, \dots, z_{k-1}, z_{k+1}, \dots, z_{n}; \tau)\nn
&\quad=\sum_{r\in \mathbb{R}}
\overline{F}^{\phi}_{\mathcal{Y}_{1}, \dots, \mathcal{Y}_{k-1},
\widehat{\Y}_{k+1}, \mathcal{Y}_{k+2}, \dots,
\mathcal{Y}_{n}}(w_{1}, 
\dots, w_{k-1}, \pi_{r}(\widehat{\Y}(w_{k}, z_{k}-z_{k+1})
w_{k+1}), \nn
&\quad\quad\quad\quad\quad\quad\quad\quad\quad\quad\quad\quad\quad\quad
w_{k+2}, \dots, w_{n};
z_{1}, \dots, z_{k-1}, z_{k+1}, \dots, z_{n}; \tau)
\end{align*}
is absolutely convergent in the region $1>|q_{z_{1}}|>\cdots 
>|q_{z_{k-1}}|>|q_{z_{k+1}}|>
\dots >|q_{z_{n}}|>|q_{\tau}|>0$ and 
$1>|q_{(z_{k}-z_{k+1})}-
1|>0$
and is convergent to 
$$\overline{F}^{\phi}_{\mathcal{Y}_{1}, \dots, \mathcal{Y}_{n}}(w_{1}, 
\dots, w_{n};
z_{1}, \dots, z_{n}; \tau)$$ 
in the region $1>|q_{ z_{1}}|>\cdots 
>|q_{z_{n}}|>|q_{\tau}|>0$,
$|q_{(z_{k}-z_{k+1})}|>1>|q_{(z_{k}-z_{k+1})}-
1|>0$.
\end{thm}

\renewcommand{\theequation}{\thesection.\arabic{equation}}
\renewcommand{\thethm}{\thesection.\arabic{thm}}
\setcounter{equation}{0} \setcounter{thm}{0} 

\section{Associative algebras, lower-bounded generalized $V$-modules
and intertwining operators}

In \cite{H-mod-inv-int}, for a vertex operator algebra $V$, 
an associative algebra 
$\widetilde{A}(V)$ isomorphic to the Zhu algebra $A(V)$, $\widetilde{A}(V)$-modules
and $\widetilde{A}(V)$-bimodules are introduced and used
in the proof in the same paper 
of the modular invariance conjecture of Moore and Seiberg for 
rational conformal field theories. In \cite{H-aa-va} and 
\cite{H-aa-int-op}, the associative algebras 
$A^{\infty}(V)$ and 
$A^{N}(V)$ for $N\in \N$, their graded modules and 
their bimodules associated to a lower-bounded generalized $V$-module $W$ 
are introduced and studied. 
In this section, we first prove more results on 
the $A^{\infty}(V)$-bimodule $A^{\infty}(W)$ and 
the $A^{N}(V)$-bimodules $A^{N}(W)$ for $N\in \N$, which will 
be needed in Section 4.
We then introduce associative algebras $\widetilde{A}^{\infty}(V)$ and 
$\widetilde{A}^{N}(V)$ for $N\in \N$ isomorphic to the associative algebras 
$A^{\infty}(V)$ and 
$A^{N}(V)$ for $N\in \N$. 
As in \cite{H-mod-inv-int}, these algebras and their modules can be 
obtained by using some operators corresponding 
to a canonical conformal transformation 
from an annulus to a parallelogram
on $V$ and on 
lower-bounded generalized $V$-modules, respectively. 
We then transport the results obtained in \cite{H-aa-va}
and \cite{H-aa-int-op} using these operators to results 
on $\widetilde{A}^{\infty}(V)$,
$\widetilde{A}^{N}(V)$ for $N\in \N$, their modules and their bimodules.
We refer the reader to \cite{H-aa-va} and \cite{H-aa-int-op} 
for the basic material and notations on these associative algebras,
their modules and bimodules. 

In this section, we in general do not assume that $V$ is $C_{2}$-cofinite.
But we will prove some results needed in later sections when $V$ is $C_{2}$-cofinite. 

Let $W=\coprod_{m\in \C}W_{[m]}$
be a lower-bounded generalized $V$-module, where $W_{[m]}$ for 
$m\in \C$ are generalized eigenspaces for $L_{W}(0)$ with eigenvalues $m$. 
Let 
$$\Gamma(W) = \{\mu\in\C/\Z \mid \text{there exist nonzero elements 
of $W$ of weights in $\mu$}\},$$
$h^{\mu}\in  \mu$ such that $W_{[h^{\mu}]}\ne  0$ but 
$W_{[h^{\mu}-n]} = 0$ for $n \in \Z_{+}$, and
$$W_{\l n\r} =\coprod_{\mu\in \Gamma(W)}
W_{[h^{\mu}+n]}.$$
Then we have
$$W =\coprod_{\mu\in\Gamma(W)}\coprod_{n\in \N}
W_{[h^{\mu}+n]} =\coprod_{n\in \N}W_{\l n\r}.$$
See \cite{H-aa-int-op} for more details.

We first prove some results on $A^{\infty}(W)$. 
Recall from \cite{H-aa-int-op} that $U^{\infty}(W)$ is
the space of  column-finite infinite
matrices with entries in $W$, but doubly indexed by $\N$ instead of $\Z_{+}$. 
Recall also that for $w\in W$ and $k, l\in \N$, $[w]_{kl}$ is 
the matrix with the $(k, l)$-entry being $w$ and all the other entries 
being $0$. Elements of $U^{\infty}(W)$ are suitable (possibly infinite) 
sums of elements of the form 
$[w]_{kl}$ for $w\in W$, $k, l\in \N$. 

Let $O^{\infty}(W)$ be the subspace of $U^{\infty}(W)$ 
spanned by infinite linear combinations of elements of the form 
$$ \res_{x}
x^{-k-l-p-2}(1+x)^{l}[Y_{W}((1+x)^{L_{V}(0)}v, x)w]_{kl}$$
for $v\in V$, $w\in W$,  $k, l, p\in \N$, 
with each pair $(k, l)$ appearing in the 
linear combinations only finitely many times. 

Let $W_{2}$ and $W_{3}$ be lower-bounded generalized $V$-modules
and $\Y$ an intertwining operator of type $\binom{W_{3}}{WW_{2}}$.
Then as we discussed above, 
\begin{align*}
W_{2} &=\coprod_{\mu\in\Gamma(W_{2})}\coprod_{n\in \N}
(W_{2})_{[h^{\mu}_{2}+n]} =\coprod_{n\in \N}(W_{2})_{\l n\r},\\
W_{3} &=\coprod_{\nu\in\Gamma(W_{3})}\coprod_{n\in \N}
(W_{3})_{[h^{\nu}_{3}+n]} =\coprod_{n\in \N}(W_{3})_{\l n\r}.
\end{align*}
For $w\in W$, let $\Y^{0}(w, x)$ be the constant term in $\Y(w, x)$
when $\Y(w, x)$ is viewed as a power series in $\log x$. 
Then we have a linear map $\vartheta_{\Y}: U^{\infty}(W) \to 
\hom(W_{2}, W_{3})$ defined by 
$$\vartheta_{\Y}([w]_{kl})w_{2}
=\sum_{\nu\in \Gamma(W_{3})}\res_{x}x^{h_{2}^{\mu}-h_{3}^{\nu}+l-k-1}
\Y^{0}(x^{L_{W}(0)_{S}}w, x)w_{2}$$
for $k, l\in \N$, $w\in W$ and $w_{2}\in W_{2}$.
Let $Q^{\infty}(W)$ be the intersection of $\ker \vartheta_{\Y}$ for all 
lower-bounded generalized $V$-modules $W_{2}$ and $W_{3}$ and all intertwining operators $\Y$ of type 
$\binom{W_{3}}{WW_{2}}$. 

\begin{prop}\label{o-infty-q-infty}
We have $O^{\infty}(W)\subset Q^{\infty}(W)$.
\end{prop}
\pf
From the definition of $Q^{\infty}(W)$, we need to prove 
$\vartheta_{\Y}(O^{\infty}(V))=0$ for every pair of
lower-bounded generalized $V$-modules $W_{2}$, $W_{3}$ and 
every intertwining operator $\Y$ of type $\binom{W_{3}}{W W_{2}}$.
For 
$$\res_{x_{0}}
x_{0}^{-k-l-p-2}(1+x_{0})^{l}[Y_{W}((1+x_{0})^{L(0)}v, x_{0})w]_{kl}
\in O^{\infty}(W),$$ 
where $v\in V$, $w\in W$, $k, l, p\in \N$, 
and for $\mu \in \Gamma(W_{2})$,
$w_{2}\in (W_{2})_{[h_{2}^{\mu}+l]}\subset (W_{2})_{\l l\r}$,  we have
\begin{align}\label{O-infty-mod-1}
\vartheta_{\Y}&(\res_{x_{0}}
x_{0}^{-k-l-p-2}(1+x_{0})^{l}[Y_{V}((1+x_{0})^{L(0)}v, x_{0})w]_{kl})
[w_{2}]_{l}\nn
&=\sum_{\nu\in \Gamma(W_{3})}
\text{Coeff}_{\log x_{2}}^{0}
\res_{x_{2}}x_{2}^{h_{2}^{\mu}-h_{3}^{\nu}+l-k-1}\res_{x_{0}}
x_{0}^{-k-l-p-2}(1+x_{0})^{l}\cdot\nn
&\quad\quad\quad\quad\quad\cdot [\Y(x_{2}^{L_{V}(0)}
Y_{W}((1+x_{0})^{L(0)}v, x_{0})w, x_{2})w_{2}]_{k}\nn
&=\sum_{\nu\in \Gamma(W_{3})}
\text{Coeff}_{\log x_{2}}^{0}
\res_{x_{2}}x_{2}^{h_{2}^{\mu}-h_{3}^{\nu}+l-k-1}\res_{x_{0}}
x_{0}^{-k-l-p-2}(1+x_{0})^{l}\cdot\nn
&\quad\quad\quad\quad\quad\cdot [\Y(
Y_{W}(x_{2}^{L_{V}(0)}(1+x_{0})^{L(0)}v, x_{0}x_{2})
x_{2}^{L_{V}(0)}w, x_{2})
w_{2}]_{k}\nn
&=\sum_{\nu\in \Gamma(W_{3})}
\res_{x_{0}}\text{Coeff}_{\log x_{2}}^{0}\res_{x_{2}}
x_{0}^{-k-l-p-2}x_{2}^{h_{2}^{\mu}-h_{3}^{\nu}-k-1}
\res_{x_{1}}x_{1}^{l}x_{1}^{-1}
\delta\left(\frac{x_{2}+x_{0}x_{2}}{x_{1}}\right)\cdot\nn
&\quad\quad\quad\quad\quad\cdot 
[\Y(Y_{W}(x_{1}^{L_{V}(0)}v, x_{0}x_{2})
x_{2}^{L_{V}(0)}w, x_{2})w_{2}]_{k}\nn
&=\sum_{\nu\in \Gamma(W_{3})}
\res_{x_{0}}\text{Coeff}_{\log x_{2}}^{0}\res_{x_{2}}
x_{0}^{-k-l-p-2}x_{2}^{h_{2}^{\mu}-h_{3}^{\nu}-k-1}
\res_{x_{1}}x_{1}^{l}x_{0}^{-1}x_{2}^{-1}
\delta\left(\frac{x_{1}-x_{2}}{x_{0}x_{2}}\right)\cdot\nn
&\quad\quad\quad\quad\quad\cdot 
[Y_{W_{3}}(x_{1}^{L_{V}(0)}v, x_{1})
\Y(x_{2}^{L_{W}(0)}w, x_{2})w_{2}]_{k}\nn
&\quad -\sum_{\nu\in \Gamma(W_{3})}
\res_{x_{0}}\text{Coeff}_{\log x_{2}}^{0}\res_{x_{2}}
x_{0}^{-k-l-p-2}x_{2}^{h_{2}^{\mu}-h_{3}^{\nu}-k-1}
\res_{x_{1}}x_{1}^{l}x_{0}^{-1}x_{2}^{-1}
\delta\left(\frac{x_{2}-x_{1}}{-x_{0}x_{2}}\right)\cdot\nn
&\quad\quad\quad\quad\quad\cdot 
[\Y(x_{2}^{L_{W}(0)}w, x_{2})Y_{W_{2}}(x_{1}^{L_{V}(0)}v, 
x_{1})w_{2}]_{k}\nn
&=\sum_{\nu\in \Gamma(W_{3})}\res_{x_{1}}
\text{Coeff}_{\log x_{2}}^{0}\res_{x_{2}}
x_{1}^{-k-p-2}(1-x_{1}^{-1}x_{2})^{-k-l-p-2}
x_{2}^{h_{2}^{\mu}-h_{3}^{\nu}+l+p}
\cdot\nn
&\quad\quad\quad\quad\quad\cdot 
[Y_{W_{3}}(x_{1}^{L_{V}(0)}v, x_{1})
\Y(x_{2}^{L_{W}(0)}w, x_{2})w_{2}]_{k}\nn
&\quad -\sum_{\nu\in \Gamma(W_{3})}\res_{x_{1}}
\text{Coeff}_{\log x_{2}}^{0}\res_{x_{2}}
(-1+x_{1}x_{2}^{-1})^{-k-l-p-2}x_{1}^{l}
x_{2}^{h_{2}^{\mu}-h_{3}^{\nu}-k-2}
\cdot\nn
&\quad\quad\quad\quad\quad\cdot 
[\Y(x_{2}^{L_{W}(0)}w, x_{2})
Y_{W_{2}}(x_{1}^{L_{V}(0)}v, x_{1})w_{2}]_{k}.
\end{align}
Since $w_{2}\in (W_{2})_{[h_{2}^{\mu}+l]}$ and 
the series $(1-x_{1}^{-1}x_{2})^{-k-l-p-2}$ 
contains only nonnegative powers of $x_{2}$,
$$\res_{x_{2}}(1-x_{1}^{-1}x_{2})^{-k-l-p-2}
x_{2}^{h_{2}^{\mu}-h_{3}^{\nu}+l+p}
\Y(x_{2}^{L_{V}(0)}w, x_{2})w_{2}=0.$$
So the first term in the right-hand side of (\ref{O-infty-mod-1}) is $0$. 
Since $w_{2}\in (W_{2})_{[h_{2}^{\mu}+l]}$ and the series 
$(-1+x_{1}x_{2}^{-1})^{-k-l-p-2}$ contains only nonnegative 
powers of $x_{1}$, 
$$\res_{x_{1}}(-1+x_{1}x_{2}^{-1})^{-k-l-p-2}x_{1}^{l}
Y_{W_{2}}(x_{1}^{L_{V}(0)}v, x_{1})
w_{2}=0.$$
So the second term in the right-hand side of (\ref{O-infty-mod-1}) is also $0$. 
Thus we have $\vartheta_{\Y}(O^{\infty}(W))=0$. 
\epfv

Recall
$$U^{NN}(V)=\left\{\sum_{k=0}^{N}[v]_{kk}\;\lbar
\;v\in V\right\}
\subset U^{\infty}(V)$$
and the subalgebra
$$A^{NN}(V)=\left\{\sum_{k=0}^{N}[v]_{kk}+Q^{\infty}(V)\;\lbar
\;v\in V\right\}$$
of $A^{\infty}(V)$ in Subsection 4.1 \cite{H-aa-va}.
Let 
$$U^{NN}(W)=\left\{\sum_{k=0}^{N}[w]_{kk}\;\lbar
\;w\in W\right\}
\subset U^{\infty}(W)$$
and
$$A^{NN}(W)=\left\{\sum_{k=0}^{N}[w]_{kk}+Q^{\infty}(V)\;\lbar
\;w\in V\right\}
\subset A^{\infty}(W).$$
Also let
$Q^{NN}(W)=Q^{\infty}(W)\cap U^{NN}(W)$ and
$O^{NN}(W)=O^{\infty}(W)
\cap U^{NN}(W)$. Then $A^{NN}(W)$ is linearly isomorphic to 
$U^{NN}(W)/Q^{NN}(W)$. 

In the proof of the next result, 
we shall use the results 
on the $A_{N}(V)$-bimodule $A_{N}(W)$ introduced in \cite{HY}. 
Here we briefly recall some basic definitions and results.
For  $u \in V$ and $w \in
W$, we define
\begin{align*}
u *_N w &= \sum_{m = 0}^N (-1)^m {m+N\choose N}
\res_x x^{-N - m - 1}Y_{W}((1+x)^{L(0)+N}u, x)w,\\
w *_N u &= \sum_{m = 0}^N (-1)^m {m+N\choose N}
\res_x x^{-N - m - 1}\cdot\nn
&\quad\quad\quad\quad \cdot (1+x)^{-(L_{W}(-1)+L_{W}(0))}
Y_{WV}^{W}((1+x)^{L(0)+N}w, x)u.
\end{align*}
Let $O_N(W)$ be the subspace of $W$ spanned by elements of the form
$$\res_x x^{-2N -2}Y_{W}((1+x)^{L(0)+N}u, x)w$$
for
$u \in V$ and $w \in W$. Let $A_N(W) = W/O_N(W)$.
Then it is proved in \cite{HY} that $A_N(W)$ is an $A_{N}(V)$-bimodule
with the left and right actions induced from $*_{N}$ above. 
For a lower-bounded generalized 
$V$-module 
$$W=\coprod_{\mu\gamma(W)}\coprod_{n\in \N}
W_{[h^{\mu}+n]}=\coprod_{n\in \N}W_{\l n\r},$$
let $\Omega_{N}^{0}(W)=\coprod_{n=0}^{N}W_{\l n\r}$
and $G_{N}(W)=W_{\l N\r}$. Then $\Omega_{N}^{0}(W)$ 
is a left $A_{N}(V)$-module and 
$G_{N}(W)$ is a left $A_{N}(V)$-submodule of $\Omega_{N}^{0}(W)$.
Theorem 6.1 in \cite{HY} give a construction of a lower-bounded 
generalized $V$-module $S_{N}(G_{N}(W))$
from the left $A_{N}(V)$-module 
$G_{N}(W)$ such that 
 $G_{N}(S_{N}(G_{N}(W)))$ is equivalent to $G_{N}(W)$
 as $A_{N}(V)$-modules. 
See \cite{HY} 
for details.

\begin{prop}\label{Q-NN=O-NN}
For $N\in \N$, $Q^{NN}(W)=O^{NN}(W)$. 
\end{prop}
\pf
By Proposition \ref{o-infty-q-infty}, $O^{NN}(W)\subset Q^{NN}(W)$.
So we need only prove $Q^{NN}(W)\subset O^{NN}(W)$. 

For a lower-bounded generalized 
$V$-module $W$, by Theorem 6.1 in \cite{HY}, we have a 
lower-bounded generalized $V$-module 
$S_{N}(G_{N}(W))$ such that $G_{N}(S_{N}(G_{N}(W)))$ is equivalent to $G_{N}(W)$
as modules for $A_{N}(V)$. In fact, in 
the construction in \cite{HY}, $G_{N}(W)$ can be an arbitrary 
$A_{N}(V)$-module $M$ 
and we obtain a lower-bounded generalized $V$-module 
$S_{N}(M)$ such that 
$G_{N}(S_{N}(M))$ is equivalent to $M$ as 
$A_{N}(V)$-modules and satisfies the 
following universal property: For any lower-bounded generalized
$V$-module  $W$ and any $A_{N}(V)$-module map $\phi : M\to  G_{N}(W)$, 
there is a unique $V$-module map  $\bar{\phi}: S_{N} (M) \to W$ 
such that  $\bar{\phi}|_{G_{N}(W)} = \phi$. Note that $S_{N}(M)$
 can also be constructed 
using the method in Section 5 of \cite{H-const-twisted-mod}. 

We view $A_{N}(V)$ as a left $A_{N}(V)$-module. Then we obtain a
lower-bounded generalized $V$-module $S_{N}(A_{N}(V))$. 
Let $W_{2}=S_{N}(A_{N}(V))$.
From the construction,  we have $G_{N}(W_{2})=A_{N}(V)$. 

We now construct a lower-bounded generalized $V$-module $W_{3}$
such that $G_{N}(W_{3})=A_{N}(W)\otimes_{A_{N}(V)}A_{N}(V)$
and an intertwining operator $\Y$ of type $\binom{W_{3}}{WW_{2}}$. 
We have 
$$W=\sum_{\mu\in \Gamma(W)}
\coprod_{n\in \N}W_{[h^{\mu}+n]}$$ 
and 
$$W_{2}=\sum_{\nu\in \Gamma(W_{2})}
\coprod_{n\in \N}(W_{2})_{[h_{2}^{\nu}+n]}.$$ 
Fix $h_{3}\in \C$. 
Consider the space 
$$M=\coprod_{\mu\in \Gamma(W), \nu\in \Gamma(W_{2})}
W\otimes t^{-(h_{3}-h^{\mu}-h_{2}^{\nu})}\C[t, t^{-1}]\otimes W_{2}.$$
We shall write an element of $M$ of the form
$w\otimes t^{n}\otimes w_{2}$
as $w(n, 0)w_{2}$ for $w\in W$, $n\in 
-(h_{3}-h^{\mu}-h_{2}^{\nu})+\Z$ and $w_{2}\in W_{2}$. 
We define the weight of the element $w(n, 0)w_{2}$ to be 
$\wt w-n-1+\wt w_{2}$ when $w$ and $w_{2}$ are homogeneous. 
We define $L_{M}(0)_{S}$ be the operator on 
$M$ given by this weight grading. 
We also define an operator $L_{M}(0)_{N}$ on $M$
by 
$$L_{M}(0)_{N}w(n, 0)w_{2}=(L_{W}(0)_{N}w)(n, 0)w_{2}
+w(n, 0)(L_{W_{2}}(0)_{N}w_{2}).$$
Then we have an operator $L_{M}(0)=L_{M}(0)_{S}+L_{M}(0)_{N}$.

From $M$ with the grading 
defined above and the operator $L_{M}(0)$, $g=1_{V}$ and 
$B=h_{3}$, we obtain a universal generalized lower-bounded $V$-module 
$\widehat{W}_{3}=\widehat{M}_{h_{3}}^{1_{V}}$
using the construction given in Subsection 4.2 of \cite{H-affine-twisted-mod}
based on the construction in Section 5 of \cite{H-const-twisted-mod}. 
Let 
$$\widehat{\Y}^{0}(w, x)w_{2}
=\sum_{n\in -(h_{3}-h^{\mu}-h_{2}^{\nu})+\N}
w(n, 0)w_{2}x^{-n-1}$$
for $w\in \coprod_{n\in \N}W_{[h^{\mu}+n]}$,
$w_{2}\in  \coprod_{n\in \N}W_{[h_{2}^{\nu}+n]}$
and 
$$\widehat{\Y}(w, x)w_{2}=x^{L_{\widehat{W}_{3}}(0)}
\widehat{\Y}^{0}(x^{-L_{W}(0)}w, 1)x^{-L_{W_{2}}(0)}w_{2}$$
for $w\in W$ and $w_{2}\in W_{2}$.
We define $o_{\widehat{\Y}}(w)w_{2}$ to be 
$w(\wt w + \wt w_{2} -
h_3 -N- 1, 0)w_{2}$
for homogeneous $w\in W$ and $w_{2}\in G_{N}(W_{2})=A_{N}(V)$
and extend the definition linearly to general $w\in W$ and $w_{2}\in 
A_{N}(V)$. 
Let $J_{1}$ be the 
generalized $V$-submodule of 
$\widehat{W}_{3}$ generated by the elements of the 
forms
$o_{\widehat{\Y}}(w)w_{2}$ for $w\in O_{N}(W)$ 
and $w_{2}\in A_{N}(V)$ and
$o_{\widehat{\Y}}(w\circ_{N}v)w_{2}
-o_{\widehat{\Y}}(w)((v+O_{N}(V))\circ_{N}w_{2})$
for $w\in W$, $v\in V$ and $w_{2}\in A_{N}(V)$
and the coefficients of the formal series of the following 
form 
\begin{align*}
&\frac{d}{dx}\widehat{\Y}(w, x)w_{2}
-\widehat{\Y}(L_{W}(-1)w, x)w_{2},\\
Y_{\widehat{W}_{3}}(v, x_{1})\widehat{\Y}(w, x_{2})w_{2}-&
\widehat{\Y}(w, x_{2})Y_{W_{2}}(v, x_{1})w_{2}
-\res_{x_{0}}x_{1}^{-1}\delta\left(\frac{x_{2}+x_{0}}{x_{1}}\right)
\widehat{\Y}(Y_{W}(v, x_{0})w, x_{2})w_{2}
\end{align*}
for $v\in V$, $w\in W$ and $w_{2}\in W_{2}$. Then the 
lower-bounded generalized $V$-module 
$\widehat{W}_{3}/J_{1}$ is generated by the coefficients of 
formal series of the form
$\widehat{\Y}(w, x)w_{2}+J_{1}$
for $w\in W$ and $w_{2}\in A_{N}(V)$. Moreover, 
elements of $\widehat{W}_{3}/J_{1}$ 
of the form $o_{\widehat{\Y}}(w)(v+O_{N}(V))$ for $w\in W$ and 
$v\in V$ can be written uniquely as 
$o_{\widehat{\Y}}(w\circ_{N}v)(\one +O_{N}(V))$.
From the definition of $J_{1}$, we see that 
$o_{\widehat{\Y}}(w)(\one+O_{N}(V))$ is not in $J_{1}$ if and only 
if  $w\in W
\setminus O_{N}(W)$, or equivalently, 
$o_{\widehat{\Y}}(w)(\one+O_{N}(V))$ is in $J_{1}$ 
if and only if $w\in O_{N}(W)$. 

We define a linear map from $A_{N}(W)=A_{N}(W)\otimes_{A_{N}(V)}
A_{N}(V)$ to $\widehat{W}_{3}/J_{1}$ by 
$$w+O_{N}(W)\mapsto o_{\widehat{\Y}}(w)(\one+O_{N}(V))+J_{1}.$$
Then by the discussion above, we see that 
this map is injective. In particular, we can identify 
the subspace $D$ of $\widehat{W}_{3}/J_{1}$ consisting of elements 
of the form  $o_{\widehat{\Y}}(w)(\one+O_{N}(V))+J_{1}$
with $A_{N}(W)$.  Since by construction,  elements of the form 
$o_{\widehat{\Y}}(w)(v+O_{N}(V))+J_{1}$ for $w\in W$ and 
$v\in V$ is in $G_{N}(\widehat{W}_{3}/J_{1})$, we have 
$D\subset G_{N}(\widehat{W}_{3}/J_{1})$.

Let $J_{2}$ be the generalized $V$-submodule of 
$\widehat{W}_{3}/J_{1}$ generated by the coefficients of 
the formal series of the form
\begin{align*}
&(x_{0}+x_{2})^{\swt v+N}\widehat{\Y}(Y_{W}(v, x_{0})w, x_{2})w_{2}
-(x_{0}+x_{2})^{\swt v+N}Y_{\widehat{W}_{3}}(v, x_{0}+x_{2})\widehat{\Y}(w, x_{2})w_{2}+J_{1}
\end{align*}
for homogeneous $v\in V$, $w\in W$ and $w_{2}\in A_{N}(V)$.
Let $W_{3}=(\widehat{W}_{3}/J_{1})/J_{2}$. 
Then $W_{3}$ is a lower-bounded generalized $V$-module.
Let 
\begin{align*}
\Y: W\otimes W_{2}&\to W_{3}\{x\}[\log x]\nn
w\otimes w_{2}&\mapsto \Y(w, x)w_{2}
\end{align*}
be a linear map defined by 
$$\Y(w, x)w_{2}=(\widehat{\Y}(w, x)w_{2}+J_{1})+J_{2}$$
for $w\in W$ and $w_{2}\in W_{2}$. 
By the definitions of $J_{1}$ and $J_{2}$, we see that 
$\Y$ satisfies the lower-truncation property, the 
$L(-1)$ property, the commutator formula for one intertwining operator
and the weak associativity for one intertwining operator when 
acting on $A_{N}(V)$. The commutativity for 
one intertwining operator and generalized rationality for one 
intertwining operator follows from the commutator formula for one 
intertwining operator. Using this commutativity, we obtain the 
weak associativity for one intertwining operator acting on 
$W_{2}$. The weak associativity for one intertwining operator
gives the associativity for one
intertwining operator. 
Since the lower-truncation property, the $L(-1)$-derivative property,
the generalized rationality, commutativity 
and associativity for one intertwining operator holds for $\Y$, 
we see that 
$\Y$ is an intertwining operator of type 
$\binom{W_{3}}{WW_{2}}$.

We want to show that $D\cap J_{2}=0$ and then we can view
$D$ as a subspace of $G_{N}(W_{3})$.
We consider the graded dual space $(\widehat{W}_{3}/J_{1})'$ with 
respect to the $\N$-grading of $\widehat{W}_{3}/J_{1}$.
Given an element $d^*\in D^{*}$, we extend it to an element of 
$G_{N}((\widehat{W}_{3}/J_{1})')$ as follows: For $w\in W$, $v\in V$,
$w_{2}\in A_{N}(V)$ and $z_{1}, z_{2}$ satisfying 
$|z_{2}|>|z_{1}-z_{2}|>0$, 
$$\langle d^{*}, 
\widehat{\Y}(Y_{W}(v, z_{1}-z_{2})w, z)w_{2}+J_{1}\rangle$$
is well defined. On the other hand, 
since the commutator formula for 
$Y_{\widehat{W}_{3}}(v, z_{1})$ and $\widehat{\Y}(w, z_{2})$ 
holds modulo $J_{1}$, 
$Y_{\widehat{W}_{3}}(v, z_{1})\widehat{\Y}(w, z_{2})w_{2}
+J_{1}$ is absolutely convergent to an element of the algebraic completion 
of $\widehat{W}_{3}/J_{1}$. 
We define 
$$\langle d^{*}, Y_{\widehat{W}_{3}}(v, z_{1})\widehat{\Y}(w, z_{2})w_{2}
+J_{1}\rangle
=\langle d^{*}, 
\widehat{\Y}(Y_{V}(v, z_{1}-z_{2})w, z)w_{2}+J_{1}\rangle$$
for $z_{1}, z_{2}\in \C$ satisfying $|z_{1}|>|z_{2}|>|z_{1}-z_{2}|>0$.
Since the homogeneous components of 
$Y_{\arc{W}_{3}}(v, z_{1})\arc{\Y}(w, z_{2})w_{2}
+J_{1}$ for $v\in V$, $w\in W$ and
$w_{2}\in A_{N}(V)$ span $\widehat{W}_{3}/J_{1}$,
$d^{*}$ gives an element of $G_{N}((\widehat{W}_{3}/J_{1})')$. 
Thus we can identify $D^{*}$ with a subspace of 
$G_{N}((\widehat{W}_{3}/J_{1})')$. 

We define a subspace $J_{3}$ of $\widehat{W}_{3}/J_{1}$
to be the subspace annihilated by $D^{*}$, that is
$$J_{3}=\{\hat{w}_{3}+J_{1}\mid \hat{w}_{3}\in \hat{W}_{3},\;
\langle d^{*}, 
\hat{w}_{3}+J_{1}\rangle=0\;\text{for}\; d^{*}\in D^{*}\}.$$
We now show that $J_{2}\subset J_{3}$.
The space $J_{2}$ is spanned by the coefficients of the formal series
\begin{align}\label{span-J-2}
&(x_{0}+x_{2})^{\swt v+N}Y_{\widehat{W}_{3}}(u, x)
\widehat{\Y}(Y_{W}(v, x_{0})w, x_{2})w_{2}\nn
&\quad -(x_{0}+x_{2})^{\swt v+N}Y_{\widehat{W}_{3}}(u, x)
Y_{\widehat{W}_{3}}(v, x_{0}+x_{2})\widehat{\Y}(w, x_{2})w_{2}+J_{1}
\end{align}
for $u, v\in V$ (with $v$ being homogeneous), 
$w\in W$ and $w_{2}\in A_{N}(V)$. 
When we substitute $z$, $z_{1}-z_{2}$ and $z_{2}$
for $x$, $x_{0}$ and $x_{2}$, where $z$, $z_{1}$ and $z_{2}$ 
are complex numbers satisfying $|z|>|z_{1}|>|z_{2}|>|z_{1}-z_{2}|>0$,
(\ref{span-J-2}) is absolutely 
convergent to element
\begin{align}\label{span-J-2-1}
&z_{1}^{\swt v+N}Y_{\widehat{W}_{3}}(u, z)
\widehat{\Y}(Y_{W}(v, z_{1}-z_{2})w, z_{2})w_{2}
 -z_{1}^{\swt v+N}Y_{\widehat{W}_{3}}(u, z)
Y_{\widehat{W}_{3}}(v, z_{1})\widehat{\Y}(w, z_{2})w_{2}+\overline{J}_{1}
\end{align}
of the algebraic completion of 
$\widehat{W}_{3}/J_{1}$, where $\overline{J}_{1}$ is 
the algebraic completion of $J_{1}$. 
Moreover, the homogeneous components of these elements in 
the algebraic completion of 
$\widehat{W}_{3}/J_{1}$ also span $J_{2}$. 
For $d^{*}\in D^{*}$, by the definition of $d^{*}$ as an element of 
$G_{N}((\widehat{W}_{3}/J_{1})')$ 
and the associativity
for the vertex operator map $Y_{W}$ and $Y_{\widehat{W}_{3}}$, 
we have
\begin{align}\label{analytic-ext}
&\langle d^{*}, z_{1}^{\swt v+N}Y_{\widehat{W}_{3}}(u, z)
\widehat{\Y}(Y_{V}(v, z_{1}-z_{2})w, z_{2})w_{2}+\overline{J}_{1}\rangle\nn
&\quad\sim z_{1}^{\swt v+N}\langle d^{*}, 
\widehat{\Y}(Y_{W}(u, z-z_{2})Y_{W}(v, z_{1}-z_{2})w, z_{2})w_{2}
+\overline{J}_{1}\rangle\nn
&\quad\sim z_{1}^{\swt v+N}\langle d^{*}, 
\widehat{\Y}(Y_{W}(Y_{V}(u, z-z_{1})v, z_{1}-z_{2})w, z_{2})w_{2}
+\overline{J}_{1}\rangle\nn
&\quad\sim z_{1}^{\swt v+N}\langle d^{*}, 
Y_{\widehat{W}_{3}}(Y_{V}(u, z-z_{1})v, z_{1})
\widehat{\Y}(w, z_{2})w_{2}+\overline{J}_{1}\rangle\nn
&\quad\sim \langle d^{*}, z_{1}^{\swt v+N}
Y_{\widehat{W}_{3}}(u, z)Y_{\widehat{W}_{3}}(v, z_{1})
\widehat{\Y}(w, z_{2})w_{2}+\overline{J}_{1}\rangle,
\end{align}
where $\sim$ means ``can be analytically extended to.'' Note that 
$z_{2}$ can be fixed for each of the analytic extension step 
in (\ref{analytic-ext}). So the analytic extensions 
in  (\ref{analytic-ext}) do not change the value of 
$\log z_{2}$. Since the left-hand side and right-hand side of 
(\ref{analytic-ext}) are convergent absolutely in the region 
$|z|>|z_{2}|>|z_{1}-z_{2}|>0$ and $|z|>|z_{1}|>|z_{2}|>0$, 
respectively, we see that in the region
$|z|>|z_{1}|>|z_{2}|>|z_{1}-z_{2}|>0$,
we have 
\begin{align}\label{d*-assoc}
&\langle d^{*}, z_{1}^{\swt v+N}Y_{\widehat{W}_{3}}(u, z)
\widehat{\Y}(Y_{V}(v, z_{1}-z_{2})w, z_{2})w_{2}+\overline{J}_{1}
\rangle\nn
&\quad =\langle d^{*}, z_{1}^{\swt v+N}
Y_{\widehat{W}_{3}}(u, z)Y_{\widehat{W}_{3}}(v, z_{1})
\widehat{\Y}(w, z_{2})w_{2}+\overline{J}_{1}\rangle.
\end{align}
For such $z, z_{1}$ and $z_{2}$, we can rewrite (\ref{d*-assoc})
as 
\begin{align}\label{d*-assoc-1}
&\left\langle d^{*}, \left(z_{1}^{\swt v+N}Y_{\widehat{W}_{3}}(u, z)
\widehat{\Y}(Y_{V}(v, z_{1}-z_{2})w, z_{2})\right.\right.\nn
&\quad\quad\quad\quad\quad
\left.\left.- z_{1}^{\swt v+N}
Y_{\widehat{W}_{3}}(u, z)Y_{\widehat{W}_{3}}(v, z_{1})
\widehat{\Y}(w, z_{2})w_{2}+\overline{J}_{1}\right)\right\rangle=0.
\end{align}
Since $J_{2}$ is spanned by the homogeneous components of 
elements of the form (\ref{span-J-2-1}), we see from 
(\ref{d*-assoc-1}) that $J_{2}\subset J_{3}$. 
Then $D\cap J_{2}\subset D\cap J_{3}=0$. So 
we can view $D$ as a subspace of $G_{N}(W_{3})$.

From Proposition 5.7 in \cite{HY}, we have an $A_{N}(V)$-module map 
\[
\rho(\mathcal {Y}): A_{N}(W)\otimes_{A_{N}(V)} \Omega^{0}_N(W_2)
\rightarrow \Omega_{N}^{0}(W_{3})
\]
defined by 
\begin{align*}
\rho(\mathcal {Y})((w+O_{N}(W))\otimes w_{2})
& = \sum_{n=0}^{N}{\rm
Res}_x x^{-h_3 -n-
1}\mathcal{Y}^{0}(x^{L_{W}(0)_{s}}w,
x)x^{L_{W_{2}}(0)_{s}}w_{2}\nn
& = \sum_{n=0}^{N}\Y_{\swt w + \swt w_{2} -
h_3 -n - 1, 0}(w)w_{2},\nn
\end{align*}
for  homogeneous $w\in W$, $w_{(2)} \in \Omega^{0}_N(W_2)$.
Note that $G_{N}(W_{2})$ is an $A_{N}(V)$-submodule 
of $\Omega_{N}^{0}(W_{2})$ and 
$A_{N}(W)\otimes_{A_{N}(V)}G_{N}(W_{2})$ is an 
$A_{N}(V)$-submodule of 
$A_{N}(W)\otimes_{A_{N}(V)}\Omega^{0}_{N}(W_{2})$. Also 
the projection $\pi_{G_{N}(W_{3})}$ from 
$\Omega_{N}^{0}(W_{3})$ to
$G_{N}(W_{3})$ is an $A_{N}(V)$-module map. Then 
we obtain an $A_{N}(V)$-module map
$$f=\pi_{G_{N}(W_{3})}\circ \rho(\Y)\circ 
e_{G_{N}(W_{2})}:
A_{N}(W)\otimes_{A_{N}(V)}G_{N}(W_{2})\to G_{N}(W_{3}),$$
where $e_{G_{N}(W_{2})}$
is the embedding map from $A_{N}(W)\otimes_{A_{N}(V)}G_{N}(W_{2})$
 to $A_{N}(W)\otimes_{A_{N}(V)}\Omega^{0}_{N}(W_{2})$.
Since $G_{N}(W_{2})=A_{N}(V)$, we see that 
$$A_{N}(W)\otimes_{A_{N}(V)}G_{N}(W_{2})=A_{N}(W)\otimes_{A_{N}(V)}(\one +O_{N}(V))$$
is equivalent as a left $A_{N}(V)$-module to $A_{N}(W)$. In particular, 
$f$ can be viewed as an $A_{N}(V)$-module map from
$A_{N}(W)$ to $G_{N}(W_{3})$. 

If $f(w+O_{N}(W))=0$ for 
homogeneous $w\in W$, 
by the definitions of $o_{\widehat{\Y}}$, $\widehat{\Y}$,
$\Y$ and $f$, we have in $W_{3}$
\begin{align*}
&(o_{\widehat{\Y}}(w)(\one+O_{N}(V))+J_{1})+J_{2}\nn
&\quad =(w(\wt w  +\wt (\one+O_{N}(V))-
h_3  -N-1, 0)(\one+O_{N}(V))+J_{1})+J_{2}\nn
&\quad =(\widehat{\Y}_{\swt w +\swt (\one+O_{N}(V))  -
h_3 -N- 1, 0}(w)(\one+O_{N}(V))+J_{1})+J_{2}\nn
&\quad =\Y_{\swt w +\swt (\one+O_{N}(V))  -
h_3 -N- 1, 0}(w)(\one+O_{N}(V))\nn
&\quad =\pi_{G_{N}(W_{3})}\sum_{n=0}^{N}\Y_{\swt w 
+\swt (\one+O_{N}(V))  -
h_3 -n- 1, 0}(w)(\one+O_{N}(V))\nn
&\quad=(\pi_{G_{N}(W_{3})}\circ \rho(\Y)\circ 
e_{G_{N}(W_{2})})((w+O_{N}(W))\otimes_{A_{N}(V)}
(\one+O_{N}(V))\nn
&\quad =f(w+O_{N}(W))\nn
&\quad=0.
\end{align*}
Since $o_{\widehat{\Y}}(w)(\one+O_{N}(V))+J_{1}\in D$ and 
$D\cap J_{2}=0$, we obtain 
$o_{\widehat{\Y}}(w)(\one+O_{N}(V))+J_{1}=0$
in $\widehat{W}_{3}/J_{1}$ or equivalently, 
$o_{\widehat{\Y}}(w)(\one+O_{N}(V))\in J_{1}$. 
We have shown above that $o_{\widehat{\Y}}(w)(\one+O_{N}(V))\in J_{1}$
if and only if $w\in O_{N}(W)$. In summary,
we have shown that $f(w+O_{N}(W))=0$
implies  $w\in O_{N}(W)$.

By definition, $Q^{\infty}(W)\subset \ker \vartheta_{\Y}$. 
In particular, for $\sum_{k=0}^{N}[w]_{kk}\in Q^{NN}(W)$,
\begin{align*}
f(w+O_{N}(W))
&=\Y_{\swt w +\swt (\one+O_{N}(V))  -
h_3 -N- 1, 0}(w)(\one+O_{N}(V))\nn
&=\res_{x}x^{\swt (\one+O_{N}(V))-h_{3}-N-1}
\Y^{0}(x^{L_{W}(0)_{S}}w, x)(\one+O_{N}(V))\nn
&=\vartheta_{\Y}([w]_{NN})(\one+O_{N}(V))\nn
&=\vartheta_{\Y}\left(\sum_{k=0}^{N}[w]_{kk}\right)(\one+O_{N}(V))\nn
&=0.
\end{align*}
Then we have $w\in O_{N}(W)\subset O_{k}(V)$
for $k=0, \dots, N$ or equivalently, $\sum_{k=0}^{N}[w]_{kk}\in O^{NN}(W)$. 
Thus $Q^{NN}(W)=O^{NN}(W)$. 
\epfv

\begin{thm}\label{isom-thm-A-N-W}
The space $A^{NN}(W)$ is invariant 
under the left and right actions of $U^{NN}(V)$ with the 
left and right actions of $Q^{NN}(V)$ on $A^{NN}(W)$
equal to $0$  and is thus 
an $A^{NN}(V)$-bimodule. Moreover, the linear map 
$f^{NN}: U^{NN}(W)\to A_{N}(W)$
defined by $f_{NN}(\sum_{k=0}^{N}[w]_{kk})=w+O_{N}(W)$ for $w\in W$
induces an invertible linear map, still denoted by $f_{NN}$, sending 
the $A^{NN}(V)$-bimodule structure on 
$A^{NN}(W)$ to the $A_{N}(V)$-bimodule structure on $A_{N}(W)$. 
\end{thm}
\pf
By  the definitions of the left action and right action of $U^{\infty}(V)$ 
on $U^{\infty}(W)$ (see Section 4 in \cite{H-aa-int-op}), we have
\begin{align*}
&\left(\sum_{k=0}^{N}[v]_{kk}\right)\diamond 
\left(\sum_{k=0}^{N}[w]_{kk}\right)\nn
& \quad=\sum_{k, l=0}^{N}[v]_{kk}\diamond [w]_{ll}\nn
&  \quad=\sum_{k=0}^{N}[v]_{kk}\diamond [w]_{kk}
\nn
& \quad=\sum_{k=0}^{N}
\res_{x}T_{2k+1}((x+1)^{-k-1})(1+x)^{k}\left[Y_{W}((1+x)^{L_{V}(0)}v, x)w\right]_{kk}\nn
& \quad=\sum_{k=0}^{N}
\res_{x}\sum_{m=0}^{k}\binom{-k-1}{m}x^{-k-m-1}
(1+x)^{k}\left[Y_{V}((1+x)^{L(0)}u, x)v\right]_{kk}\nn
&\quad=\sum_{k=0}^{N}[v*_{k}w]_{kk}\nn
&\quad\equiv \sum_{k=0}^{N}[v*_{N}w]_{kk} \mod Q^{\infty}(W)
\nn
&\quad\in U^{NN}(W)
\end{align*}
for $v\in V$ and $w\in W$, 
where for the step given by $\equiv$, we have used 
$v*_{k} w$ is equal to $v *_{N} w$ modulo $O_{k}(W)$ which can be 
proved using the same argument as in the proof of 
Proposition 2.4 in \cite{DLM} and 
the fact that $[O_{k}(W)]_{kk}\in O^{\infty}(W)\subset Q^{\infty}(W)$.
Similarly, we have 
\begin{align*}
&\left(\sum_{k=0}^{N}[w]_{kk}\right)\diamond 
\left(\sum_{k=0}^{N}[v]_{kk}\right)\nn
& \quad=\sum_{k, l=0}^{N}[w]_{kk}\diamond [v]_{ll}\nn
& \quad =\sum_{k=0}^{N}[w]_{kk}\diamond [v]_{kk}
\nn
& \quad=\sum_{k=0}^{N}
\res_{x}T_{2k+1}((x+1)^{-k-1})(1+x)^{k}
\left[Y_{W}((1+x)^{-L_{V}(0)}u, -x(1+x)^{-1})w\right]_{kk}\nn
&\quad =\sum_{k=0}^{N}
\res_{x}\sum_{m=0}^{k}\binom{-k-1}{m}x^{-k-m-1}
(1+x)^{k}\left[Y_{W}((1+x)^{-L_{V}(0)}u, -x(1+x)^{-1})w\right]_{kk}\nn
&\quad =\sum_{k=0}^{N}[w*_{k}v]_{kk}\nn
&\quad \equiv \sum_{k=0}^{N}[w*_{N}v]_{kk}\mod Q^{\infty}(W)\nn
&\quad\in U^{NN}(W)
\end{align*}
for $v\in V$ and $w\in W$, 
where for the step given by $\equiv$, we have used 
$w*_{k} v$ is equal to $w *_{N} v$ modulo $O_{k}(W)$ which can be 
proved using a similar argument as in the proof of 
Proposition 2.4 in \cite{DLM} and 
the fact that $[O_{k}(W)]_{kk}\in O^{\infty}(W)\subset Q^{\infty}(W)$.
Thus $A^{NN}(W)$  is invariant 
under the left and right actions of $U^{NN}(V)$. By Proposition
\ref{Q-NN=O-NN}, we see that the left and right actions
of $Q^{NN}(V)=O^{NN}(V)$ on $A^{NN}(W)$ are $0$.
So $A^{NN}(W)$  is an  $A^{NN}(V)$-bimodule. 

Let $f^{NN}: U^{NN}(W)\to A_{N}(W)$
be defined by 
$$f^{NN}\left(\sum_{k=0}^{N}[w]_{kk}\right)=w+O_{N}(W)$$
for $w\in W$. 
Then by the definition of $O_{N}(W)$, we have 
\begin{align*}
&f^{N}\left(\sum_{k=0}^{N}[\res_{x}
x^{-2N-p-2}(1+x)^{N}Y_{W}((1+x)^{L_{V}(0)}v, x)w]_{kk}\right)\nn
&\quad =\res_{x}
x^{-2N-p-2}(1+x)^{N}Y_{W}((1+x)^{L_{V}(0)}v, x)w+O_{N}(W)\nn
&\quad =0.
\end{align*}
Thus by the definition of $O^{\infty}(W)$, we obtain 
$O^{NN}(W)=
O^{\infty}(W)\cap U^{NN}(W)\subset \ker f^{N}$. On the other hand, 
if $w\in O_{N}(W)$, then $w\in O_{k}(W)$ for $k=0, \dots, N$
since $O_{N}(W)\subset O_{k}(W)$. Hence $[w]_{kk}\in O^{NN}(W)$
for $k=0, \dots, N$ and thus $\sum_{k=0}^{N}[w]_{kk}\in 
O^{NN}(W)$. Equivalently, 
if $\sum_{k=0}^{N}[w]_{kk}\not\in 
O^{NN}(W)$,  then $w\not\in O_{N}(W)$.
In particular,  in this case, $f^{NN}(\sum_{k=0}^{N}[w]_{kk})
=w+O_{N}(W)\ne 0$, that is, 
$\sum_{k=0}^{N}[w]_{kk}\not\in \ker f_{NN}$. 
Thus $\ker f^{NN}=O^{NN}(W)$.
By Proposition \ref{Q-NN=O-NN}, $\ker f^{NN}=Q_{NN}(W)$.
It is clear that $f^{NN}$ is surjective. 
In particular, $f^{NN}$ induces a linear isomorphism, still denoted by $f^{NN}$,
from $A_{NN}(W)$ to $A_{N}(W)$. 

For $v\in V$ and $w\in W$, 
\begin{align*}
f&^{NN}\left(\left(\sum_{k=0}^{N}[v]_{kk}+O^{NN}(V)\right)\diamond 
\left(\sum_{k=0}^{N}[w]_{kk}+O^{NN}(W)\right)\right)\nn
&=f^{NN}\left(\sum_{k=0}^{N}[v*_{N}w]_{kk}+O^{NN}(W)\right)\nn
&=v*_{N} w+O_{N}(W)\nn
&=(v+O_{N}(V))*_{N} (w+O_{N}(W)).
\end{align*}
Therefore $f_{NN}$ sends
the $A^{NN}(V)$-bimodule structure on 
$A^{NN}(W)$ to the $A_{N}(V)$-bimodule structure on $A_{N}(W)$. 
\epfv

\begin{cor}
For $N\in \N$, the space $Q^{NN}(W)=Q^{\infty}(W)\cap U^{NN}(W)$
is spanned by elements of the form 
$$ \sum_{k=0}^{N}\res_{x}
x^{-2N-p-2}(1+x)^{l}[Y_{W}((1+x)^{L_{V}(0)}v, x)w]_{kk}$$
for $v\in V$, $w\in W$ and $p\in \N$.
\end{cor}

We now introduce the associative algebras 
$\widetilde{A}^{\infty}(V)$ and $\widetilde{A}^{N}(V)$ for $N\in \N$. 
Recall from \cite{H-aa-va} that $U^{\infty}(V)$ is
the space of  column-finite infinite
matrices with entries in $V$, but doubly indexed by $\N$ instead of $\Z_{+}$. 
Recall also that for $v\in V$, $k, l\in \N$, $[v]_{kl}$ is the element of $U^{\infty}(V)$
with the $(k, l)$-entry being $v$ and 
all the other entries being $0$. 
Elements of $U^{\infty}(V)$ are suitable (possibly infinite) sums of elements of the form 
$[v]_{kl}$ for $v\in V$, $k, l\in \N$.

In \cite{H-aa-va}, a product $\diamond$ on $U^{\infty}(V)$ is introduced. 
Here we define a new product $\bkdiamd$ on $U^{\infty}(V)$ by 
$$[u]_{km}\bkdiamd [v]_{nl}=0$$
for $k, l, m, n\in \N$ and $u, v\in V$ when $m\ne n$ and 
\begin{align*}%\label{defn-bkdiamd-1}
&[u]_{kn}\bkdiamd [v]_{nl}\nn
&\quad =\res_{x}T_{k+l+1}((x+1)^{-k+n-l-1})(1+x)^{l}
\left[Y_{V}\left(u, \frac{1}{2\pi i}\log (1+x)\right)v\right]_{kl}\nn
&\quad =\sum_{m=0}^{n}
\binom{-k+n-l-1}{m}\res_{x}
x^{-k+n-l-m-1}(1+x)^{l}\left[Y_{V}\left(u, \frac{1}{2\pi i}\log (1+x)\right)v\right]_{kl}\nn
&\quad =\sum_{m=0}^{n}
\binom{-k+n-l-1}{m}\res_{y}2\pi ie^{2pi (l+1)i}
(e^{2\pi iy}-1)^{-k+n-l-m-1}
\left[Y_{V}\left(u, y\right)v\right]_{kl}
\end{align*}
for $k, l, n\in \N$ and $u, v\in V$. Then $U^{\infty}(V)$ equipped with 
$\bkdiamd$ is a nonassociative algebra.

For $W=V$, we have the invertible linear map
$\mathcal{U}_{V}(1)=(2\pi i)^{L_{V}(0)}e^{-L^{+}_{V}(A)}: V\to V$ 
(see \cite{H-mod-inv-int} and Subsectipon 2.2).
We extend the linear isomorphism $\mathcal{U}_{V}(1): V\to V$
to a linear isomorphism $\mathcal{U}_{V}(1): U^{\infty}(V)\to U^{\infty}(V)$ by
$$\mathcal{U}_{V}(1)[v]_{kl}=[\mathcal{U}_{V}(1)v]_{kl}$$
for $k, l\in \N$ and $v\in V$. 

Recall the subspace $Q^{\infty}(V)$ of $U^{\infty}(V)$ in \cite{H-aa-va}. 
Let $\widetilde{Q}^{\infty}(V)=\mathcal{U}_{V}(1)^{-1}Q^{\infty}(V)$ and 
$\widetilde{A}^{\infty}(V)=U^{\infty}(V)/\widetilde{Q}^{\infty}(V)$. 
Recall from \cite{H-aa-va} that $A^{\infty}(V)=U^{\infty}(V)/Q^{\infty}(V)$ with the product induced by 
$\diamond$ is an associative algebra.

\begin{prop}\label{tilde-A-infty}
The linear isomorphism $\mathcal{U}_{V}(1)$ from 
$U^{\infty}(V)$ to itself 
is an isomorphism of nonassociative algebras from $U^{\infty}(V)$ equipped with the product $\bkdiamd$ 
to $U^{\infty}(V)$ equipped with the product 
$\diamond$ such that $\mathcal{U}_{V}(1)\widetilde{Q}^{\infty}(V)
=Q^{\infty}(V)$. 
In particular, $\bkdiamd$ induces an associative product, still denoted by $\bkdiamd$,
on $\widetilde{A}^{\infty}(V)$ such that $\mathcal{U}_{V}(1)$ induces an 
isomorphism of associative algebras from $\widetilde{A}^{\infty}(V)$
to $A^{\infty}(V)$.
\end{prop}
\pf 
Using (\ref{U-1-v-op}), we obtain
$$\mathcal{U}_{V}(1)([u]_{km}\bkdiamd [v]_{nl})=0=[u]_{km}\diamond [v]_{nl}$$
for $k, l, m, n\in \N$ and $u, v\in V$ when $m\ne n$ and 
\begin{align*}
&\mathcal{U}_{V}(1)([u]_{kn}\bkdiamd [v]_{nl})\nn
&\quad=\mathcal{U}_{V}(1)\left(\res_{x}T_{k+l+1}((x+1)^{-k+n-l-1})(1+x)^{l}
\left[Y_{V}\left(u, \frac{1}{2\pi i}\log (1+x)\right)v\right]_{kl}\right)\nn
&\quad=\res_{x}T_{k+l+1}((x+1)^{-k+n-l-1})(1+x)^{l}\left[Y_{V}((1+x)^{L_{V}(0)}
\mathcal{U}_{V}(1)u, x)\mathcal{U}_{V}(1)v\right]_{kl}\nn
&\quad=[\mathcal{U}_{V}(1)u]_{kn}\diamond [\mathcal{U}_{V}(1)v]_{nl}
\end{align*}
for $k, l, n\in \N$, $u, v\in V$. These show that $\mathcal{U}_{V}(1)$ is an isomorphism of the 
nonassociative algebras.
By the definition of $\widetilde{Q}^{\infty}(V)$, we have $\mathcal{U}_{V}(1)\widetilde{Q}^{\infty}(V)
=Q^{\infty}(V)$. The other conclusions follow immediately. 
\epfv

For a lower-bounded generalized $V$-module $W$, we have a graded $A^{\infty}(V)$-module
structure on $W$ given by $\vartheta_{W}: U^{\infty}(V)\to \text{End}\;W$ (see \cite{H-aa-int-op}). 
Let $\tilde{\vartheta}_{W}: U^{\infty}(V)\to \text{End}\;W$ be defined by 
$$\tilde{\vartheta}_{W}(\mathfrak{v})=\vartheta_{W}(\mathcal{U}_{V}(1)\mathfrak{v})$$
for $\mathfrak{v}\in U^{\infty}(V)$. 

\begin{prop}\label{tilde-a-mod}
Let $W$ be a lower-bounded generalized $V$-module. 
Then the linear map $\tilde{\vartheta}_{W}$ gives $W$ 
 an $\widetilde{A}^{\infty}(V)$-module structure. 
\end{prop}
\pf
For $k, l, m, n\in \N$ and $u, v\in V$,
\begin{align*}
\tilde{\vartheta}_{W}([u]_{km}\bkdiamd [v]_{nl})
&=\vartheta_{W}(\mathcal{U}_{V}(1)([u]_{km}\bkdiamd [v]_{nl}))\nn
&=\vartheta_{W}([\mathcal{U}_{V}(1)u]_{kn}\diamond [\mathcal{U}_{V}(1)v]_{nl})\nn
&=\vartheta_{W}([\mathcal{U}_{V}(1)u]_{kn})\vartheta_{W}([\mathcal{U}_{V}(1)v]_{nl})\nn
&=\tilde{\vartheta}_{W}([u]_{km})\tilde{\vartheta}_{W}([v]_{nl}).
\end{align*}
Thus $W$ equipped with $\tilde{\vartheta}_{W}$
is an $\widetilde{A}^{\infty}(V)$-module.
\epfv

From Theorem 5.1 in \cite{H-aa-int-op}, we know that the category of lower-bounded generalized $V$-modules is isomorphic 
to the category of graded $A^{\infty}(V)$-modules. We now define a {\it graded 
$\widetilde{A}^{\infty}(V)$-module} to be an $\widetilde{A}^{\infty}(V)$-module obtained 
from a lower-bounded generalized $V$-module as in Proposition \ref{tilde-a-mod}. Or equivalently,
a graded  $\widetilde{A}^{\infty}(V)$-module is an
$\widetilde{A}^{\infty}(V)$-module $W$ with the $\widetilde{A}^{\infty}(V)$-module structure given by 
a linear map $\tilde{\vartheta}_{W}: U^{\infty}(V)\to \text{End}\;W$ such that 
$\vartheta_{W}: U^{\infty}(V)\to \text{End}\;W$ defined by 
$\vartheta_{W}(\mathfrak{v})=\tilde{\vartheta}_{W}(\mathcal{U}_{V}(1)^{-1}\mathfrak{v})$ 
gives a graded $A^{\infty}(V)$-module structure to $W$. 

Recall the subalgebras $A^{NN}(V)$ for $N\in \N$ in \cite{H-aa-va} 
which are proved 
in  \cite{H-aa-va} to be isomorphic to 
the associative algebras $A_{N}(V)$ introduced in \cite{DLM}. 
In this paper, we need the following 
subalgebras $\widetilde{A}^{NN}(V)$ for $N\in \N$:
For $N\in \N$, let 
$$\widetilde{A}^{NN}(V)=\left\{\sum_{k=0}^{N}[v]_{kk}
+\widetilde{Q}^{\infty}(V)\;|\;v\in V\right\}.$$

\begin{prop}\label{tilde-A-NN}
For $N\in \N$, $\widetilde{A}^{NN}(V)$ is a subalgebra of 
$\widetilde{A}^{\infty}(V)$ and 
$\mathcal{U}_{V}(1)$ induces an isomorphism of associative algebras from  
$\widetilde{A}^{NN}(V)$ to $A^{NN}(V)$. For a 
lower-bounded generalized $V$-module $W$, $\Omega_{N}^{0}(W)
=\coprod_{\mu\in 
\Gamma(W)}\coprod_{k=0}^{N}W_{[h^{\mu}+N]}$ 
is invariant under the 
action of $\widetilde{A}^{NN}(V)$ and thus is a
$\widetilde{A}^{NN}(V)$-module. 
\end{prop}
\pf
By definition, $\mathcal{U}_{V}(1)\widetilde{A}^{NN}(V)=A^{NN}(V)$.
Since $\mathcal{U}_{V}(1)$ is an isomorphsim from the 
associative algebra $\widetilde{A}^{\infty}(V)$
to the associative algebra $A^{\infty}(V)$ and 
$A^{NN}(V)$ is a subalgebra of $A^{\infty}(V)$,
$\widetilde{A}^{NN}(V)$ is a subalgebra of $\widetilde{A}^{\infty}(V)$
and $\mathcal{U}_{V}(1)$ restricted to $\widetilde{A}^{NN}(V)$
is an isomorphism from $\widetilde{A}^{NN}(V)$ to $A^{NN}(V)$. 

By the definition of $\tilde{\vartheta}_{W}$, 
$\Omega_{N}^{0}(W)$  is invariant under the action of $\widetilde{A}_{NN}(V)$.
\epfv

We also introduce another associative algebra $\widetilde{A}_{N}(V)$ 
generalizing the associative algebra $\widetilde{A}(V)$ in 
\cite{H-mod-inv-int}. Define a product $\bullet_{N}$ of $V$ by
$$u\bullet_{N}v= \res_{x}T_{2N+1}((x+1)^{-N-1})(1+x)^{N}
Y_{V}\left(u, \frac{1}{2\pi i}\log (1+x)\right)v$$
for $u, v\in V$. Let $\widetilde{O}_{N}(V)$ be the subspace of 
$V$ spanned by elements of the form 
$$\res_{x}x^{-2N-2-n}(1+x)^{N}
Y_{V}\left(u, \frac{1}{2\pi i}\log (1+x)\right)v$$
for $n\in \N$ and $u, v\in V$ and of the form $L_{V}(-1)v$.
Let $\widetilde{A}_{N}(V)=V/\widetilde{O}_{N}(V)$. 

\begin{prop}\label{tilde-A-N}
The product $\bullet_{N}$ induces an associative algebra structure on 
$\widetilde{A}_{N}(V)$ isomorphic to $A_{N}(V)$. The operator 
$\mathcal{U}_{V}(1)$ on $V$ induces an isomorphism from 
$\widetilde{A}_{N}(V)$ to $A_{N}(V)$. In particular,  
$\mathcal{U}_{V}(1)^{-1}\omega
+\widetilde{O}^{\infty}(V)$ is in the center 
of $\widetilde{A}_{N}(V)$
\end{prop}
\pf
For $n\in \N$ and $u, v\in V$, by (\ref{U-1-v-op}),
\begin{align*}
&\mathcal{U}_{V}(1)\res_{x}x^{-2N-2-n}(1+x)^{N}
Y_{V}\left(u, \frac{1}{2\pi i}\log (1+x)\right)v\nn
&\quad =\res_{x}x^{-2N-2-n}(1+x)^{N}
Y_{V}\left((1+x)^{L_{V}(0)}\mathcal{U}_{V}(1)u, x\right)
\mathcal{U}_{V}(1)v.
\end{align*}
Also for $v\in V$, we have 
$$\mathcal{U}_{V}(1)L_{V}(-1)v=2\pi i(L_{V}(-1)+L_{V}(0))v$$
(this is (1.15) in \cite{H-mod-inv-int}). These
shows $\mathcal{U}_{V}(1)\widetilde{O}_{N}(V)=O_{N}(V)$. 
Therefore $\mathcal{U}_{V}(1)$ induces a linear isomorphism from 
$\widetilde{A}_{N}(V)$ to $A_{N}(V)$.

By (\ref{U-1-v-op}) again, we have 
\begin{align*}
\mathcal{U}_{V}(1)(u\bullet_{N}v)
&=\mathcal{U}_{V}(1) \res_{x}T_{2N+1}((x+1)^{-N-1})(1+x)^{N}
Y_{V}\left(u, \frac{1}{2\pi i}\log (1+x)\right)v\nn
& =\res_{x}T_{2N+1}((x+1)^{-N-1})(1+x)^{N}
Y_{V}\left((1+x)^{L_{V}(0)}u, x\right)v\nn
&=u*_{N}v
\end{align*}
for $u, v\in V$. So $\mathcal{U}_{V}(1)$ induces an isomorphism of associative 
algebras. 

Since $\omega+O_{N}(V)$ is in the center of $A_{N}(V)$
and $\mathcal{U}_{V}(1)^{-1}$ is an isomorphism from 
$A_{N}(V)$ to $\widetilde{A}_{N}(V)$, 
$\mathcal{U}_{V}(1)^{-1}\omega
+\widetilde{O}^{\infty}(V)$ is in the center of $\widetilde{A}_{N}(V)$.
\epfv

\begin{prop}
The associative algebras $\widetilde{A}_{N}(V)$ and 
$\widetilde{A}^{NN}(V)$ are isomorphic.
The map given by $v+\widetilde{O}_{N}(V)
\mapsto \sum_{k=0}^{N}[v]_{kk}+\widetilde{Q}^{\infty}(V)$
is an isomorphism from $\widetilde{A}_{N}(V)$ to
$\widetilde{A}^{NN}(V)$.
\end{prop}
\pf
This result follows from Propositions \ref{tilde-A-NN}, \ref{tilde-A-N}
and Theorem 4.2 in \cite{H-mod-inv-int}.
\epfv

We now introduce the main subalgebras of $\widetilde{A}^{\infty}(V)$ 
needed in this paper.
Recall from \cite{H-aa-va}
 the subspaces  $U^{N}(V)$ for $N\in \N$ of $U^{\infty}(V)$
 spanned by elements of the form 
$[v]_{kl}$ for $v\in V$, $k, l=0, \dots, N$.
These subspaces are invariant under the operator $\mathcal{U}_{V}(1)$
on $U^{\infty}(V)$.
Let
$$\widetilde{A}^{N}(V)=\{\mathfrak{v}+\widetilde{Q}^{\infty}(V)\;|\;\mathfrak{v}\in U^{N}(V)\}.$$
Note that $\widetilde{A}^{N}(V)$ is
spanned by elements of the form $[v]_{kl}+\widetilde{Q}^{\infty}(V)$ for $v\in V$ and $k, l=0, \dots, N$.

\begin{prop}\label{tilde-A^N}
For $N\in \N$, $\widetilde{A}^{N}(V)$ is a subalgebra of $\widetilde{A}^{\infty}(V)$ and 
$\mathcal{U}_{V}(1)$ induces an isomorphism of associative algebras from  
$\widetilde{A}^{N}(V)$ to $A^{N}(V)$. For a 
lower-bounded generalized $V$-module $W$, $\Omega_{N}^{0}(W)$ is invariant under the 
action of $\widetilde{A}^{N}(V)$ and thus is a
$\widetilde{A}^{N}(V)$-module. 
\end{prop}
\pf
This proposition follows from the invariance of 
$U^{N}(V)$ under $\mathcal{U}_{V}(1)$ and Proposition
\ref{tilde-A-infty}.
\epfv

We now construct an $\widetilde{A}^{\infty}(V)$-bimodule from 
a lower-bounded generalized $V$-module $W$. 
We define a left action of $U^{\infty}(V)$ on $U^{\infty}(W)$ 
by 
$$[v]_{km}\bkdiamd [w]_{nl}=0$$
for $v\in V$, $w\in W$ and $k, m, n, l\in \N$ when $m\ne n$ and 
\begin{align*}%\label{defn-bkdiamd-1}
[v]&_{kn}\bkdiamd [w]_{nl}\nn
&=\res_{x}T_{k+l+1}((x+1)^{-k+n-l-1})(1+x)^{l}
\left[Y_{W}\left(v, \frac{1}{2\pi i}\log (1+x)\right)w\right]_{kl} \nn
&=\sum_{m=0}^{n}
\binom{-k+n-l-1}{m}\res_{x}
x^{-k+n-l-m-1}(1+x)^{l}
\left[Y_{W}\left(v, \frac{1}{2\pi i}\log (1+x)\right)w\right]_{kl}\nn
& =\sum_{m=0}^{n}
\binom{-k+n-l-1}{m}\res_{y}2\pi ie^{2\pi i(l+1)y}(e^{2\pi iy}-1)^{-k+n-l-m-1}
\left[Y_{W}\left(v, y\right)w\right]_{kl}.
\end{align*}
for $v\in V$, $w\in W$ and $k, n, l\in \N$.
We define a right action of $U^{\infty}(V)$ on $U^{\infty}(W)$ 
by 
$$[w]_{km} \bkdiamd [v]_{nl}=0$$
for $v\in V$, $w\in W$ and $k, m, n, l\in \N$ when $m\ne n$ and 
\begin{align}\label{defn-bkdiamd-2}
[w]&_{kn} \bkdiamd [v]_{nl}\nn
& =\res_{x}T_{k+l+1}((x+1)^{-k+n-l-1})(1+x)^{k}
\left[Y_{W}\left(v, -\frac{1}{2\pi i}\log (1+x)\right)w\right]_{kl}\nn
& =\sum_{m=0}^{n}
\binom{-k+n-l-1}{m}\res_{x}
x^{-k+n-l-m-1}(1+x)^{k}
\left[Y_{W}\left(v, -\frac{1}{2\pi i}\log (1+x)\right)w\right]_{kl}\nn
& =\sum_{m=0}^{n}
\binom{-k+n-l-1}{m}\res_{y}2\pi ie^{2\pi i(k+1)y}(e^{2\pi iy}-1)^{-k+n-l-m-1}
\left[Y_{W}\left(v, -y\right)w\right]_{kl}.
\end{align}

We extend the invertible linear map $\mathcal{U}_{W}(1): W\to W$ 
to an invertible linear map $\mathcal{U}_{W}(1): U^{\infty}(W) \to U^{\infty}(W)$
by 
$$\mathcal{U}_{W}(1)[w]_{kl}=[\mathcal{U}_{W}(1)w]_{kl}$$
for $k, l\in \N$ and $w\in W$. Let $\widetilde{Q}^{\infty}(W)=\mathcal{U}_{W}(1)^{-1}Q^{\infty}(W)$ and 
$\widetilde{A}^{\infty}(W)=U^{\infty}(W)/\widetilde{Q}^{\infty}(W)$. 

In \cite{H-aa-int-op}, a left action and a right action (denoted 
by $\diamond$), a subspace $Q^{\infty}(W)$ of $U^{\infty}(W)$
and an $A^{\infty}(V)$-bimodule $A^{\infty}(W)$ are introduced. 

\begin{prop}\label{U-1-isom}
The linear isomorphism $\mathcal{U}_{W}(1)$ from 
$U^{\infty}(W)$ to itself sends the left and right actions 
of $U^{\infty}(V)$ on $U^{\infty}(W)$ given by $\diamond$ to the 
left and right actions given by $\bkdiamd$ such that $\mathcal{U}_{W}(1)\widetilde{Q}^{\infty}(W)
=Q^{\infty}(W)$. 
In particular, $\bkdiamd$ induces left and right actions, still denoted by $\bkdiamd$,
of $\widetilde{A}^{\infty}(V)$
on $\widetilde{A}^{\infty}(W)$ such that $\widetilde{A}^{\infty}(W)$ becomes an
$\widetilde{A}^{\infty}(V)$-bimodule and $\mathcal{U}_{W}(1)$ induces an invertible 
linear map sending the $A^{\infty}(V)$-bimodule structure on $A^{\infty}(W)$
to the $\widetilde{A}^{\infty}(V)$-bimodule structure on 
$\widetilde{A}^{\infty}(W)$.
\end{prop}
\pf 
Using (\ref{U-1-v-op}), we have 
$$\mathcal{U}_{W}(1)([v]_{km}\bkdiamd [w]_{nl})=0=[v]_{km}\diamond [w]_{nl}$$
for $k, l, m, n\in \N$, $v\in V$ and $w\in W$ when $m\ne n$ and 
\begin{align*}
&\mathcal{U}_{W}(1)([v]_{kn}\bkdiamd [w]_{nl})\nn
&\quad=\mathcal{U}_{W}(1)\left(\res_{x}T_{k+l+1}((x+1)^{-k+n-l-1})(1+x)^{l}
\left[Y_{W}\left(v, \frac{1}{2\pi i}\log (1+x)\right)w\right]_{kl}\right)\nn
&\quad=\res_{x}T_{k+l+1}((x+1)^{-k+n-l-1})(1+x)^{l}\left[Y_{W}((1+x)^{L_{V}(0)}
\mathcal{U}_{V}(1)v, x)\mathcal{U}_{W}(1)w\right]_{kl}\nn
&\quad=[\mathcal{U}_{V}(1)v]_{kn}\diamond [\mathcal{U}_{W}(1)w]_{nl}
\end{align*}
for $k, l, n\in \N$, $u, v\in V$. This shows that $\mathcal{U}_{W}(1)$ sends the 
left action of $U^{\infty}(V)$ on $U^{\infty}(W)$ given by 
$\bkdiamd$ to the 
left action given by $\diamond$. 

Using (\ref{U-1-v-op}) again, we have 
\begin{align*}
&\mathcal{U}_{W}(1)\res_{x}x^{-n}
Y_{W}\left(v, -\frac{1}{2\pi i}\log(1+x)\right)w\nn
&\quad =\mathcal{U}_{W}(1)\res_{x}x^{-n}
Y_{W}\left(v, \frac{1}{2\pi i}\log(1-x(1+x)^{-1})\right)w\nn
&\quad =\res_{x}x^{-n}
Y_{W}\left((1-x(1+x)^{-1})^{L_{V}(0)}\mathcal{U}_{V}(1)v, -x(1+x)^{-1}\right)\mathcal{U}_{W}(1)w\nn
&\quad =\res_{x}x^{-n}Y_{W}((1+x)^{-L_{V}(0)}\mathcal{U}_{V}(1)v, -x(1+x)^{-1})\mathcal{U}_{W}(1)w
\end{align*}
for $n\in \Z$, $v\in V$ and $w\in W$. Then we have 
$$\mathcal{U}_{W}(1)([w]_{km}\bkdiamd [v]_{nl})=0=[w]_{km}\diamond [v]_{nl}$$
for $k, l, m, n\in \N$, $v\in V$ and $w\in W$ when $m\ne n$ and 
\begin{align*}
&\mathcal{U}_{W}(1)([w]_{kn}\bkdiamd [v]_{nl})\nn
&\quad=\mathcal{U}_{W}(1)\left(\res_{x}T_{k+l+1}((x+1)^{-k+n-l-1})(1+x)^{l}
\left[Y_{W}\left(v, -\frac{1}{2\pi i}\log (1+x)\right)w\right]_{kl}\right)\nn
&\quad=\res_{x}T_{k+l+1}((x+1)^{-k+n-l-1})(1+x)^{l}\left[Y_{W}((1+x)^{-L_{V}(0)}
\mathcal{U}_{V}(1)v, -x(1+x)^{-1})\mathcal{U}_{W}(1)w\right]_{kl}\nn
&\quad=[\mathcal{U}_{W}(1)w]_{kn}\diamond [\mathcal{U}_{V}(1)v]_{nl}
\end{align*}
for $k, l, n\in \N$, $u, v\in V$. This shows that $\mathcal{U}_{W}(1)$ sends the 
right action of $U^{\infty}(V)$ on $U^{\infty}(W)$ given by 
$\bkdiamd$ to the 
right action given by $\diamond$. 
By the definition of $\widetilde{Q}^{\infty}(W)$, we have $\mathcal{U}_{W}(1)\widetilde{Q}^{\infty}(W)
=Q^{\infty}(W)$. 

The other conclusions follow from these. 
\epfv

For $N\in \N$, let 
\begin{align*}
\widetilde{A}^{NN}(W)
&=\left\{\sum_{k=0}^{N}[w]_{kk}+\widetilde{Q}^{\infty}(W)\;|\;w\in W
\right\}.
\end{align*}

\begin{prop}\label{tilde-A-NN-W}
For $N\in \N$, 
$\widetilde{A}^{NN}(W)$ is closed under the left and right actions of 
the subalgebra $\widetilde{A}^{NN}(V)$ of $\widetilde{A}^{\infty}(V)$ 
and thus has an $\widetilde{A}^{NN}(V)$-bimodule structure. 
The map 
$\mathcal{U}_{W}(1)$ induces an invertible linear map sending 
the $\widetilde{A}^{NN}(V)$-bimodule structure on 
$\widetilde{A}^{NN}(W)$ 
to the $A^{NN}(V)$-bimodule structure on $A^{NN}(W)$. 
\end{prop}
\pf
This result follows immediately from the definitions
of the left and right actions of $\widetilde{A}^{\infty}(V)$ on 
$\widetilde{A}^{\infty}(W)$ and Proposition \ref{U-1-isom}. 
\epfv

\begin{cor}
The map given by $w+O_{N}(W)\mapsto \sum_{k=0}^{N}[w]_{kk}
+\widetilde{Q}^{\infty}(W)$
is an invertible linear map sending the 
$A_{N}(V)$-bimodule structure on 
$A_{N}(W)$ 
to the $\widetilde{A}^{NN}(V)$-bimodule structure on 
$\widetilde{A}^{NN}(W)$. 
\end{cor}

\begin{cor}\label{tilde-Q-NN}
The space $\widetilde{Q}^{\infty}(W)\cap U^{NN}(W)$ is spanned 
by the coefficients of 
$$\sum_{k=0}\res_{x}e^{2\pi i(N+1)x}(e^{2\pi ix}-1)^{-2N-p-2}
[Y_{W}(v, x)w]_{kk}$$
for $p\in \N$, $v\in V$ and $w\in W$. 
\end{cor}

Recall from \cite{H-aa-int-op} the subspaces $U^{N}(W)$ for $N\in \N$
of $U^{\infty}(W)$ spanned by elements of the form 
$[w]_{kl}$ for $w\in W$, $k, l=0, \dots, N$.
Let
$$\widetilde{A}^{N}(W)=\{\mathfrak{w}+\widetilde{Q}^{\infty}(W)\;|\;\mathfrak{w}\in U^{N}(W)\}.$$
Note that $\widetilde{A}^{N}(W)$ is
spanned by elements of the form $[w]_{kl}+\tilde{Q}^{\infty}(W)$ for $w\in W$ and $k, l=0, \dots, N$.
Also note that $\widetilde{A}^{N}(W)$ is an $\widetilde{A}^{N}(V)$-subbimodule of 
$\widetilde{A}^{\infty}(W)$ viewed as a $\widetilde{A}^{N}(V)$-bimodule.

\begin{prop}\label{U-1-isom-N}
The subspace $\widetilde{A}^{N}(W)$ of $\widetilde{A}^{\infty}(W)$
is in variant under the left and right actions of 
$\widetilde{A}^{N}(V)\subset \widetilde{A}^{\infty}(V)$.
In particular, $\widetilde{A}^{N}(W)$ is 
an $\widetilde{A}^{N}(V)$-bimodule. Moreover, 
$\mathcal{U}_{W}(1)$ induces an invertible 
linear map sending the $A^{N}(V)$-bimodule structure on $A^{N}(W)$
to the $\widetilde{A}^{N}(V)$-bimodule structure on 
$\widetilde{A}^{N}(W)$.
\end{prop}

Let $W_{1}$, $W_{2}$ and $W_{3}$ be lower-bounded generalized $V$-modules
and $\Y$ an intertwining operator of type $\binom{W_{3}}{W_{1}W_{2}}$. 
In \cite{H-aa-int-op}, a linear map $\vartheta_{\Y}:
U^{N}(W_{1})\to \hom(W_{2}, W_{3})$
is defined. Since $Q^{\infty}(V)\subset \ker \vartheta_{\Y}$,
we can view $\vartheta_{\Y}$ as a map from $A^{\infty}(V)$
to $\hom(W_{2}, W_{3})$. 
Now we define 
$\tilde{\vartheta}_{\Y}: U^{N}(W_{1})\to \hom(W_{2}, W_{3})$ 
by 
$$\tilde{\vartheta}_{\Y}(\mathfrak{w}_{1})w_{2}
=\vartheta_{\Y}(\mathcal{U}_{W_{1}}(1)\mathfrak{w}_{1})w_{2}$$
for $\mathfrak{w}_{1}\in U(W_{1})$ and $w_{2}\in W_{2}$. 
In \cite{H-aa-int-op}, for $N\in \N$, a linear map 
$$\rho^{N}: \mathcal{V}_{W_{1}W_{2}}^{W_{3}}\to 
\hom(A^{N}(W_{1})\otimes \Omega_{N}^{0}(W_{2}), 
\Omega_{N}^{0}(W_{3}))$$
is defined by
$$(\rho^{N}(\Y))((\mathfrak{w}_{1}
+Q^{\infty}(W_{1}))\otimes_{A^{N}(V)} w_{2})
=\vartheta_{\Y}(\mathfrak{w}_{1})w_{2}$$
for $\mathfrak{w}_{1}\in U^{N}(W_{1})$ and $w_{2}\in 
\Omega_{N}^{0}(W_{2})$.  
Here we define a linear map
\begin{align*}
\tilde{\rho}^{N}: \mathcal{V}_{W_{1}W_{2}}^{W_{3}}&\to 
\hom(\widetilde{A}^{N}(W_{1})\otimes \Omega_{N}^{0}(W_{2}), 
\Omega_{N}^{0}(W_{3}))\\
\Y&\mapsto \tilde{\rho}^{N}(\Y)
\end{align*}
by
\begin{align*}
&(\tilde{\rho}^{N}(\Y))(\mathfrak{w}_{1}+\widetilde{Q}^{\infty}(W_{1}))\otimes_{A^{N}(V)} w_{2})\nn
&\quad =\tilde{\vartheta}_{\Y}(\mathfrak{w}_{1})w_{2}\nn
&\quad =\vartheta_{\Y}(\mathcal{U}_{W_{1}}(1)\mathfrak{w}_{1})w_{2}\nn
&\quad =(\rho^{N}(\Y))((\mathcal{U}_{W_{1}}(1)\mathfrak{w}_{1}
+Q^{\infty}(W_{1}))\otimes_{A^{N}(V)} w_{2})
\end{align*}
for $\mathfrak{w}_{1}\in U^{N}(W_{1})$ and $w_{2}\in 
\Omega_{N}^{0}(W_{2})$. The definition of $\tilde{\rho}^{N}(\Y)$
can also be simply written as
\begin{equation}\label{tilde-rho-N}
\tilde{\rho}^{N}(\Y)=(\rho^{N}(\Y))\circ (\mathcal{U}_{W_{1}}(1)
\otimes 1_{\Omega_{N}^{0}(W_{2})}),
\end{equation}
where we use
the same notation $\mathcal{U}_{W_{1}}(1)$ to denote the 
map from  $\widetilde{A}^{N}(W_{1})$ to $A^{N}(W_{1})$
induced from the operator $\mathcal{U}_{W_{1}}(1)$ on 
$U^{N}(W_{1})$.

\begin{prop}
The linear map $\tilde{\rho}^{N}$ is in fact from $\mathcal{V}_{W_{1}W_{2}}^{W_{3}}$ to 
$\hom_{\widetilde{A}^{N}(V)}(\widetilde{A}^{N}(W_{1})\otimes_{\widetilde{A}^{N}(V)} 
\Omega_{N}^{0}(W_{2}), 
\Omega_{N}^{0}(W_{3}))$.
\end{prop}
\pf
Let $\Y$ be an intertwining operator of type $\binom{W_{3}}{W_{1}W_{2}}$. 
For $\mathfrak{v}\in U^{\infty}(V)$, $\mathfrak{w}\in U^{\infty}(W)$ and 
$w_{2}\in \Omega_{N}^{0}(W_{2})$,
\begin{align*}
&(\tilde{\rho}^{N}(\Y))((\mathfrak{v}\bkdiamd \mathfrak{w}_{1}+\widetilde{Q}^{\infty}(W_{1}))\otimes w_{2})\nn
&\quad =\tilde{\vartheta}^{N}_{\Y}(\mathfrak{v}\bkdiamd \mathfrak{w}_{1})w_{2}\nn
&\quad =\vartheta_{\Y}(\mathcal{U}_{W_{1}}(1)(\mathfrak{v}\bkdiamd \mathfrak{w}_{1}))w_{2}\nn
&\quad =\vartheta_{\Y}((\mathcal{U}_{V}(1)\mathfrak{v})\diamond 
(\mathcal{U}_{W_{1}}(1)\mathfrak{w}_{1}))w_{2}\nn
&\quad =\vartheta_{W_{3}}(\mathcal{U}_{V}(1)\mathfrak{v})
\vartheta_{\Y}(\mathcal{U}_{W_{1}}(1)\mathfrak{w}_{1})w_{2}\nn
&\quad =\tilde{\vartheta}_{W_{3}}(\mathfrak{v})
\tilde{\rho}^{N}_{\Y}((\mathfrak{w}_{1}+\widetilde{Q}^{\infty}(W_{1}))\otimes w_{2}).
\end{align*}
This shows that the image of $\tilde{\rho}^{N}$ is in
$\hom_{\widetilde{A}^{N}(V)}(\widetilde{A}^{N}(W_{1})\otimes 
\Omega_{N}^{0}(W_{2}), 
\Omega_{N}^{0}(W_{3}))$.

On the other hand, for $\mathfrak{v}\in U^{\infty}(V)$, $\mathfrak{w}\in U^{\infty}(W)$ and 
$w_{2}\in \Omega_{N}^{0}(W_{2})$, 
\begin{align*}
&(\tilde{\rho}^{N}(\Y))((\mathfrak{w}_{1}\bkdiamd \mathfrak{v}+\widetilde{Q}^{\infty}(W_{1}))\otimes w_{2})\nn
&\quad =\tilde{\vartheta}^{N}_{\Y}(\mathfrak{w}_{1}\bkdiamd \mathfrak{v})w_{2}\nn
&\quad =\vartheta_{\Y}(\mathcal{U}_{W_{1}}(1)(\mathfrak{w}_{1}\bkdiamd \mathfrak{v}))w_{2}\nn
&\quad =\vartheta_{\Y}((\mathcal{U}_{W_{1}}(1)\mathfrak{w}_{1})\diamond 
(\mathcal{U}_{V}(1)\mathfrak{v}))w_{2}\nn
&\quad =\vartheta_{\Y}(\mathcal{U}_{W_{1}}(1)\mathfrak{w}_{1})
\vartheta_{W_{2}}(\mathcal{U}_{V}(1)\mathfrak{v})w_{2}\nn
&\quad =\tilde{\rho}^{N}_{\Y}((\mathfrak{w}_{1}+\widetilde{Q}^{\infty}(W_{1}))\otimes 
\tilde{\vartheta}_{W_{2}}(\mathfrak{v})w_{2}).
\end{align*}
This shows that the image of $\tilde{\rho}^{N}$ is in
$\hom_{\widetilde{A}^{N}(V)}(\widetilde{A}^{N}(W_{1})\otimes_{\widetilde{A}^{N}(V)} 
\Omega_{N}^{0}(W_{2}), 
\Omega_{N}^{0}(W_{3}))$.
\epfv

In Proposition 6.3 in \cite{H-aa-int-op}, given a graded $A^{N}(V)$-module $M$,
a lower-bounded generalized $V$-module 
$S^{N}_{\rm voa}(\Omega_{N}^{0}(M)$ satisfying a universal property  is constructed.

\begin{thm}\label{main-N}
Let $V$ be a vertex operator algebra. 
Assume that $W_{2}$ and $W_{3}'$ are equivalent to 
 $S^{N}_{\rm voa}(\Omega_{N}^{0}(W_{2}))$ and $S^{N}_{\rm voa}(\Omega_{N}^{0}(W_{3}'))$,
respectively.
Then the linear map $\tilde{\rho}^{N}$ is a linear isomorphism.
\end{thm}
\pf
We know that $\mathcal{U}_{W_{1}}(1)$ from 
$\widetilde{A}^{N}(W_{1})$ to $A^{N}(W_{1})$ is an isomorphism.
By Theorem 6.7 in \cite{H-aa-int-op}, 
$\rho^{N}$ is a linear isomorphism.
Then by (\ref{tilde-rho-N}),
we see that $\tilde{\rho}^{N}$ is also a linear isomorphism. 
\epfv

Note that in Theorem \ref{main-N}, one condition is that 
$W_{3}'$ is equivalent to $S^{N}_{\rm voa}(\Omega_{N}^{0}(W_{3}'))$. But in applications,
for example in the proof of the modular invariance that we shall give later, what we have is often that $W_{3}$ is 
equivalent to $S^{N}_{\rm voa}(\Omega_{N}^{0}(W_{3}))$. 
In these cases, we need the following result to use Theorem \ref{main-N}:

\begin{prop}
Let $W$ be a lower-bounded generalized $V$-module of finite length equivalent to 
$S^{N}_{\rm voa}(\Omega_{N}^{0}(W))$. Assume that $W$ is of finite length such that 
the lowest weight vectors of the irreducible $V$-modules in the composition series 
are given by cosets containing elements of $\Omega_{N}^{0}(W)$. Then 
$W'$ is equivalent to $S^{N}_{\rm voa}(\Omega_{N}^{0}(W'))$. 
\end{prop}
\pf
By definition, $S^{N}_{\rm voa}(\Omega_{N}^{0}(W'))$ is generated by 
$\Omega_{N}^{0}(W')$. By the universal property of $S^{N}_{\rm voa}(\Omega_{N}^{0}(W'))$,
there is a unique $V$-module map $f: S^{N}_{\rm voa}(\Omega_{N}^{0}(W'))\to W'$ such that 
$f|_{\Omega_{N}^{0}(W')}=1_{\Omega_{N}^{0}(W')}$. Since $W$ is of finite length such that the lowest weight vectors of the irreducible $V$-modules in the composition series 
are given by cosets containing elements of $\Omega_{N}^{0}(W)$, 
the same is also true for $W'$.
Then $W'$ is also generated by 
$\Omega_{N}^{0}(W')$. So $f$ is surjective. 

We still need to prove that $f$ is injective. If $\ker f\ne 0$, then it is a nonzero generalized 
$V$-submodule of $S^{N}_{\rm voa}(\Omega_{N}^{0}(W'))$. Since $W'$ is of finite length,
$\Gamma(W')$ must be a finite set. Then since $S^{N}_{\rm voa}(\Omega_{N}^{0}(W'))$
is generated by $\Omega_{N}^{0}(W')$, $\Gamma(S^{N}_{\rm voa}(\Omega_{N}^{0}(W')))$ 
is a finite set. Since $\ker f$ is a generalized $V$-submodule of 
$S^{N}_{\rm voa}(\Omega_{N}^{0}(W'))$, $\Gamma(\ker f)$ is also a finite set. 
In particular, there exists a lowest weight vector $w_{0}'\ne 0$ of $\ker f$. Let $W_{0}$ be the 
generalized $V$-module generated by $w_{0}'$. 
Applying Zorn's lemma to the 
generalized $V$-submodules of $W_{0}$ not containing $w_{0}'$,
we know that $W_{0}$ has a maximal generalized $V$-submodule $M_{0}$ not containing 
$w_{0}'$. Thus the quotient $W_{0}/M_{0}$ 
is an irreducible lower-bounded generalized $V$-module with a lowest weight vector $w_{0}'+M_{0}$. 
Since $W'$ is also of finite length, by Property 6 in Proposition \ref{properties},
$\Gamma(W')$ is finite and for each $\mu\in \Gamma(W')$,
there exist $h^{\mu}\in \mu$ equal to a lowest weight of an irreducible $V$-module such that 
$$\Omega_{N}^{0}(W')=\coprod_{\mu\in \Gamma(W')} \coprod_{n=0}^{N}W^{*}_{[h^{\mu}+n]}.$$
Since $S^{N}_{\rm voa}(\Omega_{N}^{0}(W'))$ is generated by $\Omega_{N}^{0}(W')$,
the weight of every element of $S^{N}_{\rm voa}(\Omega_{N}^{0}(W'))$ is also 
congruent to $h^{\mu}$ modulo $\Z$ for some $\mu\in \Gamma(W')$. 
In particular, $\wt w_{0}'$ is congruent to $h^{\mu_{0}}$ modulo $\Z$ for some $\mu_{0}\in \Gamma(W')$.
But $h^{\mu_{0}}$ is also the weight of a lowest weight vector of an irreducible $V$-module. 
Since the difference between lowest weights of two irreducible $V$-modules is less than 
or equal to $N$, $\wt w_{0}'-h^{\mu_{0}}\le N$. Hence 
$$w_{0}'\in \coprod_{n=0}^{N}W^{*}_{[h^{\mu_{0}}+n]}\subset \Omega_{N}^{0}(W').$$

Since $w_{0}'\in \ker f$, we have $f(w_{0}')=0$. On the other hand, since $w_{0}'\in \Omega_{N}^{0}(W')$,
we have $f(w_{0}')=1_{\Omega_{N}^{0}(W')}(w_{0}')=w_{0}'$. This contradicts with $w_{0}'\ne 0$.
Thus $\ker f=0$ and $f$ is injective. 
\epfv

Finally, we give several results on 
$A^{N}(V)$, $A^{N}(W)$, $\widetilde{A}^{N}(W)$
 and suitable generalized or ordinary $V$-modules 
when $V$ has no nonzero elements of negative weights
and is $C_{2}$-cofinite.

A lower-bounded generalized $V$-module $W$ is said to be of finite length 
if there are generalized $V$-submodules
$W=W_{1}\supset \cdots \supset W_{l+1}=0$ such that $W_{i}/W_{i+1}$ for
$i=1, \dots, l$ are irreducible lower-bounded generalized $V$-modules. 
A generalized $V$-module $W$ is said to be quasi-finite dimensional if 
$\coprod_{\Re(n)\le N}W_{[n]}$ is finite dimensional for $N\in \Z$. 

Our first proposition is in fact a collection of results from \cite{H-cofiniteness}
and \cite{H-aa-va} together with a well known result derived 
here as a consequence of 
these results.

\begin{prop}\label{properties}
Let $V$ be a $C_{2}$-cofinite vertex operator algebra 
without nonzero elements of negative weights.  
Then we have the following properties:

\begin{enumerate}

\item For $N\in \N$, $A^{N}(V)$ is finite dimensional. 

\item Every lower-bounded generalized $V$-module $W$ generated by 
a finite-dimensional subspace of $\Omega_{N}^{0}(W)$ is 
quasi-finite dimensional. 

\item Every irreducible lower-bounded generalized $V$-module is an ordinary 
irreducible $V$-module. 

\item The set of equivalence classes of (ordinary) irreducible  $V$-modules is 
in bijection with the set of equivalence classes of irreducible nondegenerate graded $A^{N}(V)$-modules. 

\item The category of lower-bounded  generalized $V$-module of finite length,
the category of grading-restricted generalized $V$-modules 
and the category of quasi-finite-dimensional 
generalized $V$-modules are the same. 

\item There are only finitely many irreducible $V$-modules up to equivalence.

\end{enumerate}
\end{prop}
\pf
Property 1 is Theorem 4.6 in \cite{H-aa-va}. 
Property 2 is Proposition 3.8 in \cite{H-cofiniteness}
in the case that $V$ is $C_{2}$-cofinite. 
Property 3 is Theorem 5.8 in \cite{H-aa-va}.
Property 4 follows from 
Theorems 5.6  and 5.8 in \cite{H-aa-va}.

Property 5 is Proposition 4.3 in \cite{H-cofiniteness}. Note that 
though Proposition 4.3 in \cite{H-cofiniteness} assumes that $V$
satisfies in addition $V_{(0)}=\C\one$, the only paper 
quoted there that needs this condition in some results is \cite{ABD}. 
But the result needed in the proof of Proposition 4.3 in \cite{H-cofiniteness}
is only Proposition 5.2 in \cite{ABD}, which does not need this condition. 
So Proposition 4.3 in \cite{H-cofiniteness} or Property 5 is true 
for $C_{2}$-cofinite $V$ without nonzero elements of negative weights. 

Property 6 is well known and follows immediately from the finite dimension
property of Zhu algebra for $C_{2}$-cofinite vertex operator algebra 
(Proposition 3.6 in \cite{DLM2} in the case $g=1_{V}$) and 
the correspondence between irreducible modules for Zhu algebra
and irreducible ($\N$-gradable weak) $V$-modules (Theorem 2.2.2 in \cite{Z}). 
Here we derive this from Properties 1 and 4 above.
Since $A^{N}(V)$ is finite dimensional, there are only finitely many inequivalent 
irreducible 
$A^{N}(V)$-modules. In particular, there are only finitely many inequivalent 
irreducible nondegenerate graded
$A^{N}(V)$-modules. By Property 4, we obtain Property 6 that 
there are only finitely many inequivalent irreducible $V$-modules.
\epfv

The following result is about the relation between a lower-bounded generalized 
$V$-module $W$ and $A^{N}(V)$-module 
$\Omega_{N}^{0}(W)$. 

\begin{prop}\label{more-properties}
Let $V$ be a $C_{2}$-cofinite vertex operator algebra 
without nonzero elements of negative weights.  
Let $N$ be a nonnegative integer such that 
the differences of the real parts of 
the lowest weights of the finitely many irreducible $V$-modules 
(up to equivalence)
are less than or equal to $N$. Then we have: 
\begin{enumerate}

\item For a lower-bounded generalized $V$-module $W$ of finite length,
$\Gamma(W)$ is a finite set and for each $\mu\in \Gamma(W)$, 
there exists $h^{\mu}\in \mu$ equal to a lowest 
weight of an irreducible $V$-module
such that 
\begin{align}
W&=\coprod_{\mu\in \Gamma(W)}\coprod_{n\in \N}W_{[h^{\mu}+n]}, \label{properties-1}\\
\Omega_{N}^{0}(W)&=\coprod_{\mu\in \Gamma(W)} \coprod_{n=0}^{N}W_{[h^{\mu}+n]}.
\label{properties-2}
\end{align}

\item A lower-bounded generalized $V$-module $W$ of finite length is 
generated by 
$\Omega_{N}^{0}(W)$ and is equivalent to 
$S_{\rm voa}^{N}(\Omega_{N}^{0}(W))$. 
In particular, $W'$ is equivalent to $S_{\rm voa}^{N}(W')$.

\item Let $N'\in N+\N$ and 
$$M=\coprod_{n=0}^{N'}M_{[n]}$$
a finite-dimensional graded $A^{N'}(V)$-module. 
Let 
$$M^{N}=\coprod_{n=0}^{N}M_{[n]}.$$
Then $M^{N}$ is a graded $A^{N}(V)$-module and 
$M$ is equivalent to the graded $A^{N'}(V)$-module 
$\Omega_{N'}^{0}(S_{\rm voa}^{N}(M^{N}))$.

\end{enumerate}
\end{prop}
\pf 
Let $W=W_{1}\supset \cdots \supset W_{l+1}=0$ be a composition sereis 
of a lower-bounded generalized $V$-module $W$ of finite length. Let $w_{i}\in W_{i}$ 
for $i=1, \dots, l$ be homogeneous 
such that $w_{i}+W_{i+1}$ is a lowest weight vector of $W_{i}/W_{i+1}$.
Since as a graded vector space, $W$ is isomorphic to $\coprod_{i=1}^{l}W_{i}/W_{i+1}$,
$\Gamma(W)$ is the set of all congruence classes in $\C/\Z$ containing 
the weights of at least one $w_{i}$. For each element $\mu\in \Gamma(W)$, 
let $h^{\mu}$ be the smallest of all $\wt w_{i}\in \mu$. 
Then (\ref{properties-1}) and (\ref{properties-2}) hold. 

Let $h_{W}$ be the smallest of all $h^{\mu}$ for $\mu\in \Gamma(W)$. 
Then $h_{W}$ is a lowest weight of $W$ and we have 
$$\coprod_{\Re(h_{W})\le \Re(m)\le \Re(h_{W})+N}W_{[m]}
\subset \coprod_{\mu\in \Gamma(W)} \coprod_{n=0}^{N}W_{[h^{\mu}+n]}=\Omega_{N}^{0}(W).$$
By Proposition 5.11 in \cite{H-aa-va}, 
$$\coprod_{\Re(h_{W})\le \Re(m)\le \Re(h_{W})+N}W_{[m]}$$ 
generates $W$. Thus $\Omega_{N}^{0}(W)$ also generates $W$. 
By the universal property of $S_{\rm voa}^{N}(\Omega_{N}^{0}(W))$, 
there exists a unique $V$-module map $f: 
S_{\rm voa}^{N}(\Omega_{N}^{0}(W))\to W$ such that 
$f|_{\Omega_{N}^{0}(W)}=1_{\Omega_{N}^{0}(W)}$.
Since $W$ is generated by $\Omega_{N}^{0}(W)$, $f$ is 
surjective. We now show that $f$ is also injective, that is, 
the kernel $\ker f$ of $f$ is $0$. In fact, $W$ is equivalent to the quotient of 
$S_{\rm voa}^{N}(\Omega_{N}^{0}(W))$ by $\ker f$. 
Since  $f|_{\Omega_{N}^{0}(W)}=1_{\Omega_{N}^{0}(W)}$,
we have
\begin{equation}\label{ker-f}
\ker f\subset \coprod_{\mu\in \Gamma(W)} \coprod_{n\in N+\Z_{+}}
W_{[h^{\mu}+n]}.
\end{equation}
But $\ker f$ is a lower-bounded generalized $V$-submodule of $W$,
if $\ker f\ne 0$,
its lowest weight must be equal to the lowest weight of an irreducible 
$V$-module. But the real part of the lowest weight of the irreducible 
$V$-module must be less than or equal to the real part of 
the lowest weight of $W$ plus $N$. So a lowest weight vector 
of $\ker f$ must be in $\Omega_{N}^{0}(W)$, contradictory to 
(\ref{ker-f}). Since $f$ is both injective and surjective, 
$f$ is an equivalence of generalized $V$-modules. 

For a finite-dimensional graded $A^{N'}(V)$-module $M$,
it is clear that $M^{N}$ is a graded $A^{N}(V)$-module.
Since $M$ is finite dimensional, 
by Property 2 in Proposition \ref{properties}, 
$S_{\rm voa}^{N'}(M)$ is quasi-finite dimensional. By 
Property 5 in Proposition \ref{properties}, $S_{\rm voa}^{N'}(M)$ 
is of finite length. By Property 2 above, $S_{\rm voa}^{N'}(M)$ 
is generated by $\Omega_{N}^{0}(S_{\rm voa}^{N'}(M))=M^{N}$ 
and is equivalent to $S_{\rm voa}^{N}(M^{N})$.
Since $\Omega_{N'}^{0}(S_{\rm voa}^{N'}(M))=M$,
we see that $M$ is equivalent to 
$\Omega_{N'}^{0}(S_{\rm voa}^{N}(M^{N}))$.
\epfv

In \cite{H-aa-va}, it is proved that $A^{N}(V)$ is finite dimensional 
when $V$ has no nonzero elements of negative weights,
$V_{(0)}=\C\one$ and is $C_{2}$-cofinite. In fact, the condition 
$V_{(0)}=\C\one$ is not needed. Below we prove that 
$A^{N}(W)$ is finite dimensional without this condition for 
a grading-restricted generalized $V$-module $W$. 

\begin{thm}\label{finite-d-A-N-W}
Let $V$ be a $C_{2}$-cofinite vertex operator algebra 
without nonzero elements of negative weights
and  $W$ a grading-restricted generalized $V$-module. 
Then $A^{N}(W)$ is finite dimensional. 
\end{thm}
\pf
By Proposition 5.2 in \cite{ABD}, 
every irreducible $V$-module is $C_{2}$-cofinite. 
By Proposition \ref{properties},  we see that $W$ is of finite length.
Since every irreducible $V$-module is $C_{2}$-cofinite, $W$ as a 
generalized $V$-module of finite length must also be $C_{2}$-cofinite. 
If we assume in addition that $V_{(0)}=\C\one$, then 
by Theorem 11 in 
\cite{GN} (see Proposition 5.5 in \cite{AbN}), 
$V$ is also $C_{n}$-cofinite for $n\ge 2$. By Proposition 5.1
in \cite{AbN}, $W$ is $C_{n}$-cofinite for $n\ge 2$. 

On the other hand, since Lemma 2.4 in \cite{M}
gives a spanning set of $W$ without
the condition $V_{(0)}=\C\one$, the arguments in \cite{GN} and 
\cite{AbN} can also be used to show that $W$ is $C_{n}$-cofinite for $n\ge 2$
without this condition. For reader's convenience, here we give 
a direct proof of this fact observed by McRae using Lemma 2.4 in \cite{M}. 
In fact, since $W$ is of finite length, it is finitely generated. 
By Lemma 2.4 in \cite{M}, $W$ is spanned by elements of the form
\begin{equation}\label{spanning-set}
(Y_{W})_{i_{1}}(v_{1})\cdots (Y_{W})_{i_{k}}(v_{k})w_{j}
\end{equation}
for homogeneous $v_{1}, \dots, v_{k}$ in a finite set, homogeneous 
$w_{j}\in W$ for $j=1, \dots, l$ and 
$i_{1}, \dots, i_{k}\in \Z$ satisfying $i_{1}<\cdots <i_{k}$. 
Using  the lower-truncation property of the vertex operator map $Y_{W}$
and the fact that $v_{k}$ can change only in a finite set and 
there are only finitely many $w_{j}$, we see that there exists $m\in\Z$ such that 
for $i_{k}>m$, $(Y_{W})_{i_{k}}(v_{k})w_{j}=0$ and hence 
(\ref{spanning-set}) is $0$. Since in (\ref{spanning-set}),
we have $i_{1}<\cdots <i_{k}$,  there are only finitely many 
elements of the form (\ref{spanning-set}) satisfying $i_{1}>-n$.
Since elements of the form \ref{spanning-set}) span $W$, 
we see that $W$ is $C_{n}$-cofinite for $n\ge 2$. 

Take $n=k+l+2$ for $k, l=0, \dots, N$. 
Then $W$ is $C_{k+l+2}$-cofinite. 
By definition, $C_{k+l+2}(W)$ are spanned by elements of the 
form $(Y_{W})_{-k-l-2}(v)w$ for $v\in V$ and $w\in W$. 
Since $W$ is $C_{k+l+2}$-cofinite,
there exists a finite dimensional subspace $X_{k+ l}$ of $W$ such that 
$X_{k+l}+C_{k+l+2}(W)=W$. Let $U^{N}(X)$ be 
the subspace of $U^{N}(W)$ consisting 
matrices in $U^{N}(W)$ whose $(k, l)$-th entries are in $X_{k+l}$
for $k, l=0, \dots, N$. 
Since $X_{k+l}$
for $k, l=0, \dots, N$ are finite dimensional, 
$U^{N}(W)$ is also finite dimensional. We now prove 
$U^{N}(X)+(O^{\infty}(W)\cap U^{N}(W))=U^{N}(W)$. 
To prove this, we need only prove that every element of $U^{N}(W)$ 
of the form 
$[w]_{kl}$ for $w\in W$ and $0\le k, l\le N$, can be written as $[w]_{kl}=[w_{1}]_{kl}+[w_{2}]_{kl}$, where $w_{1}\in X_{k+l}$ and 
$w_{2}\in V$ such that $[w_{2}]_{kl}\in O^{\infty}(W)$. We shall denote the subspace of $W$ consisting of elements
$w$ such that $[w]_{kl}\in O^{\infty}(W)$ by $O^{\infty}_{kl}(W)$. 
Then what we need to prove is $W=X_{k+l}+O^{\infty}_{kl}(W)$.

Since $W=\coprod_{n\in \N}W_{\l n\r}$ is of finite length, 
$W_{\l n\r}$ for $n\in \N$ is finite dimensional. 
We can always take $X_{k+l}$ to be a subspace of $W$ containing $W_{\l 0\r}$. 
We use
induction on $p$ for $w\in W_{\l p\r}$. 
When $w\in W_{\l 0\r}$, $w$
can be written as $w=w+0$, where $v\in X_{k+l}$ 
and $0\in O^{\infty}_{kl}(V)$.

Assume that  when $v\in W_{\l p\r}$ for $p<q$,  $w=w_{1}+w_{2}$, where  $w_{1}\in X_{k+l}$ and 
$w_{2}\in O_{kl}^{\infty}(V)$. 
Then since $W$ is $C_{k+l+2}$-cofinite, for $w\in W_{\l q\r}$, there exists
homogeneous $w_{1}\in X_{k+l}$ and $v^{i}\in V$ and 
$w^{i}\in W$ for $i=1, 
\dots, m$ such that $v=u_{1}+\sum_{i=1}^{m}
(Y_{W})_{-k-l-2}(v^{i})w^{i}$.
Moreover, we can always find such $w_{1}$ and $v^{i}, w^{i}\in V$ for $i=1, 
\dots, m$ such that $w_{1}, (Y_{W})_{-k-l-2}(v^{i})w_{i}\in W_{\l q\r}$. 
Since 
$(Y_{W})_{n-k-l-2}(v^{i})w^{i}\in W_{\l q-n\r}$, where 
$q-n<q$
for $i=1, 
\dots, m$ and
$n\in \Z_{+}$, by induction assumption, 
$(Y_{W})_{n-k-l-2}(v^{i})w^{i}\in X_{k+l}+O^{\infty}_{kl}(V)$
for $i=1, \dots, m$ and $k\in \Z_{+}$.
Thus 
\begin{align*}
w&=w_{1}+\sum_{i=1}^{m}(Y_{W})_{-k-l-2}(v^{i})w^{i}\nn
&=w_{1}+\sum_{i=1}^{m}
\res_{x}x^{-k-l-2}(1+x)^{l}Y_{W}((1+x)^{L(0)}v^{i}, x)w^{i}\nn
& \quad-\sum_{i=1}^{m}\sum_{n\in \Z_{+}}\binom{\wt v^{i}+l}{n}
(Y_{W})_{n-k-l-2}(v^{i})w^{i}.
\end{align*}
By definition, 
$$[\res_{x}x^{-k-l-2}(1+x)^{l}Y_{W}((1+x)^{L(0)}v^{i}, x)w^{i}]_{kl}\in O^{\infty}(W).$$
Thus 
$$\res_{x}x^{-k-l-2}(1+x)^{l}Y_{W}((1+x)^{L(0)}v^{i}, x)w^{i}\in O^{\infty}_{kl}(W).$$
Thus we have $w=w_{1}+w_{2}$, where  $w_{1}\in X_{k+l}$ and 
$w_{2}\in O_{kl}^{\infty}(W)$. 
By induction principle, we have $W=X_{k+l}+O^{\infty}_{kl}(W)$. 

We now have proved $U^{N}(X)+(O^{\infty}(W)\cap U^{N}(W))=U^{N}(W)$. 
Since $O^{\infty}(W)\cap U^{N}(W)\subset Q^{\infty}(W)\cap U^{N}(W)$, we also have 
$U^{N}(X)+(Q^{\infty}(W)\cap U^{N}(W))=U^{N}(W)$. Since $U^{N}(X)$ is finite
dimensional, $A^{N}(W)$ is finite dimensional. 
\epfv

By Theorems \ref{U-1-isom-N} and \ref{finite-d-A-N-W}, we obtain
immediately the following result:

\begin{cor}\label{finite-d-til-A-N-W}
Let $V$ be a $C_{2}$-cofinite vertex operator algebra 
without nonzero elements of negative weights
 and  $W$ a grading-restricted generalized $V$-module. 
Then $\widetilde{A}^{N}(W)$ is finite dimensional. 
\end{cor}

\renewcommand{\theequation}{\thesection.\arabic{equation}}
\renewcommand{\thethm}{\thesection.\arabic{thm}}
\setcounter{equation}{0} \setcounter{thm}{0} 

\section{Symmetric linear functions 
on $\widetilde{A}^{N}(V)$-bimodules}

In this section,  we give two constructions of symmetric linear functions
on the $\widetilde{A}^{N}(V)$-bimodule $\widetilde{A}^{N}(W)$ for $N\in \N$ 
and a  grading-restricted generalized 
$V$-module $W$. We first give symmetric 
linear functions on $\widetilde{A}^{N}(W)$ using shifted
pseudo-$q$-traces 
of intertwining operators. Then we also construct 
symmetric linear functions on $\widetilde{A}^{N}(W)$ starting 
from linear maps from a lower-bounded generalized $V$-module
to $\C\{q\}[\log q]$ satisfying suitable conditions corresponding to 
the conditions for genus-one $1$-point conformal blocks (see Definition 
\ref{1-pt-cb} in the
next section).
This second construction is the main technically difficult part of the present paper. 
In fact, the first construction can also be obtained using the 
second construction. But we still give the first construction to show that 
it is much easier to obtain 
such linear symmetric functions on  $\widetilde{A}^{N}(W)$ 
using the properties of shifted pseudo-$q$-traces of intertwining operators
than using only the properties of genus-one $1$-point conformal blocks. 

Let $W$ be a grading-restricted generalized 
$V$-module, $\widetilde{W}$ a grading-restricted 
generalized $V$-$P$-bimodule
which is projective as a right $P$-module
 and $\Y$ an intertwining operator 
of type $\binom{\widetilde{W}}{W\widetilde{W}}$ and 
compatible with $P$. 

Recall from \cite{H-aa-int-op} and the preceding section that we can write
$$\widetilde{W}=\coprod_{\mu\in \Gamma(\widetilde{W})}
\coprod_{n\in \N}\widetilde{W}_{[h^{\mu}+n]}
=\coprod_{\mu\in \Gamma(\widetilde{W})}\widetilde{W}^{\mu}
=\coprod_{n\in \N}\widetilde{W}_{\l n\r},$$
where $\Gamma(\widetilde{W})\subset \C/\Z$,
$h^{\mu}\in \mu$ for $\mu\in \Gamma(\widetilde{W})$ and
$$\widetilde{W}^{\mu}=\coprod_{n\in \N}\widetilde{W}_{[h^{\mu}+n]},\;\;
\widetilde{W}_{\l n\r}=\coprod_{\mu\in \Gamma(\widetilde{W})}
\widetilde{W}_{[h^{\mu}+n]}.$$
Then from the definition of shifted pseudo-$q$-traces, for $w\in W$, 
we have 
\begin{align}\label{pseudo-tr-int}
&\tr_{\widetilde{W}}^{\phi}\mathcal{Y}(\mathcal{U}_{W}(q_{z})w, q_{z})
q_{\tau}^{L(0)-\frac{c}{24}}\nn
&\quad=\sum_{k=0}^{K}\frac{(2\pi i)^{k}}{k!}
\sum_{\mu\in \Gamma(\widetilde{W})}\sum_{n\in \N}
\phi_{\widetilde{W}_{[h^{\mu}+n]}}\pi_{h^{\mu}+n}
\Y(\mathcal{U}_{W}(q_{z})w, q_{z})
L_{\widetilde{W}}(0)_{N}^{k}\lbar_{\widetilde{W}_{[h^{\mu}+n]}}
\tau^{k}q_{\tau}^{h^{\mu}+n-\frac{c}{24}}.
\end{align}

Note that $[L_{\widetilde{W}}(0)_{N}, Y_{\widetilde{W}}(v, x)]=0$
for $v\in V$, that is, $L_{\widetilde{W}}(0)_{N}$ is in fact a 
$V$-module map from $W$ to itself. Assume that there exists $K\in \N$ such that
$L_{\widetilde{W}}(0)_{N}^{K+1}\tilde{w}=0$ for $\tilde{w}\in \widetilde{W}$.
This is always true if $\widetilde{W}$ is of finite length. 
By Proposition \ref{properties},
this is always true when $V$ has no nonzero elements of negative weights
and is $C_{2}$-cofinite. 
In this case, for each $k=0, \dots, K$,
$\Y^{k}=\Y\circ (1_{W}\otimes L_{\widetilde{W}}(0)_{N}^{k})$
is an intertwining operators of the type 
$\binom{\widetilde{W}}{W\widetilde{W}}$. Note that $\Y^{0}=\Y$.

Let $\mu \in \Gamma(\widetilde{W})$ and $k=0, \dots, K$. We define 
$$\psi_{\Y^{k}, \phi}^{\mu}([w]_{mn})=0$$
for $w\in W$ and $m, n\in \N$ such that $m\ne n$ and 
\begin{align*}
\psi_{\Y^{k}, \phi}^{\mu}([w]_{nn})
&=\phi_{\widetilde{W}_{[h^{\mu}+n]}}
\pi_{h^{\mu}+n}
\Y^{k}(\mathcal{U}_{W}(q_{z})w, q_{z})
\lbar_{\widetilde{W}_{[h^{\mu}+n]}}\nn
&=\phi_{\widetilde{W}_{[h^{\mu}+n]}}
\pi_{h^{\mu}+n}
\Y(\mathcal{U}_{W}(q_{z})w, q_{z})
L_{\widetilde{W}}(0)_{N}^{k}\lbar_{\widetilde{W}_{[h^{\mu}+n]}}
\end{align*}
for $w\in W$ and $n\in \N$.
We then define
$$\psi_{\Y^{k}, \phi}([w]_{mn})
=\sum_{\mu\in \Gamma(\widetilde{W})}
\psi_{\Y^{k}, \phi}^{\mu}([w]_{mn})$$
for $w\in W$ and $m, n\in \N$.
For $m, n\in \N$ such that $m\ne n$, 
$$\psi_{\Y^{k}, \phi}([w]_{mn})=0$$
and for $n\in \N$, 
\begin{align*}
\psi_{\Y^{k}, \phi}([w]_{nn})
&=\sum_{\mu\in \Gamma(\widetilde{W})}
\phi_{\widetilde{W}_{[h^{\mu}+n]}}
\pi_{h^{\mu}+n}
\Y^{k}(\mathcal{U}_{W}(q_{z})w, q_{z})
\lbar_{\widetilde{W}_{[h^{\mu}+n]}}\nn
&=\sum_{\mu\in \Gamma(\widetilde{W})}\phi_{\widetilde{W}_{[h^{\mu}+n]}}
\pi_{h^{\mu}+n}
\Y(\mathcal{U}_{W}(q_{z})w, q_{z})
L_{\widetilde{W}}(0)_{N}^{k}\lbar_{\widetilde{W}_{[h^{\mu}+n]}}.
\end{align*}

Using the $L(-1)$-derivative property and $L(-1)$-commutator formula for the 
intertwining operator $\Y$, we have
\begin{align*}
&\frac{\partial}{\partial z}\tr_{\widetilde{W}}^{\phi}
\mathcal{Y}(\mathcal{U}_{W}(q_{z})w, q_{z})
q_{\tau}^{L_{\widetilde{W}}(0)-\frac{c}{24}}\nn
&\quad =\tr^{\phi}_{\widetilde{W}}\mathcal{Y}((2\pi iL_{W}(0)+
2\pi iq_{z}L_{W}(-1))
\mathcal{U}_{W}(q_{z})w, 
q_{z})q_{\tau}^{L_{\widetilde{W}}(0)-\frac{c}{24}}\nn
&\quad =2\pi i\tr^{\phi}_{\widetilde{W}}[L_{\widetilde{W}}(0), 
\mathcal{Y}(\mathcal{U}_{W}(q_{z})w, 
q_{z})]q_{\tau}^{L_{\widetilde{W}}(0)-\frac{c}{24}}\nn
&\quad=0.
\end{align*}
So $\tr_{\widetilde{W}}^{\phi}
\mathcal{Y}(\mathcal{U}_{W}(q_{z})w, q_{z})
q_{\tau}^{L_{\widetilde{W}}(0)-\frac{c}{24}}$ is 
independent of $z$. Then the coefficients of this series in powers 
of $q_{\tau}$ and 
$\tau$ are also independent of $z$. In particular, 
for $k=0, \dots, K$, $\mu\in \Gamma(\widetilde{W})$
and $m, n\in \N$
$\psi_{\Y^{k}, \phi}^{\mu}([w]_{mn})$ and
$\psi_{\Y^{k}, \phi}([w]_{mn})$ are independent of $z$. 
Then we obtain linear functions $\psi_{\Y^{k}, \phi}^{\mu}$ for $\mu \in 
\Gamma(\widetilde{W})$ and $\psi_{\Y^{k}, \phi}$ on 
$U^{\infty}(W)$. 

Since $\tr_{\widetilde{W}}^{\phi}
\mathcal{Y}(\mathcal{U}_{W}(q_{z})w, q_{z})
q_{\tau}^{L_{\widetilde{W}}(0)-\frac{c}{24}}$
is independent of $z$, we have 
\begin{align}\label{pseudo-q-trace-residue}
&\tr_{\widetilde{W}}^{\phi}
\mathcal{Y}(\mathcal{U}_{W}(q_{z})w, q_{z})
q_{\tau}^{L_{\widetilde{W}}(0)-\frac{c}{24}}\nn
&\quad =\text{coeff}_{\log q_{z}}^{0}\res_{q_{z}}q_{z}^{-1}
\tr_{\widetilde{W}}^{\phi}
\mathcal{Y}(\mathcal{U}_{W}(q_{z})w, q_{z})
q_{\tau}^{L_{\widetilde{W}}(0)-\frac{c}{24}}\nn
&\quad =
\tr_{\widetilde{W}}^{\phi}
\text{coeff}_{\log x}^{0}\res_{x}x^{-1}\mathcal{Y}(\mathcal{U}_{W}(x)w, x)
q_{\tau}^{L_{\widetilde{W}}(0)-\frac{c}{24}},
\end{align}
where as in \cite{H-aa-int-op}, we use 
$\text{coeff}_{\log q_{z}}^{0}$ to denote the operation to take the 
constant term of a polynomial in $\log q_{z}$. 

Let $n, N\in \N$ such that $n\le N$.
Using (\ref{pseudo-q-trace-residue}), $\mathcal{U}_{W}(x)
=x^{L_{W}(0)}\mathcal{U}(1)$, the definition of 
$\tilde{\vartheta}_{\Y}$ and $\pi_{h^{\mu}+m}\tilde{\vartheta}_{\Y}
([w]_{nn})=0$ for $m\ne n$, we obtain 
\begin{align}\label{psi-vartheta}
\psi_{\Y, \phi}^{k}([w]_{nn})&
=\sum_{\mu\in \Gamma(\widetilde{W})}
\phi_{\widetilde{W}_{[h^{\mu}+n]}}
\pi_{h^{\mu}+n}\text{coeff}_{\log x}^{0}\res_{x}x^{-1}
\Y^{k}(\mathcal{U}_{W}(x)w, x)
\lbar_{\widetilde{W}_{[h^{\mu}+n]}}\nn
&=\sum_{\mu\in \Gamma(\widetilde{W})}
\phi_{\widetilde{W}_{[h^{\mu}+n]}}
\pi_{h^{\mu}+n}\tilde{\vartheta}_{\Y^{k}}([w]_{nn})
\lbar_{\widetilde{W}_{[h^{\mu}+n]}}\nn
&=\sum_{\mu\in \Gamma(\widetilde{W})}\sum_{m=0}^{N}
\phi_{\widetilde{W}_{[h^{\mu}+m]}}
\pi_{h^{\mu}+m}\tilde{\vartheta}_{\Y^{k}}([w]_{nn})
\lbar_{\widetilde{W}_{[h^{\mu}+m]}}\nn
&=\phi_{\Omega_{N}^{0}(\widetilde{W})}
\tilde{\vartheta}_{\Y^{k}}([w]_{nn})\lbar_{\Omega_{N}^{0}(\widetilde{W})}.
\end{align}
for $n\in \N$. 

\begin{lemma}
For $m, n\in \N$, $k=0, \dots, K$, $v\in V$ and $w\in W$, 
we have 
\begin{equation}\label{psi-Y-symmetry-0}
\psi_{\Y^{k}, \phi}([v]_{mn}\bkdiamd [w]_{nm})
=\psi_{\Y^{k}, \phi}([w]_{nm}\bkdiamd [v]_{mn}).
\end{equation}
\end{lemma}
\pf
By definition, 
$$[v]_{mn}\bkdiamd [w]_{nm}
=\res_{x}T_{2m+1}((x+1)^{-2m+n-1})(1+x)^{m}
\left[Y_{W}\left(v, \frac{1}{2\pi i}\log (1+x)\right)w\right]_{mm}.$$
Let $N\in \N$ be larger than or equal to both $m$ and $n$. 
Using (\ref{psi-vartheta}) and 
\begin{align*}
\tilde{\vartheta}_{\Y}([v]_{mn}\bkdiamd
[w]_{nm})&=\tilde{\vartheta}_{\widetilde{W}}([v]_{mn})
\tilde{\vartheta}_{\Y}([w]_{nm}),\nn
\tilde{\vartheta}_{\Y}([w]_{nm}\bkdiamd
[v]_{mn})&=\tilde{\vartheta}_{\Y}([w]_{nm})
\tilde{\vartheta}_{\widetilde{W}}([v]_{mn}),
\end{align*}
we have 
\begin{align*}
\psi_{\Y^{k}, \phi}&([v]_{mn}\bkdiamd [w]_{nm})\nn
&=\phi_{\Omega_{N}^{0}(\widetilde{W})}\tilde{\vartheta}_{\Y^{k}}
([v]_{mn}\bkdiamd [w]_{nm})\lbar_{\Omega_{N}^{0}(\widetilde{W})}\nn
&=\phi_{\Omega_{N}^{0}(\widetilde{W})}
\tilde{\vartheta}_{\widetilde{W}}([v]_{mn})
\tilde{\vartheta}_{\Y^{k}}([w]_{nm})\lbar_{\Omega_{N}^{0}(\widetilde{W})}\nn
&=\phi_{\Omega_{N}^{0}(\widetilde{W})}
\tilde{\vartheta}_{\Y^{k}}([w]_{nm})
\tilde{\vartheta}_{\widetilde{W}}([v]_{mn})\lbar_{\Omega_{N}^{0}(\widetilde{W})}\nn
&=\phi_{\Omega_{N}^{0}(\widetilde{W})}
\tilde{\vartheta}_{\Y^{k}}([w]_{nm}\bkdiamd [v]_{mn})
\lbar_{\Omega_{N}^{0}(\widetilde{W})}\nn
&=\psi_{\Y^{k}, \phi}([w]_{nm}\bkdiamd [v]_{mn}).
\end{align*}
\epfv

Fix $N\in \N$. Then the restriction of the linear function $\psi_{\Y^{k}, \phi}$ on 
$U^{\infty}(W)$ to $U^{N}(W)$ is a linear function $\psi_{\Y^{k}, \phi}^{N}$
on $U^{N}(W)$. 

\begin{prop}\label{slf-int-op}
The linear function $\psi_{\Y^{k}, \phi}^{N}$ on 
$U^{\infty}(W)$ is in fact a symmetric linear function on 
$\widetilde{A}^{N}(W)$, that is, 
$\psi_{\Y^{k}, \phi}^{N}(\widetilde{Q}^{\infty}(W)\cap U^{N}(W))=0$
and 
\begin{equation}\label{psi-Y-symmetry}
\psi_{\Y^{k}, \phi}^{N}(\mathfrak{v}\bkdiamd \mathfrak{w})
=\psi_{\Y^{k}, \phi}^{N}(\mathfrak{w}\bkdiamd \mathfrak{v})
\end{equation}
for $\mathfrak{v}\in U^{N}(V)$ and $\mathfrak{w}\in U^{N}(W)$.
\end{prop}
\pf
Since $\psi_{\Y^{k}, \phi}([w]_{mn})=0$ for $m\ne n$, to prove
$\psi_{\Y^{k}, \phi}^{N}(\widetilde{Q}^{\infty}(W)\cap U^{N}(W))=0$,
we need only prove $\psi_{\Y^{k}, \phi}([w]_{nn})=0$ if 
$[w]_{nn}\in \widetilde{Q}^{\infty}(W)$. 
From (\ref{psi-vartheta}), we have 
$$\psi_{\Y^{k}, \phi}([w]_{nn})=\phi_{\Omega_{N}^{0}(\widetilde{W})}
\tilde{\vartheta}_{\Y^{k}}([w]_{nn})\lbar_{\Omega_{N}^{0}(\widetilde{W})}.$$
But from the definition of $\widetilde{Q}^{\infty}(W)$,
$\tilde{\vartheta}_{\Y^{k}}([w]_{nn})=0$ 
if $[w]_{nn}\in \widetilde{Q}^{\infty}(W)$. This proves
$\psi_{\Y^{k}, \phi}^{N}(\widetilde{Q}^{\infty}(W)\cap U^{N}(W))=0$.

For $\mathfrak{v}=[v]_{mn}$ and $\mathfrak{w}=[w]_{kl}$,
if $n\ne k$, we have 
$$\mathfrak{v}\bkdiamd \mathfrak{w}=[v]_{mn}\bkdiamd [w]_{kl}=0$$
and 
$$\psi_{\Y^{k}, \phi}(\mathfrak{w}\bkdiamd \mathfrak{v})=
\psi_{\Y^{k}, \phi}([w]_{kl}\bkdiamd [v]_{mn})=0.$$
So (\ref{psi-Y-symmetry}) holds in this case. The same argument shows that 
if $m\ne l$, (\ref{psi-Y-symmetry}) holds.
The case $n=k$ and $m=l$ is given by (\ref{psi-Y-symmetry-0}).
Thus (\ref{psi-Y-symmetry}) holds for all 
$\mathfrak{v}\in U^{N}(V)$ and $\mathfrak{w}\in U^{N}(W)$.
\epfv

Proposition \ref{slf-int-op} gives the first construction of linear symmetric functions 
on $\widetilde{A}^{N}(W)$. We now give the second construction.

In the formulations, discussions and proofs below, for $m\in \Z_{+}$,
 $(e^{2\pi ix}-1)^{-m}$ is always understood as the 
formal Laurent series obtained by expanding $(e^{2\pi iz}-1)^{-m}$ 
as a Laurent series near $z=0$ and then replacing $z$ by $x$. 

The following theorem gives our second construction of 
symmetric linear functions on $\widetilde{A}^{N}(W)$ 
and is the main result of this section:

\begin{thm}\label{wp-fn-gr-sym-fn}
Let $V$ be a $C_2$-cofinite vertex operator algebra without nonzero elements
of negative weights. 
Let $W$ be a grading-restricted generalized $V$-module and 
\begin{align*}
S: W&\to \C\{q\}[\log q]\nn
w&\mapsto S(w; q)
\end{align*} 
a linear map satisfying  
\begin{align}
&S(w; q)=\sum_{k=0}^{K}\sum_{j=1}^{J}
\sum_{m\in \N}S_{k, j, m}(w)(\log q)^{k}q^{r_{j}+m},\nn
&S(\res_{x}Y_{W}(v, x)w; q)=0,\label{main-lemma-1}\\
&S(\res_{x}\tilde{\wp}_{2}(x; q)Y_{W}(v, x)w; q)=0,
\label{main-lemma-2}\\
&(2\pi i)^{2}q\frac{\partial}{\partial q}
S(w; q)=S(\res_{x}(\tilde{\wp}_{1}(x; q)-\widetilde{G}_{2}(q)x)Y_{W}
(\omega, x)w; q)\label{2nd-main-lemma-1}
\end{align}
for $u\in V$ and $w\in W$, where $K\in \N$, $r_{1}, \dots, r_{J}\in \C$ 
are independent of $w$. 
For $N\in \N$, $k=0, \dots, K$ and $j=1, \dots, J$, let 
$\psi^{N}_{S; k, j}: U^{N}(W)\to  \C$ be the linear map defined by 
$\psi_{S; k, j}^{N}([w]_{mn})=0$ for $w\in W$ and $m, n\in \N$ satisfying
$0\le m, n\le N$ and $m\ne n$
and 
$\psi_{S; k, j}^{N}([w]_{nn})=S_{k, j, n}(w)$ 
for $w\in W$ and 
$n\in \N$ satisfying $n\le N$.
Then $\psi^{N}_{S; k, j}(\widetilde{Q}^{\infty}(W))=0$ 
so that $\psi^{N}_{S; k, j}$ for $k=0, \dots, K$ and $j=1, \dots, J$
induce linear maps from $\widetilde{A}^{N}(W)$ to $\C$,
still denoted by $\psi^{N}_{S; k, j}$. Moreover, these induced 
linear maps $\psi^{N}_{S; k, j}$
are in fact symmetric linear functions on $\widetilde{A}^{N}(W)$, 
that is, for $\mathfrak{v}\in \widetilde{A}^{N}(V)$ and 
$\mathfrak{w}\in \widetilde{A}^{N}(W)$, 
$$\psi^{N}_{S; k, j}(\mathfrak{v}\bkdiamd \mathfrak{w})
=\psi^{N}_{S; k, j}(\mathfrak{w}\bkdiamd \mathfrak{v})$$
satisfying 
\begin{equation}\label{grading-sym-linear-fn-1}
\psi_{S; k, j}^{N}(([\omega]_{nn}
-(r_{j}+n)[\one]_{nn})^{\bkdiamd (K-k+1)}
\bkdiamd \mathfrak{w})=0
\end{equation}
for $n=0, \dots, N$ and $\mathfrak{w}\in U(W)$, 
where 
$$([\omega]_{nn}
-(r_{j}+n)[\one]_{nn})^{\bkdiamd (K-k+1)}
=\overbrace{([\omega]_{nn}
-(r_{j}+n)[\one]_{nn})\bkdiamd \cdots 
\bkdiamd ([\omega]_{nn}
-(r_{j}+n)[\one]_{nn})}^{K-k+1}.$$
\end{thm}

We will prove this result after we prove a number of lemmas and propositions.
We first prove a lemma for the vertex operators for a generalized $V$-module. 

\begin{lemma}
Let $W$ be a generalized $V$-module. 
For $m, n\in \N$ satisfying $n\ge m$, we have 
\begin{align}\label{lemma-1-1}
&\res_{x}e^{2(m+1)\pi ix}
(e^{2\pi ix}-1)^{-n-2}Y_{W_{1}}(v, \pm x)w\nn
&\quad =\sum_{j=0}^{m}\frac{1}{(n-j+1)}\binom{m}{j}
\res_{x}
e^{2\pi ix}(e^{2\pi ix}-1)^{-2}
Y_{W}\left( \binom{\pm\frac{1}{2\pi i}L_{V}(-1)-1}{n-j}
v, \pm x\right)w.
\end{align}
\end{lemma}
\pf
We first prove (\ref {lemma-1-1}) in the special case $m=0$, that is,
\begin{align}\label{lemma-1-1.5}
&\res_{x}e^{2\pi ix}
(e^{2\pi ix}-1)^{-n-2}Y_{W_{1}}(v, \pm x)w\nn
&\quad =\frac{1}{n+1}\res_{x}
e^{2\pi ix}(e^{2\pi ix}-1)^{-2}
Y_{W}\left(\binom{\pm \frac{1}{2\pi i}L_{V}(-1)-1}{n}
v, \pm x\right)w
\end{align}
for $n\in \N$. In this case,
\begin{align}\label{recurrence}
&\res_{x}e^{2\pi ix}(e^{2\pi ix}-1)^{-n-2}
Y_{W}(v, \pm x)w\nn
&\quad =\res_{x}e^{4\pi ix}(e^{2\pi ix}-1)^{-n-2}
Y_{W}(v, \pm x)w-\res_{x}e^{2\pi ix}(e^{2\pi ix}-1)^{-n-1}
Y_{W}(v, \pm x)w\nn
&\quad = -\frac{1}{2\pi i(n+1)}\res_{x}
\left(\frac{d}{dx}e^{2\pi ix}(e^{2\pi ix}-1)^{-n-1}\right)
Y_{W}(v, \pm x)w\nn
&\quad \quad +\frac{1}{n+1}\res_{x}
e^{2\pi ix}(e^{2\pi ix}-1)^{-n-1}
Y_{W}(v, \pm x)w\nn
&\quad \quad -\res_{x}e^{2\pi ix}(e^{2\pi ix}-1)^{-n-1}
Y_{W}(v, \pm x)w\nn
&\quad =\frac{1}{2\pi i(n+1)}\res_{x}
e^{2\pi ix}(e^{2\pi ix}-1)^{-n-1}
\frac{d}{dx}Y_{W}(v, \pm x)w\nn
&\quad \quad -\frac{n}{n+1}\res_{x}e^{2\pi ix}(e^{2\pi ix}-1)^{-n-1}
Y_{W}(v, \pm x)w\nn
&\quad =\frac{1}{(n+1)}\res_{x}
e^{2\pi ix}(e^{2\pi ix}-1)^{-n-1}
Y_{W}\left(\left(\pm \frac{1}{2\pi i}L_{V}(-1)-n\right)v, \pm x\right)w.
\end{align}
Using (\ref{recurrence}) repeatedly, we have
\begin{align*}
&\res_{x}e^{2\pi ix}(e^{2\pi ix}-1)^{-n-2}Y_{W}(v, \pm x)w\nn
&\quad =\frac{1}{(n+1)!}\res_{x}
e^{2\pi ix}(e^{2\pi ix}-1)^{-2}\cdot\nn
&\quad\quad\quad\quad\quad\cdot 
Y_{W}\left(\left(\pm \frac{1}{2\pi i}L_{V}(-1)-1\right)
\cdots \left(\pm \frac{1}{2\pi i}L_{V}(-1)-n\right)v, \pm x\right)w\nn
&\quad=\frac{1}{(n+1)}\res_{x}
e^{2\pi ix}(e^{2\pi ix}-1)^{-2}
Y_{W}\left(\binom{\pm \frac{1}{2\pi i}L_{V}(-1)-1}{n}
v, \pm x\right)w.
\end{align*}

In the general case,
using the binomial expansion of $(1+(e^{2\pi x}-1))^{m}$,
we have 
\begin{align}\label{main-lemma-4.5}
&\res_{x}e^{2(m+1)\pi ix}(e^{2\pi ix}-1)^{-n-2}Y_{W}(v, \pm x)w\nn
&\quad =\sum_{j=0}^{m}\binom{m}{j}\res_{x}
e^{2\pi ix}(e^{2\pi ix}-1)^{-(n-j)-2}Y_{W}(v, \pm x)w.
\end{align}
Since $n\ge m$, we have $n\ge j$ for $j=0, \dots, m$.
Then (\ref{lemma-1-1.5}) gives
\begin{align}\label{main-lemma-5}
&\res_{x}e^{2\pi ix}(e^{2\pi ix}-1)^{-(n-j)-2}Y_{W}(v, \pm x)w\nn
&\quad=\frac{1}{(n-j+1)}\res_{x}
e^{2\pi ix}(e^{2\pi ix}-1)^{-2}
Y_{W}\left(\binom{\pm \frac{1}{2\pi i}L_{V}(-1)-1}{n-j}
v, \pm x\right)w.
\end{align}
Using (\ref{main-lemma-4.5}) and \ref{main-lemma-5}), 
we obtain (\ref{lemma-1-1}). 
\epfv

We now prove a lemma giving some identities 
for maps from a generalized $V$-module to 
the space of power series in a formal variable $q$ satisfying 
(\ref{main-lemma-1}) and (\ref{main-lemma-2}).

\begin{lemma}
Let $W$ be a generalized $V$-module. 
Assume that a linear map
\begin{align*}
S: W&\to \C[[q]]\nn
w&\mapsto S(w; q)=\sum_{m\in \N}S_{m}(w)q^{m}
\end{align*} 
satisfies (\ref{main-lemma-1}) and (\ref{main-lemma-2}). 
Then for $m, n\in \N$ satisfying $n\ge m$, 
\begin{equation}\label{mod-inv-10.3}
S_{0}(\res_{x}e^{2(m+1)\pi ix}(e^{2\pi ix}-1)^{-n-2}
Y_{W_{1}}(u, \pm x)w_{1})=0
\end{equation}
and 
\begin{align}\label{mod-inv-10.4}
&S_{p}(\res_{x}e^{2(m+1)\pi ix}
(e^{2\pi ix}-1)^{-n-2}Y_{W_{1}}(v, \pm x)w_{1})\nn
&\quad =-\sum_{j=0}^{m}\sum_{s=1}^{p}\sum_{l | s}
\frac{l}{(n-j+1)}\binom{m}{j}\cdot\nn
&\quad\quad\quad\quad\quad\cdot 
 S_{p-s}\left(\res_{x}\left(\binom{- l-1}{n-j}
e^{2l\pi ix}+\binom{ l-1}{n-j}e^{-2l\pi ix}\right)
Y_{W}(v, \pm x)w\right).
\end{align}
for $p\in \Z_{+}$. 
\end{lemma}
\pf 
We first prove (\ref{mod-inv-10.3})
and (\ref{mod-inv-10.4}) in the special case $m=n=0$, that is, 
\begin{equation}\label{mod-inv-10.3-1}
S_{0}(\res_{x}e^{2\pi ix}(e^{2\pi ix}-1)^{-2}
Y_{W_{1}}(u, \pm x)w_{1})=0
\end{equation}
and 
\begin{equation}\label{mod-inv-10.4-1}
S_{p}(\res_{x}e^{2\pi ix}
(e^{2\pi ix}-1)^{-2}Y_{W_{1}}(v, \pm x)w_{1})
=-\sum_{s=1}^{p}\sum_{l | s}
l S_{p-s}(\res_{x}(
e^{2l\pi ix}+e^{-2l\pi ix}))
Y_{W}(v, \pm x)w)
\end{equation}
for $p\in \Z_{+}$
From (\ref{tilde-wp-2}) 
and (\ref{main-lemma-2}), 
we have
\begin{align}\label{lemma-2-1}
0&=S(\res_{x}\tilde{\wp}_{2}(x; q)Y(u, \pm x)w_{1}; q)\nn
&=(2\pi i)^{2}\sum_{n\in \mathbb{N}}
S_{n}(\res_{x}e^{2\pi ix}(e^{2\pi ix}-1)^{-2}
Y_{W_{1}}(u, \pm x)w_{1})
q^{n}\nn
&\quad +(2\pi i)^{2}\sum_{n\in \mathbb{N}}\sum_{s\in \Z_{+}}\sum_{l | s}l
S_{n}(\res_{x}
(e^{2l\pi ix}+e^{-2l\pi ix})Y_{W_{1}}(u, \pm x)w_{1})
q^{n+s}\nn
&\quad -\frac{\pi^{2}}{3}\sum_{n\in \mathbb{N}}
S_{n}(\res_{x}Y_{W_{1}}(u, \pm x)w_{1})
q^{n}\nn
&\quad -2(2\pi i)^{2}\sum_{n\in \mathbb{N}}\sum_{s\in \Z_{+}}
\sigma(s)
S_{n}(\res_{x}Y_{W_{1}}(u, \pm x)w_{1})
q^{n+s}.
\end{align}
Using (\ref{main-lemma-1}), we see that (\ref{lemma-2-1})
gives
\begin{align}\label{mod-inv-10.2}
&\sum_{n\in \mathbb{N}}
S_{n}(\res_{x}e^{2\pi ix}(e^{2\pi ix}-1)^{-2}
Y_{W_{1}}(u, \pm x)w_{1})
q^{n}\nn
&\quad =-\sum_{n\in \mathbb{N}}\sum_{s\in \Z_{+}}\sum_{l | s}l
S_{n}(\res_{x}
(e^{2l\pi ix}+e^{-2l\pi ix})Y_{W_{1}}(u, \pm x)w_{1})
q^{n+s}.
\end{align}
Taking the coefficients of the power of $q^{p}$ for $p\in \N$
of both sides of 
(\ref{mod-inv-10.2}), we obtain (\ref{mod-inv-10.3-1})
and (\ref{mod-inv-10.4-1}) for $p\in \Z_{+}$.

For $m, n\in \N$ satisfying $n\ge m$, using (\ref{lemma-1-1})
and (\ref{mod-inv-10.3-1}), we obtain
\begin{align*}
&S_{0}(\res_{x}e^{2(m+1)\pi ix}(e^{2\pi ix}-1)^{-n-2}
Y_{W_{1}}(u, \pm x)w_{1})\nn
&\quad =\sum_{j=0}^{m}\frac{1}{(n-j+1)}\binom{m}{j}
S_{0}\left(\res_{x}
e^{2\pi ix}(e^{2\pi ix}-1)^{-2}
Y_{W}\left(\binom{\pm \frac{1}{2\pi i}L_{V}(-1)-1}{n-j}
v, \pm x\right)w\right)\nn
&\quad =0,
\end{align*}
proving (\ref{mod-inv-10.3}). 
For $m, n\in \N$ satisfying $n\ge m$ and $p\in \Z_{+}$, 
using (\ref{lemma-1-1}) and 
(\ref{mod-inv-10.4-1}), we obtain
\begin{align}\label{main-lemma-7-1}
&S_{p}(\res_{x}e^{2(m+1)\pi ix}
(e^{2\pi ix}-1)^{-n-2}Y_{W_{1}}(v, \pm x)w_{1})\nn
&\quad =-\sum_{j=0}^{m}\sum_{s=1}^{p}\sum_{l | s}
\frac{l}{(n-j+1)}\binom{m}{j}\cdot\nn
&\quad\quad\quad\quad\quad\cdot 
 S_{p-s}\left(\res_{x}
(e^{2l\pi ix}+e^{-2l\pi ix})
Y_{W}\left(\binom{\pm \frac{1}{2\pi i}L_{V}(-1)-1}{n-j}
v, \pm x\right)w\right)\nn
&\quad =-\sum_{j=0}^{m}\sum_{s=1}^{p}\sum_{l | s}
\frac{l}{(n-j+1)}\binom{m}{j}\cdot\nn
&\quad\quad\quad\quad\quad\cdot 
 S_{p-s}\left(\res_{x}
(e^{2l\pi ix}+e^{-2l\pi ix})
\binom{\frac{1}{2\pi i}\frac{d}{dx}-1}{n-j}Y_{W}(
v, \pm x)w\right)\nn
&\quad =-\sum_{j=0}^{m}\sum_{s=1}^{p}\sum_{l | s}
\frac{l}{(n-j+1)}\binom{m}{j}\cdot\nn
&\quad\quad\quad\quad\quad\cdot 
 S_{p-s}\left(\res_{x}\left(\binom{-\frac{1}{2\pi i}\frac{d}{dx}-1}{n-j}
(e^{2l\pi ix}+e^{-2l\pi ix})\right)
Y_{W}(v, \pm x)w\right)\nn
&\quad =-\sum_{j=0}^{m}\sum_{s=1}^{p}\sum_{l | s}
\frac{l}{(n-j+1)}\binom{m}{j}\cdot\nn
&\quad\quad\quad\quad\quad\cdot 
 S_{p-s}\left(\res_{x}\left(\binom{-l-1}{n-j}
e^{2l\pi ix}+\binom{l-1}{n-j}e^{-2l\pi ix}\right)
Y_{W}(
v, \pm x)w\right),
\end{align}
proving (\ref{mod-inv-10.4}). 
\epfv

We are now ready to prove our main technical result.

\begin{prop}\label{main-lemma}
Let $W$ be a grading-restricted generalized $V$-module. 
Assume that a linear map
\begin{align*}
S: W&\to \C[[q]]\nn
w&\mapsto S(w; q)=\sum_{m\in \N}S_{m}(w)q^{m}
\end{align*} 
satisfies (\ref{main-lemma-1}) and (\ref{main-lemma-2})
for 
$u\in \V$ and $w\in W$. Then we have:

\begin{enumerate}

\item For $m, n, p\in \N$ satisfying 
$0\le p\le m$ and $2m\le n$, $v\in V$ and $w\in W$, 
\begin{equation}\label{main-lemma-3}
S_{p}(\res_{x}e^{2(m+1)\pi ix}(e^{2\pi ix}-1)^{-n-2}Y_{W}(v, x)w)=0.
\end{equation}

\item For $m, n\in \N$, $v\in V$ and $w\in W$, 
\begin{align}\label{main-lemma-4}
&\sum_{k=0}^{n}\binom{-2m+n-1}{k}
S_{m}(\res_{x}e^{2\pi i(m+1)x}(e^{2\pi ix}-1)^{-2m+n-k-1}
Y_{W}(v, x)w)\nn
&\quad =\sum_{k=0}^{m}\binom{-2n+m-1}{k}
S_{n}(\res_{x}e^{2\pi i(n+1)x}(e^{2\pi ix}-1)^{-2n+m-k-1}
Y_{W}(v, -x)w).
\end{align}
\end{enumerate}
\end{prop}
\pf 
Note that (\ref{main-lemma-3})  in the case $m=0$
is (\ref{mod-inv-10.3}) in the case $m=0$ with $\pm$ being $+$. 
So we need only prove 
(\ref{main-lemma-3}) in the case  $m, n\in \Z_{+}$.
From (\ref{mod-inv-10.4}), 
we have 
\begin{align*}
&S_{p}(\res_{x}e^{2(m+1)\pi ix}
(e^{2\pi ix}-1)^{-n-2}Y_{W_{1}}(v, x)w_{1})\nn
&\quad =-\sum_{j=0}^{m}\sum_{s=1}^{p}\sum_{l | s}
\frac{l}{(n-j+1)}\binom{m}{j}\cdot\nn
&\quad\quad\quad\quad\quad\cdot 
 S_{p-s}\left(\res_{x}\left(\binom{-l-1}{n-j}
e^{2l\pi ix}+\binom{l-1}{n-j}e^{-2l\pi ix})\right)
Y_{W}(
v, x)w\right)\nn
&\quad =-\sum_{s=1}^{p}\sum_{l | s}\sum_{j=0}^{m}
\frac{l}{(n-j+1)}\binom{m}{j}\binom{-l-1}{n-j}
 S_{p-s}\left(\res_{x}
e^{2l\pi ix}
Y_{W}(v, x)w\right)\nn
&\quad=0,
\end{align*}
where in the last two steps, we have used 
$\binom{l-1}{n-j}=0$ (since $0\le l-1<s\le p\le m\le 2m-j\le n-j$) and 
(\ref{lemma-3-1}). This proves
(\ref{main-lemma-3}).

We divide the proof of (\ref{main-lemma-4}) into three cases: 
$m=n$, $m>n$ and $m<n$. 

We first prove (\ref{main-lemma-4})  in the case $m=n$. 
In this case, the difference between
left-hand side and the right-hand side of (\ref{main-lemma-4})
is 
\begin{align}\label{main-lemma-20}
&\sum_{k=0}^{n}\binom{-n-1}{k}
S_{n}(\res_{x} e^{2\pi i(n+1)x}(e^{2\pi ix}-1)^{-n-k-1}Y_{W}(v, x)w)\nn
&\quad\quad -\sum_{k=0}^{n}\binom{-n-1}{k}
S_{n}(\res_{x} e^{2\pi i(n+1)x}(e^{2\pi ix}-1)^{-n-k-1}Y_{W}(v, -x)w)\nn
&\quad =\sum_{k=0}^{n}\binom{-n-1}{k}
S_{n}(\res_{x} e^{2\pi i(n+1)x}(e^{2\pi ix}-1)^{-n-k-1}Y_{W}(v, x)w)\nn
&\quad\quad +\sum_{k=0}^{n}\binom{-n-1}{k}
S_{n}(\res_{x} e^{-2\pi i(n+1)x}(e^{-2\pi ix}-1)^{-n-k-1}Y_{W}(v, x)w)\nn
&\quad =\sum_{k=0}^{n}\binom{-n-1}{k}
S_{n}(\res_{x} (e^{2\pi i(n+1)x}+(-1)^{-n-k-1}e^{2\pi ikx})
(e^{2\pi ix}-1)^{-n-k-1}Y_{W}(v, x)w).
\end{align}
Using (\ref{main-lemma-21}) and (\ref{main-lemma-1}),
we see that the right-hand side of (\ref{main-lemma-20})
is equal to $0$, proving (\ref{main-lemma-4}) in this case.

Swapping $m$ and $n$ in
(\ref{main-lemma-4}) in the case of $m<n$, we obtain 
\begin{align*}
&\sum_{k=0}^{n}\binom{-2m+n-1}{k}
S_{m}(\res_{x}e^{2\pi i(m+1)x}(e^{2\pi ix}-1)^{-2m+n-k-1}
Y_{W}(v, -x)w)\nn
&\quad =\sum_{k=0}^{m}\binom{-2n+m-1}{k}
S_{n}(\res_{x}e^{2\pi i(n+1)x}(e^{2\pi ix}-1)^{-2n+m-k-1}
Y_{W}(v, x)w),
\end{align*}
which differs from (\ref{main-lemma-4}) only by the sign 
of the variable x in the vertex operators. Thus we can prove 
(\ref{main-lemma-4}) in the cases $m>n$ and $m<n$ together 
by proving 
\begin{align}\label{main-lemma-4.1}
&\sum_{k=0}^{n}\binom{-2m+n-1}{k}
S_{m}(\res_{x}e^{2\pi i(m+1)x}(e^{2\pi ix}-1)^{-2m+n-k-1}
Y_{W}(v, \pm x)w)\nn
&\quad =\sum_{k=0}^{m}\binom{-2n+m-1}{k}
S_{n}(\res_{x}e^{2\pi i(n+1)x}(e^{2\pi ix}-1)^{-2n+m-k-1}
Y_{W}(v, \mp x)w).
\end{align}
in the case of $m>n$. 

For $m, n, k\in \N$ satisfying $m> n$,  $v\in V$ and $w\in W$,  
using (\ref{mod-inv-10.4}) and 
(\ref{lemma-3-1}) with $n$ there replaced by $2m-n+k-1$
and noting that $2m-n+k-1\ge m$, 
we have
\begin{align}\label{main-lemma-10}
&\sum_{k=0}^{n}\binom{-2m+n-1}{k}
S_{m}(\res_{x}e^{2\pi i(m+1)x}(e^{2\pi ix}-1)^{-2m+n-k-1}
Y_{W}(v, \pm x)w)\nn
&\quad=-\sum_{k=0}^{n}\sum_{j=0}^{m}
\sum_{s=1}^{m}\sum_{l | s}\binom{-2m+n-1}{k}
\frac{l}{2m-j-n+k}\binom{m}{j}\cdot\nn
&\quad\quad\quad\quad\quad\quad\cdot 
S_{m-s}\Biggl(\res_{x}
\Biggl(\binom{-l-1}{2m-j-n+k-1}e^{2l\pi ix}\nn
&\quad\quad\quad\quad\quad\quad\quad\quad\quad\quad\quad\quad\quad
+\binom{l-1}{2m-j-n+k-1}e^{-2l\pi ix}\Biggr)
Y_{W}\left(v, \pm x\right)w\Biggr)\nn
&\quad=-\sum_{k=0}^{n}\sum_{j=0}^{m}
\sum_{s=1}^{m}\sum_{l | s}\binom{-2m+n-1}{k}
\frac{l}{2m-j-n+k}\binom{m}{j}\binom{l-1}{2m-j-n+k-1}\cdot\nn
&\quad\quad\quad\quad\quad\quad\cdot 
S_{m-s}(\res_{x}
e^{-2l\pi ix}
Y_{W}(v, \pm x)w).
\end{align}

From  (\ref{main-lemma-13}), we see that in this case ($m>n$), 
the right-hand side of (\ref{main-lemma-10}) and thus also the 
left-hand side of (\ref{main-lemma-4.1}) are equal to 
\begin{equation}\label{main-lemma-14}
-\sum_{(m-n) | s, \; 1\le s\le m}
S_{m-s}(\res_{x}e^{-2(m-n)\pi ix}Y_{W}(v, \pm x)w).
\end{equation}

We now prove (\ref{main-lemma-4})
in the case $m>n$.
We give the proof in the three cases $m>2n$, $m=2n>0$ 
and $2n>m>n$ separately.

In the case $m> 2n$, we have $p(m-n)\ge 
2(m-n)=m+(m-2n)>m$ for $p\in \Z_{+}+1$. So the only 
integer $s$ satisfying $1\le s\le m$ and $(m-n)| s$ is 
$m-n$. Thus in this case, 
(\ref{main-lemma-14}) is equal to 
\begin{equation}\label{main-lemma-14.5}
-S_{n}(\res_{x}e^{-2(m-n)\pi ix}Y_{W}(v, \pm x)w).
\end{equation}
On the other hand, in this case ($m> 2n$), 
we have $m\ge -2n+m-1\ge 0$. Then
for $v\in V$ and $w\in W$,  we have 
\begin{align}\label{main-lemma-15}
&\sum_{k=0}^{m}\binom{-2n+m-1}{k}
S_{n}(\res_{x}e^{2\pi i(n+1)x}(e^{2\pi ix}-1)^{-2n+m-k-1}
Y_{W}(v, \mp x)w)\nn
&\quad =\sum_{k=0}^{-2n+m-1}\binom{-2n+m-1}{k}
S_{n}(\res_{x}e^{2\pi i(n+1)x}(e^{2\pi ix}-1)^{-2n+m-k-1}
Y_{W}(v, \mp x)w)\nn
&\quad =S_{n}(\res_{x}e^{2\pi i(m-n)x}
Y_{W}(v, \mp x)w)\nn
&\quad=-S_{n}(\res_{x}e^{-2\pi i(m-n)x}
Y_{W}(v, \pm x)w).
\end{align}
 So (\ref{main-lemma-4.1}) holds when 
$m>2n$. 

In the case $m=2n>0$, we have $m-n=n>0$. Then the only integers 
$s$ between $1$ and $m=2n$ containing a factor $n$ are 
$n$ and $2n$. So in this case, (\ref{main-lemma-14}) becomes
\begin{equation}\label{main-lemma-16}
-S_{n}(\res_{x}e^{-2n\pi ix}Y_{W}(v, \pm x)w)
-S_{0}(\res_{x}e^{-2n\pi ix}Y_{W}(v, \pm x)w).
\end{equation}
In this case, the right-hand side of (\ref{main-lemma-4.1})
is equal to 
\begin{align}\label{main-lemma-16.1}
&\sum_{k=0}^{2n}(-1)^{k}
S_{n}(\res_{x}e^{2\pi i(n+1)x}(e^{2\pi ix}-1)^{-k-1}
Y_{W}(v, \mp x)w)\nn
&\quad =\sum_{k=0}^{2n}(-1)^{k}
S_{n}(\res_{x}e^{2\pi inx}(e^{2\pi ix}-1)^{-k}
Y_{W}(v, \mp x)w)\nn
&\quad\quad +\sum_{k=0}^{2n}(-1)^{k}
S_{n}(\res_{x}e^{2\pi inx}(e^{2\pi ix}-1)^{-k-1}
Y_{W}(v, \mp x)w)\nn
&\quad =S_{n}(\res_{x}e^{2\pi inx}
Y_{W}(v, \mp x)w)+S_{n}(\res_{x}e^{2\pi inx}(e^{2\pi ix}-1)^{-2n-1}
Y_{W}(v, \mp x)w)\nn
&\quad =-S_{n}(\res_{x}e^{-2\pi inx}
Y_{W}(v, \pm x)w)-S_{n}(\res_{x}e^{-2\pi inx}(e^{-2\pi ix}-1)^{-2n-1}
Y_{W}(v, \pm x)w)\nn
&\quad =-S_{n}(\res_{x}e^{-2\pi inx}
Y_{W}(v, \pm x)w)+S_{n}(\res_{x}e^{2\pi i(n+1)x}(e^{2\pi ix}-1)^{-2n-1}
Y_{W}(v, \pm x)w). 
\end{align}
By (\ref{mod-inv-10.4}), 
we see that the second term in the right-hand side of 
(\ref{main-lemma-16.1}) is equal to
\begin{align}\label{main-lemma-16.2}
&-\sum_{j=0}^{n}\sum_{s=1}^{n}\sum_{l | s}
\frac{l}{(2n-j)}\binom{n}{j}\cdot\nn
&\quad\quad\quad\quad\quad\cdot 
 S_{n-s}\left(\res_{x}\left(\binom{-l-1}{2n-1-j}
e^{2l\pi ix}+\binom{l-1}{2n-1-j}e^{-2l\pi ix}\right)
Y_{W}(
v, \pm x)w\right)\nn
&\quad =-\sum_{s=1}^{n}\sum_{l | s}\sum_{j=0}^{n}
\frac{l}{(2n-j)}\binom{n}{j}\binom{-l-1}{2n-1-j}
 S_{n-s}(\res_{x}e^{2l\pi ix}Y_{W}(v, \pm x)w)\nn
 &\quad \quad -\sum_{s=1}^{n}\sum_{l | s}\sum_{j=0}^{n}
\frac{l}{(2n-j)}\binom{n}{j}\binom{l-1}{2n-1-j}
 S_{n-s}(\res_{x}e^{-2l\pi ix}
Y_{W}(v, \pm x)w).
\end{align}
Taking $m$ and $n$ be $n$ and $2n-1$, respectively,  in (\ref{lemma-3-1}),
we obtain 
$$\sum_{j=0}^{n}
\frac{l}{(2n-j)}\binom{n}{j}\binom{-l-1}{2n-1-j}=0$$
and thus the first term in the right-hand side of 
(\ref{main-lemma-16.2}) is equal to $0$. 
Note that $\binom{l-1}{2n-1-j}=0$ when $2n-j-1>l-1$ or equivalently 
when $2n>l+j$. But $l, j\le n$. So only when $l=j=n$, 
$\binom{l-1}{2n-1-j}$ is not $0$.  
Thus the second term in the right-hand side of 
(\ref{main-lemma-16.2}) is equal to 
$$-\sum_{s=1}^{n}\sum_{n | s}
\frac{n}{n}\binom{n}{n}\binom{n-1}{n-1}
 S_{n-s}(\res_{x}e^{-2n\pi ix}
Y_{W}(v, \pm x)w)=- S_{0}(\res_{x}e^{-2n\pi ix}
Y_{W}(v, \pm x)w).$$
From these calculations, we see that 
the right-hand side of 
(\ref{main-lemma-16.1}) is equal to (\ref{main-lemma-16}),
proving (\ref{main-lemma-4.1}) in the case $m=2n>0$.

In the case $2n>m>n$, 
we have
\begin{align}\label{main-lemma-17}
&\sum_{k=0}^{m}\binom{-2n+m-1}{k}
S_{n}(\res_{x}e^{2\pi i(n+1)x}(e^{2\pi ix}-1)^{-2n+m-k-1}
Y_{W}(v, \mp x)w)\nn
&\quad =\sum_{k=0}^{m}\binom{-2n+m-1}{k}
S_{n}(\res_{x}e^{2\pi i(m-n)x}e^{2\pi i(2n-m+1)x}
(e^{2\pi ix}-1)^{-2n+m-k-1}
Y_{W}(v, \mp x)w)\nn
&\quad=\sum_{k=0}^{m}\sum_{j=0}^{2n-m+1}\binom{-2n+m-1}{k}
\binom{2n-m+1}{j}\cdot\nn
&\quad\quad\quad \quad\quad\quad\quad\quad\quad
\cdot 
S_{n}(\res_{x}e^{2\pi i(m-n)x}(e^{2\pi ix}-1)^{-j-k}
Y_{W}(v, \mp x)w)\nn
&\quad=-\sum_{k=0}^{m}\sum_{j=0}^{2n-m+1}\binom{-2n+m-1}{k}
\binom{2n-m+1}{j}\cdot\nn
&\quad\quad\quad \quad\quad\quad\quad\quad\quad
\cdot 
S_{n}(\res_{x}e^{-2\pi i(m-n)x}(e^{-2\pi ix}-1)^{-j-k}
Y_{W}(v, \pm x)w)\nn
&\quad=-\sum_{k=0}^{m}\sum_{j=0}^{2n-m+1}(-1)^{-j-k}
\binom{-2n+m-1}{k}
\binom{2n-m+1}{j}\cdot\nn
&\quad\quad\quad \quad\quad\quad\quad\quad\quad
\cdot 
S_{n}(\res_{x}e^{2\pi i(-m+n+j+k)x}(e^{2\pi ix}-1)^{-j-k}
Y_{W}(v, \pm x)w)\nn
&\quad= -\sum_{r=0}^{2n+1}(-1)^{-r}\sum_{k=0}^{\min(m, r)}
\binom{-2n+m-1}{k}
\binom{2n-m+1}{r-k}\cdot\nn
&\quad\quad\quad \quad\quad\quad\quad\quad\quad
\cdot 
S_{n}(\res_{x}e^{2\pi i(-m+n+r)x}(e^{2\pi ix}-1)^{-r}
Y_{W}(v, \pm x)w)\nn
&\quad= -\sum_{r=0}^{m}(-1)^{-r}\sum_{k=0}^{r}
\binom{-2n+m-1}{k}
\binom{2n-m+1}{r-k}\cdot\nn
&\quad\quad\quad \quad\quad\quad\quad\quad\quad
\cdot 
S_{n}(\res_{x}e^{2\pi i(-m+n+r)x}(e^{2\pi ix}-1)^{-r}
Y_{W}(v, \pm x)w)\nn
&\quad\quad  -\sum_{r=m+1}^{2n+1}(-1)^{-r}\sum_{k=0}^{m}
\binom{-2n+m-1}{k}
\binom{2n-m+1}{r-k}\cdot\nn
&\quad\quad\quad \quad\quad\quad\quad\quad\quad
\cdot 
S_{n}(\res_{x}e^{2\pi i(-m+n+r)x}(e^{2\pi ix}-1)^{-r}
Y_{W}(v, \pm x)w).
\end{align}

Using (\ref{s-p-binomial-1}), we see that the first term 
in the right-hand side of (\ref{main-lemma-17}) is equal to 
$$-S_{n}(\res_{x}e^{2\pi i(m-n)x}
Y_{W}(v, \pm x)w).$$
Since $m>n\ge l$ and $r\ge m+1$, we have $r-2\ge -m+n+r-1\ge n\ge l$. 
Then by (\ref{lemma-3-1}) with $m, n$ replaced by 
$-m+n+r-1, r-2$, respectively, we have
\begin{equation}\label{main-lemma-18}
\sum_{j=0}^{-m+n+r-1}\frac{l}{(r-j-1)}\binom{-m+n+r-1}{j}
\binom{-l-1}{r-2-j}=0.
\end{equation}
Using (\ref{mod-inv-10.4}) with $m, n$ replaced by 
$-m+n+r-1, r-2$, respectively, and (\ref{main-lemma-18}), 
we see that the second term in
the right-hand side of (\ref{main-lemma-17}) is equal to 
\begin{align}\label{main-lemma-18.5}
& \sum_{r=m+1}^{2n+1}(-1)^{-r}\sum_{k=0}^{m}
\binom{-2n+m-1}{k}
\binom{2n-m+1}{r-k}\cdot\nn
&\quad\quad\quad \quad\quad\quad\quad\quad
\cdot 
\sum_{j=0}^{-m+n+r-1}\sum_{s=1}^{n}\sum_{l | s}
\frac{l}{(r-j-1)}\binom{-m+n+r-1}{j}\cdot\nn
&\quad\quad\quad\quad\quad\quad\quad\quad\cdot 
 S_{n-s}\left(\res_{x}\left(\binom{-l-1}{r-2-j}
e^{2l\pi ix}+\binom{l-1}{r-2-j}e^{-2l\pi ix}\right)
Y_{W}(
v, \pm x)w\right)\nn
&\quad=\sum_{s=1}^{n}\sum_{l | s} \sum_{r=m+1}^{2n+1}(-1)^{-r}\sum_{k=0}^{m}
\sum_{j=0}^{-m+n+r-1}
\binom{-2n+m-1}{k}\binom{2n-m+1}{r-k}\cdot\nn
&\quad\quad\quad\quad\quad\quad\quad\quad\cdot 
\binom{-m+n+r-1}{j}\binom{l}{r-1-j}
 S_{n-s}\left(\res_{x}e^{-2l\pi ix}
Y_{W}(
v, \pm x)w\right)\nn
&\quad=  -\sum_{s=1}^{n}\sum_{l | s} \sum_{p=0}^{2n-m}(-1)^{-p-m}\sum_{k=0}^{m}
\binom{-2n+m-1}{k}\binom{2n-m+1}{p+m+1-k}\cdot\nn
&\quad\quad\quad\quad\quad\quad\quad\quad\cdot 
\sum_{j=0}^{m+p}\binom{n+p}{j}\binom{l}{p+m-j}
 S_{n-s}\left(\res_{x}e^{-2l\pi ix}
Y_{W}(
v, \pm x)w\right).
\end{align}
Using (\ref{s-p-binomial-1}) with $\alpha, \beta, m$ replaced by 
$n+p, l, p+m$, respectively, we see that the right-hand side of 
(\ref{main-lemma-18.5}) is equal to 
\begin{align}\label{main-lemma-19}
& -\sum_{s=1}^{n}\sum_{l | s} \sum_{p=0}^{2n-m}(-1)^{-p-m}\sum_{k=0}^{m}
\binom{-2n+m-1}{k}\binom{2n-m+1}{p+m+1-k}\cdot\nn
&\quad\quad\quad\quad\quad\quad\quad\quad\cdot 
\binom{n+p+l}{p+m}
 S_{n-s}\left(\res_{x}e^{-2l\pi ix}
Y_{W}(
v, \pm x)w\right).
\end{align}

Using (\ref{main-comb-formu}), we
see that (\ref{main-lemma-19})
is equal to 
\begin{align*}
&-\sum_{s=1}^{n}\sum_{l | s}\delta_{l, m-n}
 S_{n-s}\left(\res_{x}e^{-2l\pi ix}
Y_{W}(
v, \pm x)w\right)\nn
&\quad =-\sum_{(m-n) | s, \; 1\le s\le n}
 S_{n-s}\left(\res_{x}e^{-2(m-n)\pi ix}
Y_{W}(
v, \pm x)w\right).
\end{align*}
From the calculations above, we obtain that the right-hand side of 
(\ref{main-lemma-17}) is equal to 
\begin{align*}
&-S_{n}(\res_{x}e^{2\pi i(m-n)x}
Y_{W}(v, \pm x)w)-\sum_{(m-n) | s, \; 1\le s\le n}
 S_{n-s}\left(\res_{x}e^{-2(m-n)\pi ix}
Y_{W}(
v, \pm x)w\right)\nn
&\quad=\sum_{(m-n) | s, \; 1\le s\le m}
 S_{m-s}\left(\res_{x}e^{-2(m-n)\pi ix}
Y_{W}(
v, \pm x)w\right).
\end{align*}
Thus (\ref{main-lemma-4.1}) holds in the case $2n>m>n$.
This finishes the proof of (\ref{main-lemma-4.1}) in the case $m> n$.
The proof of  (\ref{main-lemma-4}) is now complete.
\epfv

\begin{rema}\label{miyamoto-prop-4.4}
{\rm The special case $n=m$ and $W=V$ of (\ref{main-lemma-3}) 
is the same as Proposition 4.4 
in \cite{M}. In fact, this proposition is not proved
in \cite{M}. The proof of Proposition 4.4 
in \cite{M} uses some formulas obtained by Zhu for shifted
$q$-traces of vertex operators. It is claimed in \cite{M}
that these formulas give some properties of $O_{q}(V)$. 
But actually these formulas of Zhu show that 
such properties hold for the kernels of suitable linear maps 
constructed from shifted $q$-traces of
vertex operators. Although
$O_{q}(V)$ is a subspace of these kernels, 
one cannot conclude from only these facts
that $O_{q}(V)$ also satisfies the same properties. 
The arguments  
in \cite{M} indeed give strong evidence 
 that the conclusion of Proposition 4.4 
in \cite{M} must be true. But these arguments do not give a proof. 
A proof is given by the proof of the special case $n=m$ and $W=V$ 
of (\ref{main-lemma-3}). In fact, 
the first proof of this proposition was given by
McRae \cite{Mc}. }
\end{rema}

\begin{prop}\label{2nd-main-lemma}
Let $W$ be a grading-restricted generalized $V$-module. 
Assume that a linear map
\begin{align*}
S: W&\to q^{r}\C[[q]][\log q]\nn
w&\mapsto S(w; q)=\sum_{k=0}^{K}\sum_{n\in \N}S_{k, n}(w)
(\log q)^{k}q^{r+n}
\end{align*} 
satisfies (\ref{main-lemma-1}), (\ref{main-lemma-2}) and
(\ref{2nd-main-lemma-1}). Then
\begin{align}\label{grading-slf}
&\sum_{m=0}^{n}\binom{-n-1}{m}S_{k, n}(\res_{x}e^{2\pi (n+1)ix}
(e^{2\pi ix}-1)^{-n-m-1}Y_{W}(\omega, x)w)
-2\pi i(r+n)S_{k, n}(w)\nn
&\quad =
2\pi i(k+1)S_{k+1, n}(w)
\end{align}
for $k=0, \dots, K$, $n\in \N$ and $w\in W$, where $S_{K+1, n}(w)=0$
for $n\in \N$ and $w\in W$. 
\end{prop}
\pf
From  (\ref{2nd-main-lemma-1}), (\ref{tilde-wp-g-2}) and 
(\ref{main-lemma-1}), 
we obtain
\begin{align}\label{tilde-wp-expansion}
&(2\pi i)^{2}\sum_{k=0}^{K}
\sum_{n\in \N}kS_{k, n}(w)(\log q)^{k-1}
q^{r+n}+(2\pi i)^{2}\sum_{k=0}^{K}
\sum_{n\in \N}(r+n)S_{k, n}(w)
(\log q)^{k}q^{r+n}\nn
&\quad=2\pi i\sum_{k=0}^{K}\sum_{n\in \N}S_{k, n}(\res_{x}
e^{2\pi ix}(e^{2\pi ix}-1)^{-1}
Y_{W_{1}}(\omega, x)w)(\log q)^{k}q^{r+n}\nn
&\quad\quad -2\pi i\sum_{k=0}^{K}\sum_{n\in \N}\sum_{s\in \Z_{+}}
\sum_{l|s}S_{k, n}(\res_{x}
(e^{2l\pi ix}-e^{-2l\pi ix})
Y_{W_{1}}(\omega, x)w)(\log q)^{k}q^{r+n+s}\nn
&\quad=2\pi i\sum_{k=0}^{K}\sum_{n\in \N}S_{k, n}(\res_{x}
e^{2\pi ix}(e^{2\pi ix}-1)^{-1}
Y_{W_{1}}(\omega, x)w)(\log q)^{k}q^{r+n}\nn
&\quad\quad -2\pi i\sum_{k=0}^{K}
\sum_{n\in \Z_{+}}\sum_{s=1}^{n}
\sum_{l|s}S_{k, n-s}(\res_{x}
(e^{2l\pi ix}-e^{-2l\pi ix})
Y_{W_{1}}(\omega, x)w)(\log q)^{k}q^{r+n}.
\end{align}
Taking the coefficients of $(\log q)^{k}q^{r+n}$ in 
(\ref{tilde-wp-expansion}), we obtain
\begin{equation}\label{omega-00}
S_{k,  0}(\res_{x}
e^{2\pi ix}(e^{2\pi ix}-1)^{-1}
Y_{W_{1}}(\omega, x)w)-2\pi irS_{k, 0}(w)
=2\pi i(k+1)S_{k+1, 0}(w)
\end{equation}
for $k=0, \dots, K$, which is (\ref{grading-slf}) in the case of $n=0$, 
and 
\begin{align}\label{2nd-main-lemma-2}
&S_{k, n}(\res_{x}
e^{2\pi ix}(e^{2\pi ix}-1)^{-1}
Y_{W_{1}}(\omega, x)w)-\sum_{s=1}^{n}
\sum_{l|s}S_{k, n-s}(\res_{x}
(e^{2l\pi ix}-e^{-2l\pi ix})
Y_{W_{1}}(\omega, x)w)\nn
&\quad=2\pi i(k+1)S_{k+1, n}(w)+2\pi i(r+n)S_{k, n}(w)
\end{align}
for $k=0, \dots, K$ and $m\in \Z_{+}$. 

Note that for $k=0, \dots, K$, 
$\sum_{n\in \N}S_{k, n}(w)q^{r+n}$
satisfy (\ref{main-lemma-1}) and (\ref{main-lemma-2}).
Thus for fixed $k$, all the results above hold for $S_{k, n}$ for $n\in \N$. 
Expanding $e^{2\pi nix}=(1+(e^{2\pi ix}-1))^{n}$
as a polynomial in $e^{2\pi ix}-1$ and then using 
(\ref{mod-inv-10.4}) for $S_{k, n}$, we have
\begin{align}\label{2nd-main-lemma-3}
&\sum_{m=0}^{n}\binom{-n-1}{m}S_{k, n}(\res_{x}e^{2\pi (n+1)ix}
(e^{2\pi ix}-1)^{-n-m-1}Y_{W}(\omega, x)w)\nn
&\quad=\sum_{m=0}^{n}\sum_{j=0}^{n}\binom{-n-1}{m}\binom{n}{j}
S_{k, n}(\res_{x}e^{2\pi ix}
(e^{2\pi ix}-1)^{-n+j-m-1}Y_{W}(\omega, x)w)\nn
&\quad=S_{k, n}(\res_{x}e^{2\pi ix}
(e^{2\pi ix}-1)^{-1}Y_{W}(\omega, x)w)\nn
&\quad\quad -\sum_{j=0}^{n-1}
\sum_{s=1}^{n}
\sum_{l | s}\binom{n}{j}\frac{l}{n-j}\cdot \nn
&\quad\quad\quad \quad\cdot
S_{k, n-s}\left(\res_{x}\left(\binom{-l-1}{n-j-1}e^{2\pi i lx}
+\binom{l-1}{n-j-1}e^{-2\pi i lx}\right)Y_{W}(\omega, x)w\right)\nn
&\quad\quad -\sum_{m=1}^{n}\sum_{j=0}^{n}
\sum_{s=1}^{n}
\sum_{l | s}\binom{-n-1}{m}\binom{n}{j}\frac{l}{n-j+m}\cdot \nn
&\quad\quad\quad \quad\cdot
S_{k, n-s}\left(\res_{x}\left(\binom{-l-1}{n-j+m-1}e^{2\pi i lx}
+\binom{l-1}{n-j+m-1}e^{-2\pi i lx}\right)Y_{W}(\omega, x)w\right).
\end{align}

Using (\ref{2nd-main-lemma-6}), (\ref{2nd-main-lemma-7}),
(\ref{2nd-main-lemma-8}) and (\ref{2nd-main-lemma-9}), 
we see that the right-hand side of (\ref{2nd-main-lemma-3})
is equal to the left-hand side of  (\ref{2nd-main-lemma-2}). Then by 
(\ref{2nd-main-lemma-2}), we obtain (\ref{grading-slf}).
\epfv

\begin{rema}\label{miyamoto-prop-4.5}
{\rm The special case $W=V$ of (\ref{grading-slf}) implies immediately 
Proposition 4.5
in \cite{M}. In fact, this proposition is not proved
in \cite{M}. Note that
Proposition 4.5 in \cite{M} is based on the same arguments as what
Proposition 4.4 in \cite{M} is based on. As we have discussed in 
Remark \ref{miyamoto-prop-4.4}, though
these arguments indeed give strong evidence 
 that the conclusion of Proposition 4.5  in \cite{M} 
must be true, they do not give a proof. A proof of this proposition 
is now given by the proof of the special case $W=V$ of (\ref{grading-slf}). 
We also note that \cite{M} does not have a proof of
Proposition 4.6 in \cite{M}. Since  
we do not need anything equivalent to Proposition 4.6 in \cite{M}
in this paper,
we do not give a proof. On the other hand, without a proof 
of Proposition 4.6 in \cite{M}, the proof of the modular invariance 
theorem in \cite{M} is not complete. Certainly the 
modular invariance theorem in \cite{M} is a special case of 
Theorem \ref{mod-inv} below and thus we do obtain a complete
proof in this paper. But it is still interesting to see whether the 
method in \cite{M} can indeed give a complete proof of 
the modular invariance theorem in \cite{M}. To see this, 
we need to find a proof
of Proposition 4.6 in \cite{M}. }
\end{rema}

We now prove the main result of this section. 

\vspace{1em}

\noindent {\it Proof of Theorem \ref{wp-fn-gr-sym-fn}}. 
Since $\psi_{S}^{N; k, j}([v]_{mn})=0$ for $m, n\in \N$ satisfying
$0\le m, n\le N$ and $m\ne n$, to show 
$\psi^{N}_{S; k, j}(\widetilde{Q}^{\infty}(W))=0$, 
we need only show $\psi_{S; k, j}^{N}([w]_{nn})=0$ for $w\in W$ 
and $0\le n\le N$
such that $[w]_{nn}\in \widetilde{Q}^{\infty}(W)$. 

Since $S$ satisfies 
(\ref{main-lemma-1}) and (\ref{main-lemma-2}), we see that
for $k=0, \dots, K$ and $j=1, \dots, J$, the linear map given by
$w\mapsto \sum_{m\in \N}S_{k, j, m}(w)q^{m}$ for $w\in W$ 
also satisfies (\ref{main-lemma-1}) and (\ref{main-lemma-2}).
In particular, (\ref{main-lemma-3}) holds for $S_{p}=S_{k, j, p}$.
Then by (\ref{main-lemma-3}) for $S_{n}=S_{k, j, n}$, we
have
$$S_{k, j, n}(\res_{x}e^{2\pi i(n+1)x}
(e^{2\pi ix}-1)^{-2n-p-2}
Y_{W}(v, x)w)=0$$
for $p\in \N$. 
Recall the invertible operator $\mathcal{U}(1): W\to W$
in Subsection 2.2. 
Changing the variable in the formal residue and using 
(\ref{U-1-v-op}), we obtain 
\begin{align*}
&S_{k, j, n}(\mathcal{U}(1)^{-1}
\res_{x}x^{-2n-2-p}
Y_{W}((1+x)^{L_{V}(0)+n}v, x)w)\nn
&\quad =S_{k, j, n}\left(
\res_{x}x^{-2n-2-p}(1+x)^{n}
Y_{W}\left(\mathcal{U}_{V}(1)^{-1}v, \frac{1}{2\pi i}
\log (1+x)\right)\mathcal{U}_{W}(1)^{-1}w\right)\nn
&\quad =2\pi iS_{k, j, n}(\res_{x}e^{2\pi i(n+1)x}
(e^{2\pi ix}-1)^{-2n-p-2}
Y_{W}(\mathcal{U}_{V}(1)^{-1}v, x)\mathcal{U}_{W}(1)^{-1}w)\nn
&\quad =0
\end{align*}
for $v\in V$, $w\in W$ and $p\in \N$. 
Since elements of $W$ the form 
$\res_{x}x^{-2n-2-p}
Y_{W}((1+x)^{L_{V}(0)+n}v, x)w$ span 
$O_{n}(W)$ (see \cite{HY}), 
the linear function $S_{k, j, n}\circ \mathcal{U}(1)^{-1}$
on $W$ is in fact a linear function on $A_{n}(W)=V/O_{n}(V)$ (see \cite{HY}). 
Changing the variable in the formal residue and 
using (\ref{U-1-v-op}), we see that in the case of $m=n$, 
(\ref{main-lemma-4}) in fact gives 
$$S_{k, j, n}( \mathcal{U}(1)^{-1}(v*_{n}w))
=S_{k, j, n}( \mathcal{U}(1)^{-1}(w*_{n}v))$$
(see \cite{HY}) for $v\in V$ and $w\in W$. This means that 
$S_{k, j, n}\circ \mathcal{U}(1)^{-1}$ is in fact a 
symmetric linear function on $A_{n}(W)$. 

Since $V$ is $C_2$-cofinite and has no nonzero elements
of negative weights, by Theorem \ref{properties},
$A^{n}(V)$ is finite dimensional. But from \cite{H-aa-va},
$A_{n}(V)$ is isomorphic to a subalgebra of 
$A^{n}(V)$ and is thus also finite dimensional. 
By Theorem \ref{finite-d-A-N-W}, $A^{n}(W)$ is also finite dimensional. 
By Theorem \ref{isom-thm-A-N-W}, $A_{n}(W)$ is linearly isomorphic to 
a subspace of $A^{n}(W)$ and is thus also finite dimensional.

By Theorem \ref{slf=>p-tr}, 
there exist finite-dimensional 
basic symmetric algebras $P_{i}$ for $i=1, \dots, l$
equipped with symmetric linear functions $\phi_{i}$, 
finite-dimensional $A_{n}(V)$-$P_{i}$-bimodules $U_{i}$ (projective as 
right $P_{i}$-modules) and $f_{i}\in \hom_{A_{n}(V), P_{i}}(A_{n}(W)
\otimes_{A_{n}(V)}U_{i}, U_{i})$
such that for $w\in W$
$$S_{k, j, n}( \mathcal{U}(1)^{-1}w)=\sum_{i=1}^{l}
(\phi_{i})^{f_{i}}_{U_{i}}(w+O_{n}(W)).$$
For a lower-bounded generalized $V$-module
$W$ and $n\in \N$, we have $G_{n}(W)=W_{\l n\r}
=\coprod_{\mu\in \Gamma(W)}W_{[h^{\mu}+n]}$. 
In \cite{HY}, a lower-bounded generalized $V$-module 
$S_{n}(G_{n}(W))$ is constructed such that 
$G_{n}(S_{n}(G_{n}(W)))$ as an $A_{n}(V)$-module is 
equivalent to $G_{n}(W)$. In fact, this construction 
does not have to start from $G_{n}(W)$ of  
a lower-bounded generalized $V$-module
$W$. We can start with an arbitrary $A_{n}(V)$-module
$M$ and the construction gives us $S_{n}(M)$ is constructed such that 
$G_{n}(S_{n}(M))$ as an $A_{n}(V)$-module is 
equivalent to $M$. One can also obtain $S_{n}(M)$ using the construction 
given in \cite{H-const-twisted-mod}  and \cite{H-affine-twisted-mod}.

For $i=1, \dots, l$, we have lower-bounded generalized 
$V$-modules $S_{n}(U_{i})$ and $S_{n}(U_{i}^{*})$. 
Since $U_{i}$ and $U_{i}^{*}$ for $i=1, \dots, n$ 
are finite dimensional, by Proposition \ref{properties},
$S_{n}(U_{i})$ and $S_{n}(U_{i}^{*})$ are grading restricted.
In particular, $S_{n}(U_{i}^{*})''=S_{n}(U_{i}^{*})$. 
Let $W_{2}^{i}=S_{n}(U_{i})$ and $W_{3}^{i}=S_{n}(U_{i}^{*})'$.
Then $W_{2}^{i}=S_{n}(U_{i})=S_{n}(G_{n}(W_{2}^{i})$ and 
$(W_{3}^{i})'=S_{n}(U_{i}^{*})=S_{n}(G_{n}((W_{3}^{i})')$.

In the proof of Theorem 6.7 in \cite{H-aa-int-op}, an element 
$f\in \hom_{A^{\infty}(V)}(A^{\infty}(W)\otimes_{A^{\infty}(V)}
W_{2}, W_{3})$ is constructed from 
an element $f^{N}\in \hom_{A^{N}(V)}(A^{N}(W)\otimes_{A^{N}(V)}
\Omega_{N}^{0}(W_{2}), \Omega_{N}^{0}(W_{3}))$ such that 
the restriction of $f$ to $A^{N}(W)\otimes_{A^{N}(V)}
\Omega_{N}^{0}(W_{2})$ gives $f^{N}$. 
The same proof, with $A^{N}(V)$, $A^{N}(W)$, $\Omega_{N}^{0}(W_{2})$,
$\Omega_{N}^{0}(W_{3})$ and $A^{N, \infty}(W)$
replaced by $A_{n}(V)$, $A_{n}(W)$, $G_{n}(W_{2})$, $G_{n}(W_{3})$
and the subspace of $A^{\infty}(W)$ consisting of sums of 
elements of the form $\sum_{i=0}^{n}[w]_{ii}$ for $w\in W$ and 
$[w]_{ij}$ for $w\in W$ and $i, j\in N$ satisfying $i<j$, respectively,
in fact gives a construction of an element 
$\tilde{f}\in \hom_{A^{\infty}(V)}(A^{\infty}(W)\otimes_{A^{\infty}(V)}
W_{2}, W_{3})$ from an element $f\in 
\hom_{A_{n}(V)}(A_{n}(W)
\otimes_{A_{n}(V)}G_{n}(W_{2}), G_{n}(W_{3}))$
such that the restriction of $\tilde{f}$ to 
$A_{n}(W)
\otimes_{A_{n}(V)}G_{n}(W_{2})$ gives $f$. 
Apply this construction to $f_{i}\in \hom_{A_{n}(V), P_{i}}(A_{n}(W)
\otimes_{A_{n}(V)}U_{i}, U_{i})$ above for $i=1, \dots, l$, we obtain 
$\tilde{f}_{i}\in \hom_{A^{\infty}(V)}(A^{\infty}(W)\otimes_{A^{\infty}(V)}
W_{2}^{i}, W_{3}^{i})$ such that the restriction of 
$\tilde{f}_{i}$ to $A_{n}(W)
\otimes_{A_{n}(V)}U_{i}$ gives $f_{i}$. 
In particular, for $w\in W$ such that 
$[w]_{nn}\in Q^{\infty}(W)$ and $w_{2}^{i}\in U_{i}$, 
\begin{align*}
f_{i}((w+O_{n}(W))\otimes_{A_{n}(V)}w_{2})
&=\tilde{f}_{i}(([w]_{nn}+Q^{\infty}(W))\otimes_{A^{\infty}(V)}
w_{2})\nn
&=0.
\end{align*}
Then by the definition of $(\phi_{i})^{f_{i}}_{U_{i}}$, we 
have
$(\phi_{i})^{f_{i}}_{U_{i}}(w+O_{n}(W))=0$ for $i=1, \dots, l$
for $w\in W$ such that $[w]_{nn}\in Q^{\infty}(W)$.
Then 
\begin{align*}
\psi_{S; k, j}^{N}([\mathcal{U}(1)^{-1}w]_{nn})
&=S_{k, j, n}(\mathcal{U}(1)^{-1}w)\nn
&=\sum_{i=1}^{l}
(\phi_{i})^{f_{i}}_{U_{i}}(w+O_{n}(W))\nn
&=0.
\end{align*}
for $w\in W$ such that $[w]_{nn}\in Q^{\infty}(W)$.
Since $\widetilde{Q}^{\infty}(W)= \mathcal{U}(1)^{-1}Q^{\infty}(W)$
and $0\le n\le N$,
we have proved $\psi_{S; k, j}^{N}([w]_{nn})=0$ for $w\in W$
and $0\le n\le N$ such that $[w]_{nn}\in \widetilde{Q}^{\infty}(W)$. 
Thus we have $\psi^{N}_{S; k, j}(\widetilde{Q}^{\infty}(W))=0$.

To prove $\psi_{S; k, j}^{N}$ is symmetric, we need only 
prove 
\begin{equation}\label{symmetry}
\psi_{S; k, j}^{N}([v]_{mn}\bkdiamd [w]_{nm})
=\psi_{S; k, j}^{N}( [w]_{nm} \bkdiamd [v]_{mn})
\end{equation}
for $v\in V$, $w\in W$ and $0\le m, n\le N$. 
Since for $k=0, \dots, K$ and $j=1, \dots, J$, 
the linear map given by
$w\mapsto \sum_{m\in \N}S_{k, j, m}(w)$ for $w\in W$ 
satisfies (\ref{main-lemma-1}) and (\ref{main-lemma-2}),
we see that (\ref{main-lemma-4}) holds for $S_{k, j, m}$ for $m, n\in \N$.
But by the definitions of 
$[v]_{mn}\bkdiamd [w]_{nm}$, $[w]_{nm} \bkdiamd [v]_{mn}$
and $\psi_{S; k, j}^{N}$, (\ref{symmetry}) is the same as 
(\ref{main-lemma-4}). 

Since $S$ satisfies (\ref{2nd-main-lemma-1}), we see that 
for $j=1, \dots, J$, the linear map given by
$w\mapsto \sum_{k=0}^{K}\sum_{m\in \N}S_{k, j, m}(w)$ for $w\in W$ 
also satisfies (\ref{2nd-main-lemma-1}). 
In particular, (\ref{grading-slf}) holds for 
$S_{k, j, n}$ for $j=1, \dots, J$. 
By the definition of $\psi_{S; k, j}^{N}$ and $[\omega]_{nn} \bkdiamd 
[w]_{nn}$, we see that the first term in the left-hand side of
(\ref{grading-slf}) multiplied by $2\pi i$ is equal to 
$\psi_{S; k, j}^{N}([\omega]_{nn} \bkdiamd 
[w]_{nn})$
and thus from (\ref{grading-slf}) and $[\one]_{nn}\bkdiamd [w]_{nn}
=[w]_{nn}$, we obtain
\begin{equation}\label{grading-sym-linear-fn-2}
\psi_{S; k, j}^{N}(([\omega]_{nn}-(2\pi i)^{2}
(r_{j}+n)[\one]_{nn}) \bkdiamd 
[w]_{nn})=(2\pi i)^{2}(k+1)\psi_{S; k+1, j}^{N}([w]_{nn}).
\end{equation}
Using (\ref{grading-sym-linear-fn-2}) repeatedly and 
note that $\psi_{S; K+1, j}^{N}([w]_{nn})=S_{K+1, j, n}(w)=0$, 
we obtain 
$$\psi_{S; k, j}^{N}(([\omega]_{nn}
-(2\pi i)^{2}(r_{j}+n)[\one]_{nn})^{\bkdiamd (K-k+1)}
\bkdiamd [w]_{nn})=0$$
for $w\in W$ and $n=0, \dots, N$. For $k, l\in \N$, $k\ne n$ and $w\in W$,
we have 
$$\psi_{S; k, j}^{N}(([\omega]_{nn}
-(2\pi i)^{2}(r_{j}+n)[\one]_{nn})^{\bkdiamd (K-k+1)}
\bkdiamd [w]_{kl})=\psi_{S; k, j}^{N}(0)=0.$$
If for $l\in \N$ and $l\ne n$, note that 
$$([\omega
-(2\pi i)^{2}
(r_{j}+n)\one]_{nn})^{\bkdiamd (K-k+1)}
\bkdiamd [w]_{nl}=[\tilde{w}]_{nl},$$
for some $\tilde{w}\in W$. Then by the definition of 
$\psi_{S; k, j}^{N}$
we have 
$$\psi_{S; k, j}^{N}(([\omega]_{nn}
-(2\pi i)^{2}(r_{j}+n)[\one]_{nn})^{\bkdiamd (K-k+1)}
\bkdiamd [w]_{nl})=\psi_{S; k, j}^{N}([\tilde{w}]_{nl})=0.$$
Thus (\ref{grading-sym-linear-fn-1}) is proved.
\epfv

\renewcommand{\theequation}{\thesection.\arabic{equation}}
\renewcommand{\thethm}{\thesection.\arabic{thm}}
\setcounter{equation}{0} \setcounter{thm}{0} 

\section{Modular invariance}

In this section, we prove the conjecture on 
the modular invariance of intertwining operators. 
See Theorem \ref{mod-inv}. 
From the genus-one associativity, as in the 
proof of the modular invariance in the semisimple case 
in \cite{H-mod-inv-int}, we see that
we need only prove that the space of all genus-one $1$-point 
correlation functions constructed from shifted pseudo-$q$-traces 
of geometrically-modified intertwining operators
are invariant under the modular transformations. We introduce a notion
of genus-one $1$-point  conformal block as a map from a 
grading-restricted generalized $V$-module $W$ to the space of analytic functions
on the upper-half plane $\H$ satisfying some properties involving the 
Weierstrass functions $\wp$-function $\wp_{2}(z; \tau)$, 
the Weierstrass $\zeta$-function $\wp_{1}(z; \tau)$
and the Eisenstein series $G_{2}(\tau)$. We prove that modular transformations of 
genus-one $1$-point 
correlation functions constructed from shifted pseudo-$q$-traces 
 of geometrically-modified intertwining operators give 
genus-one $1$-point  conformal blocks. 
Then we prove that the image of $w\in W$ under a genus-one 
$1$-point conformal block must be a sum of 
genus-one $1$-point 
correlation functions constructed from shifted pseudo-$q$-traces 
 of intertwining operators. This in particular 
proves the modular invariance conjecture. 

We now introduce the notion of genus-one $1$-point conformal block.
Let $\H$ be the open upper-half plane and $H(\H)$ the space of 
analytic functions on $\H$. Recall the formal Laurent series expansion
$\wp_{1}(x; \tau)-G_{2}(\tau)x$ of $\wp_{1}(z; \tau)-G_{2}(\tau)z$ 
given by (\ref{wp-2-x})
and the formal Laurent series expansion 
$\wp_{2}(x; \tau)$ of $\wp_{2}(z; \tau)$ given by (\ref{wp-1-g-2-x}).

\begin{defn}\label{1-pt-cb}
{\rm 
Let  $W$  be a grading-restricted generalized 
$V$-module. A {\it  
genus-one $1$-point conformal block labeled by $W$} is a linear map
\begin{align*}
F: W&\to H(\H)\nn
w&\mapsto F(w; \tau)
\end{align*}
satisfying 
\begin{align}
&F(\res_{x}Y_{W}(v, x)w; \tau)=0,\label{1-pt-cb-1}\\
&2\pi i\frac{\partial}{\partial \tau}F(w; \tau)
=F(\res_{x}(\wp_{1}(x; \tau)-G_{2}(\tau)x)
Y_{W}(\omega, x)w; \tau),\label{1-pt-cb-2}\\
&F(\res_{x}\wp_{2}(x; \tau)
Y_{W}(v, x)w;
\tau)=0\label{1-pt-cb-3}
\end{align}
for $v\in V$ and $w\in W$. }
\end{defn}

Let $W$ be grading-restricted generalized $V$-modules
and $w\in W$. Given a finite-dimension\-al 
associative algebra $P$ with a symmetric linear function $\phi$, 
a grading-restricted generalized $V$-$P$-bimodule 
$\widetilde{W}$, projective as a right $P$-module, and an
intertwining operator $\Y$ of type
$\binom{\widetilde{W}}{W\widetilde{W}}$
compatible with $P$, 
we have a genus-one $1$-point correlation function
$\overline{F}^{\phi}_{\mathcal{Y}}(w;
z; \tau)$, which is the analytic extension of 
the sum of the series
$$F^{\phi}_{\mathcal{Y}}(w; z; \tau)
=\tr_{\widetilde{W}}^{\phi}\mathcal{Y}(\mathcal{U}_{W}(q_{z})w, q_{z})
q_{\tau}^{L_{\widetilde{W}}(0)-\frac{c}{24}}.$$
Using the $L(-1)$-derivative property and $L(-1)$-commutator formula for the 
intertwining operator $\Y$, we have
\begin{align*}
&\frac{\partial}{\partial z}\tr_{\widetilde{W}}^{\phi}
\mathcal{Y}(\mathcal{U}_{W}(q_{z})w, q_{z})
q_{\tau}^{L_{\widetilde{W}_{1}}(0)-\frac{c}{24}}\nn
&\quad =\tr^{\phi}_{\tilde{W}}\mathcal{Y}_{1}((2\pi iL_{W}(0)+
2\pi iq_{z}L_{W}(-1))
\mathcal{U}_{W}(q_{z})w, 
q_{z})q_{\tau}^{L_{\widetilde{W}}(0)-\frac{c}{24}}\nn
&\quad =2\pi i\tr^{\phi}_{\tilde{W}}[L_{\widetilde{W}}(0), \mathcal{Y}(
\mathcal{U}_{W}(q_{z})w, 
q_{z})]q_{\tau}^{L_{\widetilde{W}}(0)-\frac{c}{24}}\nn
&\quad=0.
\end{align*}
Since $\overline{F}^{\phi}_{\Y}(w; z; \tau)$ is the analytic extension of 
$\tr_{\widetilde{W}}^{\phi}\mathcal{Y}(\mathcal{U}_{W}(q_{z})w, q_{z})
q_{\tau}^{L_{\widetilde{W}}(0)-\frac{c}{24}}$, 
we obtain 
\begin{equation}\label{mod-inv-0}
\frac{d}{dz}\overline{F}^{\phi}_{\Y_{1}}(w; z; \tau)=0,
\end{equation}
that is, $\overline{F}^{\phi}_{\Y}(w; z; \tau)$ is independent of $z$. 
From now on, we shall write $\overline{F}^{\phi}_{\Y}(w; z; \tau)$ simply as 
$\overline{F}^{\phi}_{\Y}(w; \tau)$.
Then we have a linear map
$\overline{F}^{\phi}_{\mathcal{Y}}: W\to H(\mathbb{H})$
given by 
$w\mapsto \overline{F}^{\phi}_{\mathcal{Y}}(w; \tau)$. 
For $a\in \C^{\times}$, we define $a^{-L_{W}(0)}$ to 
be $a^{-L_{W}(0)_{S}}e^{-(\log a)L_{W}(0)_{N}}$, where, as is explicitly 
discussed in Subsection 2.2, $\log a=
\log |a|+\arg a$, $0\le \arg a<2\pi$. 

\begin{prop}\label{mod-transf-conf-bk}
For 
$$\left(\begin{array}{cc}
\alpha&\beta\\
\gamma&\delta
\end{array}\right)\in SL(2, \mathbb{Z}),$$
the linear map from $W$ to $H(\mathbb{H})$ given by 
\begin{align}\label{mod-transf-conf-bk-1}
w\mapsto \overline{F}^{\phi}_{\mathcal{Y}}
\Biggl((\gamma\tau+\delta)^{-L_{W}(0)}w;
\frac{\alpha\tau+\beta}{\gamma\tau+\delta}\Biggr)
\end{align}
is a genus-one $1$-point conformal block labeled by $W$.
\end{prop}
\pf 
The identity  (2.19) in \cite{F} in the case $n=1$ becomes
\begin{align}\label{mod-inv-1}
&\left(2\pi i\frac{\partial}{\partial \tau}
+G_{2}(\tau)z\frac{\partial}{\partial z}\right)
F^{\phi}_{\mathcal{Y}}(w; \tau)
+G_{2}(\tau)F^{\phi}_{\mathcal{Y}}(L_{W}(0)w;
\tau)\nn
&\quad =F^{\phi}_{\mathcal{Y}}(L_{W}(-2)w;
\tau) -\sum_{k\in \mathbb{Z}_{+}}G_{2k+2}(\tau)
F^{\phi}_{\mathcal{Y}}(L_{W}(2k)w;
\tau),
\end{align}
where in the left-hand side, we have used 
$$F^{\phi}_{\mathcal{Y}}(L_{W}(0)w;
\tau)=(\wt w)F^{\phi}_{\mathcal{Y}}(w;
\tau)+F^{\phi}_{\mathcal{Y}}(L_{W}(0)_{N}w;
\tau).$$
Using (\ref{mod-inv-0}), we see that (\ref{mod-inv-1})
becomes
\begin{align}\label{mod-inv-2}
2\pi i\frac{\partial}{\partial \tau}
F^{\phi}_{\mathcal{Y}}(w;
\tau)
&=F^{\phi}_{\mathcal{Y}}(L_{W}(-2)w;
\tau) -\sum_{k\in \mathbb{N}}G_{2k+2}(\tau)
F^{\phi}_{\mathcal{Y}}(L_{W}(2k)w;
\tau)\nn
&=F^{\phi}_{\mathcal{Y}}\left(L_{W}(-2) -\sum_{k\in \mathbb{N}}G_{2k+2}(\tau)
L_{W}(2k)w;
\tau\right).
\end{align}
Let $\tau'=\frac{\alpha \tau+\beta}{\gamma \tau+\delta}$.
Then (\ref{mod-inv-2}) holds with $\tau$  replaced by $\tau'$. 
Using the modular transformation properties of $G_{2k+2}(\tau)$ for $k\in \N$ 
and the commutator formulas between $L_{W}(0)$ and
$L_{W}(n)$ 
for $n\in \Z$, we see that 
the formula obtained from (\ref{mod-inv-2}) with $\tau$ and 
$w$ replaced by $\tau'$ and $(\gamma\tau+\delta)^{-L_{W}(0)}w$
is equivalent to
\begin{align}\label{mod-inv-3}
&2\pi i\frac{\partial}{\partial \tau}
F^{\phi}_{\mathcal{Y}}((\gamma\tau+\delta)^{-L_{W}(0)}w;
\tau')\nn
&\quad=F^{\phi}_{\mathcal{Y}}((\gamma\tau+\delta)^{-L_{W}(0)}L_{W}(-2)w;
\tau') \nn
&\quad\quad-\sum_{k\in \mathbb{N}}G_{2k+2}(\tau)
F^{\phi}_{\mathcal{Y}}((\gamma\tau+\delta)^{-L_{W}(0)}L_{W}(2k)w;
\tau')\nn
&\quad=F^{\phi}_{\mathcal{Y}}\left((\gamma\tau+\delta)^{-L_{W}(0)}
\left(L_{W}(-2) -\sum_{k\in \mathbb{N}}G_{2k+2}(\tau)
L_{W}(2k)\right)w;
\tau'\right).
\end{align}

For $u\in V$, from the $n=1$ case of the identity (1.10) in \cite{F} and 
the fact that $L_{W}(0)_{N}$ commutes with vertex operators 
for $W$, we obtain 
\begin{equation}\label{mod-inv-4}
F^{\phi}_{\mathcal{Y}}((\gamma\tau+\delta)^{-L_{W}(0)}
\res_{x}Y_{W}(u, x)w;
\tau')=0
\end{equation}
and from the $n=1$, $l=2$ case of the identity 
(1.14) in \cite{F}, 
the modular transformation property of  $G_{2k+2}(\tau)$ for $k\in \N$
and the commutator formula between $L_{W}(0)$ and 
vertex operators 
for $W$,
we obtain
\begin{align}\label{mod-inv-5}
F^{\phi}_{\mathcal{Y}}
\left((\gamma\tau+\delta)^{-L_{W}(0)}\left(u_{-2}
 +\sum_{k\in \mathbb{Z}_{+}}(2k+1)G_{2k+2}(\tau)
u_{2k}\right)w;
\tau'\right)=0.
\end{align}

Formulas (7.9) and (7.10) in \cite{H-mod-inv-int} still hold, that is, we have 
\begin{equation}\label{mod-inv-6}
\left(L_{W}(-2)-\sum_{k\in \mathbb{N}}G_{2k+2}(\tau)
L_{W_{2}}(2k)\right)w=\res_{x}(\wp_{1}(x; \tau)-G_{2}(\tau)x)
Y_{W}(\omega, x)w
\end{equation}
and 
\begin{equation}\label{mod-inv-7}
\left(u_{-2}+\sum_{k\in \mathbb{Z}_{+}}(2k+1)G_{2k+2}(\tau)
u_{2k}\right)w=\res_{x}\wp_{2}(x; \tau)Y_{W}(u, x)w.
\end{equation}

Using (\ref{mod-inv-6}) and (\ref{mod-inv-7}), we see that (\ref{mod-inv-3}) and
(\ref{mod-inv-5}) become
\begin{align}\label{mod-inv-8}
&2\pi i\frac{\partial}{\partial \tau}
F^{\phi}_{\mathcal{Y}}((\gamma\tau+\delta)^{-L_{W}(0)}w;
\tau')\nn
&\quad=F^{\phi}_{\mathcal{Y}}\left((\gamma\tau+\delta)^{-L_{W}(0)}
\res_{x}(\wp_{1}(x; \tau)-G_{2}(\tau)x)Y_{W}(\omega, x)w;
\tau'\right)
\end{align}
and 
\begin{align}\label{mod-inv-9}
F^{\phi}_{\mathcal{Y}}
\left((\gamma\tau+\delta)^{-L_{W}(0)}\res_{x}\wp_{2}(x; \tau)
Y_{W}(u, x)w;
\tau'\right)=0,
\end{align}
respectively. 

Using analytic extensions, we see that (\ref{mod-inv-4}), (\ref{mod-inv-8})
and (\ref{mod-inv-9}) becomes
\begin{align*}
&\overline{F}^{\phi}_{\mathcal{Y}}
((\gamma\tau+\delta)^{-L_{W}(0)}\res_{x}Y_{W}(u, x)w;
\tau')=0,%\label{mod-inv-4-1}
\\
&2\pi i\frac{\partial}{\partial \tau}
\overline{F}
^{\phi}_{\mathcal{Y}}((\gamma\tau+\delta)^{-L_{W}(0)}w;
\tau')\nn
&\quad=\overline{F}^{\phi}_{\mathcal{Y}}
\left((\gamma\tau+\delta)^{-L_{W}(0)}
\res_{x}(\wp_{1}(x; \tau)-G_{2}(\tau)x)Y_{W}(\omega, x)w;
\tau'\right),%\label{mod-inv-8-1}
\\
&\overline{F}^{\phi}_{\mathcal{Y}}
\left((\gamma\tau+\delta)^{-L_{W}(0)}\res_{x}\wp_{2}(x; \tau)
Y_{W}(u, x)w;
\tau'\right)=0.
\end{align*}
These are exactly the conditions for the linear map given by 
(\ref{mod-transf-conf-bk-1}) to be a genus-one $1$-point 
conformal block. 
\epfv

As in \cite{H-mod-inv-int}, Theorem \ref{F-genus1-conv} 
in \cite{F1} and \cite{F} is proved by deriving 
a system of differential equations of regular singular points satisfied by 
the formal series of shifted pseudo-$q_{\tau}$-traces of products of 
geometrically-modified intertwining operators. In particular, the
genus-one $1$-point correlation functions
$\overline{F}^{\phi}_{\Y}(L_{W}(0)_{N}^{p}w; \tau)$ 
for $k=0, \dots, K$ satisfy a system of 
differential equations of a regular singular point $q_{\tau}=0$, where 
$K$ is the smallest number of $\N$
such that $L_{W}(0)_{N}^{K+1}w=0$. 
In fact, the system of differential 
equations satisfied by $\overline{F}^{\phi}_{\Y}(L_{W}(0)_{N}^{p}w; \tau)$
are derived using only the $C_{2}$-cofiniteness and the 
$q_{\tau}$-expansions of 
(\ref{1-pt-cb-1}), (\ref{1-pt-cb-2}) and (\ref{1-pt-cb-3}).
Using this fact, we first show that elements in the image of 
a genus-one $1$-point 
conformal block satisfy differential equations of a regular singular point
$q_{\tau}=0$.

\begin{prop}\label{diff-conf-bk}
Let $V$ be a $C_{2}$-cofinite vertex operator algebra and 
 $W$ a grading-restricted generalized $V$-module. 
For a  genus-one $1$-point 
conformal block $F$ labeled by $W$ and $w\in W$,
$F(L_{W}(0)_{N}^{j}w; \tau)$ for $j=0, \dots, K$
satisfy a system of $K+1$ differential equations of a
regular singular point at $q_{\tau}=0$, where 
$K$ is the smallest number of $\N$ such that $L_{W}(0)^{K+1}w=0$, 
\end{prop}
\pf
As in the case of $n=1$ in \cite{H-mod-inv-int},
\cite{F1} and \cite{F}, let 
$R=\C[G_{4}(\tau), G_{6}(\tau)]$. Then $G_{2k+2}(\tau)\in R$
for $k\in \Z_{+}$. Consider the 
$R$-module $M_{F}$ generated by functions of $\tau$ of the form 
$F(w; \tau)$. Then the linear map $F$ in fact induces an $R$-module map
$\widehat{F}$ from the $R$-module $T=R\otimes W$ to 
$M_{F}$ given by $\widehat{F}(f(\tau)\otimes w)=f(\tau) F(w; \tau)$. 

As in \cite{H-mod-inv-int}, \cite{F1} and \cite{F}, let
$J$ be an $R$-submodule of $T$ generated by elements of the form
$$(Y_{W})_{-2}(v)w+\sum_{k\in \Z_{+}}(2k+1)G_{2k+2}(\tau)
(Y_{W})_{2k}(v)w$$
for $v\in V$ and $w\in W$. 
Then by (\ref{mod-inv-7}) and (\ref{1-pt-cb-3}), we see that 
$\ker \widehat{F}\subset J$. In particular, 
$\widehat{F}$ can in fact be viewed as an $R$-module map
from $T/J$ to $M_{F}$. 

Since $V$ is $C_{2}$-cofinite, 
the same proof as those in \cite{H-mod-inv-int}, \cite{F1} and \cite{F}
shows that $T/J$ is a finitely generated $R$-module.
As in \cite{H-mod-inv-int}, \cite{F1} and \cite{F}, let 
$\mathcal{Q}: T\to T$ be the map defined by
\begin{align*}
\mathcal{Q}(f(\tau)\otimes w)&=
f(\tau)\otimes (\res_{x}\wp_{1}(x; \tau)
Y_{W}(\omega, x)w).
\end{align*}
Since $R$ is Noetherian, given $w\in W$, the $R$-submodule 
of $T$ generated by $\mathcal{Q}^{n}(1\otimes w)$ for $n\in \N$
is also finitely generated. Then there exist $m\in \Z_{+}$ and 
$b_{p}(\tau)\in R$ for $p=1, \dots, m$
such that 
$$\mathcal{Q}^{m}(1\otimes w)+\sum_{p=1}^{m}b_{p}(\tau)
\mathcal{Q}^{m-p}(1\otimes w)\in J.$$
Then by the definition of $\mathcal{Q}$, we obtain
\begin{equation}\label{diff-conf-bk-1}
1\otimes (\res_{x}\wp_{1}(x; \tau)
Y_{W}(\omega, x))^{m}w+\sum_{p=1}^{m}b_{p}(\tau)
\otimes (\res_{x}\wp_{1}(x; \tau)
Y_{W}(\omega, x))^{m-p}w)\in J.
\end{equation}
For homogeneous $w\in W$, from (\ref{1-pt-cb-2}), we have 
\begin{align*}
&\widehat{F}(1\otimes (\res_{x}\wp_{1}(x; \tau)
Y_{W}(\omega, x))w)\nn
&\quad=F(\res_{x}\wp_{1}(x; \tau)
Y_{W}(\omega, x))w)\nn
&\quad=2\pi i\frac{d}{d\tau}F(w)+
G_{2}(\tau)F(L_{W}(0)w)\nn
&\quad=\left(2\pi i\frac{d}{d\tau}+
 (\wt w)G_{2}(\tau)\right)F(w)
+ G_{2}(\tau)F(L_{W}(0)_{N}w).
\end{align*}
Then for $m\in \N$, 
\begin{align*}
&\widehat{F}(1\otimes (\res_{x}\wp_{1}(x; \tau)
Y_{W}(\omega, x))^{m}w)\nn
&\quad =\left(2\pi i\frac{d}{d\tau}+
 (\wt w)G_{2}(\tau)\right)^{m}F(w)+
\sum_{j=1}^{m}D_{m, j}\left(\frac{d}{d\tau}, \tau\right)
F(L_{W}(0)_{N}^{j}w),
\end{align*}
where for $j=1, \dots, m$, $D_{m, j}(\frac{d}{d\tau}, \tau)$
is a polynomial in $\frac{d}{d\tau}$ of degree less than $m$ with 
polynomials in $G_{2}(\tau)$ and its derivatives as coefficients. 

Applying $\widehat{F}$ to the element of $J$ in (\ref{diff-conf-bk-1})
and using $\ker \hat{F}\subset J$, 
we obtain 
\begin{align}\label{diff-1}
&\left(2\pi i\frac{d}{d\tau}+
 (\wt w)G_{2}(\tau)\right)^{m}F(w)
+\sum_{p=1}^{m}b_{p}(\tau)\left(2\pi i\frac{d}{d\tau}+
 (\wt w)G_{2}(\tau)\right)^{m-p}F(w)\nn
 &\quad\quad +\sum_{j=1}^{m}D_{m, j}\left(\frac{d}{d\tau}, \tau\right)
F(L_{W}(0)_{N}^{j}w) +\sum_{p=1}^{m}\sum_{j=1}^{m-p}
b_{p}(\tau)D_{m-p, j}\left(\frac{d}{d\tau}, \tau\right)
F(L_{W}(0)_{N}^{j}w) 
\nn
 &\quad =0.
\end{align}

Note that in (\ref{diff-1}), $b_{p}(\tau)$ and 
$D_{m-p, j}\left(\frac{d}{d\tau}, \tau\right)$ for $p=1, \dots, m$ and 
$j=0, \dots, m-p$ are independent of the genus-one $1$-point 
conformal block $F$.
Since $L_{W}(0)_{N}$ commutes with vertex operators acting 
on $W$, we see that the linear map given by 
$w\mapsto F(L_{W}(0)^{j}w)$ for $w\in W$ is also a 
genus-one $1$-point 
conformal block. Then for $j=1, \dots, K$, 
$F(L_{W}(0)^{j}w)$ also satisfy 
the equation (\ref{diff-1}). So we see that 
$F(L_{W}(0)^{j}w)$ for $j=0, \dots, K$ satisfy a system 
of $K+1$ differential equations. Moreover, since 
$\frac{d}{d\tau}=2\pi i q_{\tau}\frac{d}{dq_{\tau}}$, 
The singular point $q_{\tau}=0$ of this system of differential equations 
is regular. 
\epfv

Recall the genus-one $1$-point conormal block 
$\overline{F}^{\phi}_{\Y}(w; \tau)$ obtained by taking 
shifted pseudo-$q_{\tau}$-traces of an intertwining operator 
$\Y$ and then extending it analytically. 
For $w\in W$,
let $\mathcal{F}_{w}$ be the space spanned by 
$\overline{F}^{\phi}_{\Y}(w; \tau)$
for all finite-dimensional 
associative algebras $P$ with a symmetric linear functions $\phi$, 
all grading-restricted generalized $V$-$P$-bimodules
$\widetilde{W}$, projective as a right $P$-module, and
all intertwining operators $\Y$ of type
$\binom{\widetilde{W}}{W\widetilde{W}}$
compatible with $P$. 

We now prove the following main result of this section,
which together with Proposition \ref{mod-transf-conf-bk}
implies 
the modular invariance theorem (Theorem \ref{mod-inv}):

\begin{thm}\label{conf-bk=>pseudo-q-tr}
Let $V$ be a $C_{2}$-cofinite vertex operator algebra  
without nonzero elements of negative weights,
$W$ a grading-restricted generalized $V$-module
and $F: W\to H(\H)$ a genus-one $1$-point conformal block labeled by
$W$. Then for $w\in W$, $F(w; \tau)$ is in $\mathcal{F}_{w}$. 
\end{thm} 
\pf
By Proposition \ref{diff-conf-bk}, 
$F(L_{W}(0)^{j}w)$ for $j=0, \dots, K$ satisfy a system 
of $K+1$ differential equations and the singular point $q_{\tau}=0$ 
is regular.
Using the theory of differential equations of regular singular points, we have
an expansion
\begin{equation}\label{mod-inv-10}
F(L_{W}(0)_{N}^{p}w; \tau)
=\sum_{k=0}^{K}\sum_{j=1}^{J}\sum_{m\in \mathbb{N}}
C_{k, j, m}^{p}(w)\tau^{k}q_{\tau}^{r_{j}+m},
\end{equation}
where $r_{j}$ for $j=1, \dots, J$ are complex numbers such that 
$r_{j_{1}}-r_{j_{2}}\not \in \mathbb{Z}$ when $j_{1}\ne j_{2}$.

From (\ref{mod-inv-10}), the properties of the 
genus-one $1$-point conformal block $F$ and Lemma \ref{tilde-wp-2-and-1}, 
we obtain
\begin{align}
&\sum_{k=0}^{K}\sum_{j=1}^{J}
\sum_{m\in \N}C_{k, j, m}^{0}(\res_{x}Y_{W}(u, x)w)
(\log q)^{k}q^{r_{j}+m}=0,
\label{mod-inv-4-2}\\
&2\pi i\sum_{k=0}^{K}\sum_{j=1}^{J}
\sum_{m\in \N}(k+1)C_{k+1, j, m}^{0}(w)
(\log q)^{k}q^{r_{j}+m}
+(2\pi i)^{2}\sum_{m\in \N}(r_{j}+m)C_{k, j, m}^{0}(w)
(\log q)^{k}q^{r_{j}+m}\nn
&\quad=\sum_{m\in \N}C_{k, j, m}^{0}(\res_{x}
(\tilde{\wp}_{1}(x; q)-\widetilde{G}_{2}(q)x)
Y_{W}(\omega, x)w)(\log q)^{k}q^{r_{j}+m},
\label{mod-inv-8-2}\\
&\sum_{k=0}^{K}\sum_{j=1}^{J}\sum_{m\in \N}C_{k, j, m}^{0}(\res_{x}\tilde{\wp}_{2}(x; q)
Y_{W}(u, x)w)(\log q)^{k}q^{r_{j}+m}=0.\label{mod-inv-9-2}
\end{align}
From (\ref{mod-inv-4-2})--(\ref{mod-inv-9-2}), we see that 
the linear map 
$S: W\to \C\{q\}[\log q]$ given by 
$$S(w)=\sum_{k=0}^{K}\sum_{j=1}^{J}
\sum_{m\in \N}C_{k, j, m}^{0}(w)q^{m}$$
satisfies the conditions needed in Theorem \ref{wp-fn-gr-sym-fn}. 

Since $V$ is $C_{2}$-cofinite, there are only finitely many inequivalent 
irreducible $V$-modules. Let $N$ be a nonnegative integer larger than 
all the real parts of the differences between the finitely many 
lowest weights of the irreducible $V$-modules. Then 
any $N'\in N+\N$ is also larger than 
all the real parts of the differences between the finitely many 
lowest weights of the irreducible $V$-modules. 
By Theorem \ref{wp-fn-gr-sym-fn}, for $k=0, \dots, K$ and $
j=1, \dots, J$, the linear map 
$\psi_{S; k, j}^{N'}: U^{N'}(W)\to \C$ defined by 
$$\psi_{S; k, j}^{N'}([w]_{mn})=0$$
for $0\le m, n\le N'$, $m\ne n$, and $w\in W$
and 
$$\psi_{S; k, j}^{N'}([w]_{mm})=C_{k, j, m}^{0}(w)$$
for $0\le m\le N'$ and $w\in W$
induces a symmetric linear function, still denoted by
$\psi_{S; k, j}^{N'}$, on $\widetilde{A}^{N'}(W)$ satisfying 
\begin{equation}\label{psi-grading}
\psi_{S; k, j}^{N'}(([\omega]_{mm}
-(2\pi i)^{2}(r_{j}+m)[\one]_{mm})^{\bkdiamd (K-k+1)}
\bkdiamd \mathfrak{w})=0.
\end{equation}
By Proposition \ref{finite-d-til-A-N-W}, we know that 
$\widetilde{A}^{N'}(W)$
is finite dimensional. 

Note that 
$$\one^{N'}+\widetilde{Q}^{\infty}(V)
=\sum_{k=0}^{N'}[\one]_{kk}+\widetilde{Q}^{\infty}(V)$$
is the identity of $\widetilde{A}^{N'}(V)$.
Let 
$$\one^{N'}+\widetilde{Q}^{\infty}(V)
=\tilde{e}^{N'}_{1}+\cdots +\tilde{e}^{N'}_{n'},$$
where 
$\tilde{e}_{1}^{N'}, \dots, \tilde{e}^{N'}_{n'}\in \widetilde{A}^{N'}(V)$ 
are orthogonal primitive central idempotents in 
$\widetilde{A}^{N'}(V)$. 
Then $\widetilde{A}^{N'}(V)=\widetilde{A}^{N'}_{1}
\oplus\cdots \oplus \widetilde{A}^{N'}_{n'}$
where $\widetilde{A}^{N'}_{i}=\widetilde{A}^{N'}(V)\bkdiamd 
\tilde{e}^{N'}_{i}$ for 
$i=1, \dots, n'$, is a decomposition of 
$\widetilde{A}^{N'}(V)$ into a direct sum of indecomposable 
$\widetilde{A}^{N'}(V)$-bimodules. 
Let $\widetilde{U}^{N'}_{i}=\tilde{e}^{N'}_{i}\bkdiamd 
\widetilde{A}^{N'}(W)
\bkdiamd \tilde{e}^{N'}_{i}$ 
for $i=1, \dots, n'$. Then we have the decomposition 
$$\widetilde{A}^{N'}(W)=\widetilde{U}^{N'}_{1}\oplus \cdots 
\oplus \widetilde{U}^{N'}_{n'}$$
of $\widetilde{A}^{N'}(W)$ as a direct sum
of $\widetilde{A}_{1}^{N'}$-, $\dots$, 
$\widetilde{A}_{n'}^{N'}$-bimodules. 
Let $\widetilde{B}^{N'}=\widetilde{A}^{N'}(V)
\oplus \widetilde{A}^{N'}(W)$ be the trivial 
square-zero extension of $\widetilde{A}^{N'}(V)$ 
by $\widetilde{A}^{N'}(W)$ and 
$\widetilde{B}_{i}^{N'}=\widetilde{A}^{N'}_{i}\oplus 
\widetilde{U}^{N'}_{i}$ for $i=1, \dots, n'$ 
the trivial 
square-zero extension of $\widetilde{A}^{N'}_{i}$ by 
$\widetilde{U}^{N'}_{i}$. Then 
$$\widetilde{B}^{N'}=\widetilde{B}_{1}^{N'}\oplus \cdots\oplus 
\widetilde{B}_{n'}^{N'}.$$

For $i=1, \dots, n'$, let $\bar{\phi}^i_{k, j; N'}
=\psi_{S; k, j}^{N'}|_{\widetilde{B}^{N'}_{i}}$ and let 
\begin{equation}\label{idempodents}
\{\tilde{e}^{N'}_{i, \mu\nu}\mid \mu=1, \dots, l'_{i},\; 
\nu=1, \dots, l'_{i\mu}\}
\end{equation}
be a complete 
set of orthogonal primitive idempodents of $\widetilde{A}^{N'}_{i}$
such that 
$$1_{\widetilde{A}^{N'}_{i}}=\tilde{e}
^{N'}_{i}=\sum_{\mu=1}^{l'_{i}}\sum_{\nu=1}^{l'_{i\mu}}\tilde{e}
^{N'}_{i, \mu\nu},$$
$\tilde{e}^{N'}_{i, \mu\nu_{1}}\widetilde{A}^{N'}_{i}
\simeq \tilde{e}^{N'}_{i, \mu\nu_{2}}
\widetilde{A}^{N'}_{i}$
for $\mu=1, \dots, l'_{i}$,
$\nu_{1}, \nu_{2}=1, \dots, l'_{i\mu}$
and 
$\tilde{e}^{N'}_{i, \mu_{1}\nu_{1}}\widetilde{A}^{N'}_{i}
\not\simeq \tilde{e}^{N'}_{i, \mu_{2}\nu_{2}}
\widetilde{A}^{N'}_{i}$
for  $1\le \mu_{1}<\mu_{2}\le l_{i}'$, $\nu_{1}=1, \dots, l'_{i\mu_{1}}$,
$\nu_{2}=1, \dots, l'_{i\mu_{2}}$. Then \eqref{idempodents}
is also a complete set of orthogonal primitive idempodents 
of $\widetilde{B}^{N'}_{i}$ such that 
$$1_{\widetilde{B}^{N'}_{i}}=1_{\widetilde{A}^{N'}_{i}}
=\sum_{\mu=1}^{l'_{i}}\sum_{\nu=1}^{l'_{i\mu}}\tilde{e}
^{N'}_{i, \mu\nu},$$
$\tilde{e}^{N'}_{i, \mu\nu_{1}}\widetilde{B}^{N'}_{i}
\simeq \tilde{e}^{N'}_{i, \mu\nu_{2}}
\widetilde{B}^{N'}_{i}$
for $\mu=1, \dots, l'_{i}$,
$\nu_{1}, \nu_{2}=1, \dots, l'_{i\mu}$
and 
$\tilde{e}^{N'}_{i, \mu_{1}\nu_{1}}\widetilde{B}^{N'}_{i}
\not\simeq \tilde{e}^{N'}_{i, \mu_{2}\nu_{2}}
\widetilde{B}^{N'}_{i}$
for  $1\le \mu_{1}<\mu_{2}\le l_{i}'$, $\nu_{1}=1, \dots, l'_{i\mu_{1}}$,
$\nu_{2}=1, \dots, l'_{i\mu_{2}}$.
Let 
$$\tilde{\epsilon}^{N'}_{i, \mu\nu}=\tilde{e}^{N'}_{i, \mu\nu}
+\text{Rad}(\phi^i_{k, j; N'})\in \widetilde{B}^{N'}_{i}
/\text{Rad}(\phi^i_{k, j; N'})=\widetilde{B}^{N'}
/\text{Rad}(\phi^i_{k, j; N'})$$ 
for $\mu=1, \dots, l'_{i}$,
$\nu=1, \dots, l'_{i\mu}$. Then 
$$1_{\widetilde{B}^{N'}/\text{Rad}(\phi^i_{k, j; N'})}
=\sum_{\mu=1}^{l'_{i}}\sum_{\nu=1}^{l'_{i\mu}}\tilde{\epsilon}
^{N'}_{i, \mu\nu},$$
$$\tilde{\epsilon}^{N'}_{i, \mu\nu_{1}}(\widetilde{B}^{N'}
/\text{Rad}(\phi^i_{k, j; N'}))
\simeq \tilde{\epsilon}^{N'}_{i, \mu\nu_{2}}
(\widetilde{B}^{N'}
/\text{Rad}(\phi^i_{k, j; N'}))$$
for $\mu=1, \dots, l'_{i}$,
$\nu_{1}, \nu_{2}=1, \dots, l'_{i\mu}$
and 
$$\tilde{\epsilon}^{N'}_{i, \mu_{1}\nu_{1}}(\widetilde{B}^{N'}
/\text{Rad}(\phi^i_{k, j; N'}))
\not\simeq 
\tilde{\epsilon}^{N'}_{i, \mu_{2}\nu_{2}}
(\widetilde{B}^{N'}
/\text{Rad}(\phi^i_{k, j; N'}))$$
for  $1\le \mu_{1}<\mu_{2}\le l_{i}'$, $\nu_{1}=1, \dots, l'_{i\mu_{1}}$,
$\nu_{2}=1, \dots, l'_{i\mu_{2}}$. 
Let
$$\tilde{\epsilon}^{N'}_{i}=\sum_{\mu=1}^{l'_{i}}
\tilde{\epsilon}^{N'}_{i, \mu 1}.$$
Then by Theorem \ref{slf=>p-tr}, 
$$P^i_{k, j; N'}=\tilde{\epsilon}_{i}^{N'}(\widetilde{B}^{N'}
/\text{Rad}(\phi^i_{k, j; N'}))
\tilde{\epsilon}_{i}^{N'}$$ 
for $i=1, \dots, n'$
are basic symmetric algebras 
equipped with symmetric linear functions given by $\phi^i_{k, j; N'}$, 
and 
$$M^i_{k, j; N'}=(\widetilde{B}^{N'}/\text{Rad}(\phi^i_{k, j; N'}))
\tilde{\epsilon}_{i}^{N'}$$
are 
$\widetilde{A}^{N'}(V)$-$P^i_{k, j; N'}$-bimodules  
which are finitely generated and projective as right $P^i_{k, j; N'}$-modules.
Moreover, we define 
$$f^i_{k, j; N'}\in \hom_{\widetilde{A}^{N'}(V), P^{i}_{k, j; N'}}
(\widetilde{A}^{N'}(W)\otimes_{\widetilde{A}^{N'}(V)} 
M^i_{k, j; N'}, M^i_{k, j; N'})$$
for $i=1, \dots, n'$ by
$f^i_{k, j; N'}(\mathfrak{w}\otimes w_{i})
=(0, \mathfrak{w})w_{i}$ for $\mathfrak{w}\in 
\widetilde{A}^{N'}(W)$ and $w_{i}\in M^i_{k, j; N'}$
and then we have 
\begin{equation}\label{mod-inv-13}
\psi_{S; k, j}^{N'}(\mathfrak{w}) 
= \sum_{i=1}^{n'} (\phi^{i}_{k, j; N'})^{f^i_{k, j; N'}}_{M^i_{k, j; N'}}
(\mathfrak{w}
+\widetilde{Q}^{\infty}(W))
\end{equation}
for $\mathfrak{w}\in A^{N'}(W)$.
We use $\tilde{\vartheta}_{M^i_{k, j; N'}}: \widetilde{A}^{N'}(V)
\to (\text{End}\; M^i_{k, j; N'})$ to denote
the homomorphism of associative algebras giving the 
$\widetilde{A}^{N'}(V)$-module structure on 
$M^i_{k, j; N'}$. 
Then from (\ref{psi-grading}) and 
Theorem \ref{slf=>p-tr}, we have
\begin{equation}\label{generalized-eigenspace}
\tilde{\vartheta}_{M^i_{k, j; N'}}(([\omega]_{mm}
-(2\pi i)^{2}(r_{j}+m)[\one]_{mm})^{\bkdiamd (K-k+1)}
+\widetilde{Q}^{\infty}(V))
M^i_{k, j; N'}=0
\end{equation}
for $m=0, \dots, N'$.

From Proposition \ref{tilde-A-infty}, 
the map $\mathcal{U}(1)^{-1}$ induces
an isomorphism from $A^{\infty}(V)$ to $\widetilde{A}^{\infty}(V)$ and 
thus also induces an isomorphism from 
$A^{N'}(V)$ to $\widetilde{A}^{N'}(V)$. In particular, 
the map $\vartheta_{M^i_{k, j; N'}}: A^{N'}(V)\to
(\text{End} \; M^i_{k, j; N'})$ defined by 
$$\vartheta_{M^i_{k, j; N'}}([v]_{kl}+Q^{\infty}(V))
=\tilde{\vartheta}_{M^i_{k, j; N'}}([\mathcal{U}(1)^{-1}v]_{kl}
+\widetilde{Q}^{\infty}(V))$$
for $v\in V$ and $k, l=0, \dots, N'$ gives an $A^{N'}(V)$-module 
structure to $M^i_{k, j; N'}$. 

We now show that 
$M^{i}_{k, j; N'}$ for $i=1, \dots, n'$, $k=0, \dots, K$ and
$j=1, \dots, J$  are 
graded $A^{N'}(V)$-modules (see Definition 5.1 in \cite{H-aa-va}). 

Note that 
$$\sum_{m=0}^{N'}\vartheta_{M^i_{k, j; N'}}
([\one]_{mm}+Q^{\infty}(V))=1_{M^i_{k, j; N'}}.$$
Also we have 
\begin{align*}
([\one]_{mm}+Q^{\infty}(V))
\diamond ([\one]_{mm}+Q^{\infty}(V))
&=[\one]_{mm}+Q^{\infty}(V),\\
([\one]_{ll}+Q^{\infty}(V))\diamond ([\one]_{mm}+Q^{\infty}(V))&=0
\end{align*}
for $l, m=0, \dots, N'$, $l\ne m$. So
$$\{\vartheta_{M^i_{k, j; N'}}([\one]_{mm}
+Q^{\infty}(V))\}_{m=0}^{N'}$$
is a partition of the identity on $M^i_{k, j; N'}$. In particular, 
$$M^i_{k, j; N'}=\coprod_{m=0}^{N'}(M^i_{k, j; N'})_{[m]},$$
where for $m=0, \dots, N'$, 
$$(M^i_{k, j; N'})_{[m]}=\vartheta_{M^i_{k, j; N'}}([\one]_{mm}
+Q^{\infty}(V))M^i_{k, j; N'}.$$ 
Since $[\omega]_{mm}+Q^{\infty}(V)$ 
commutes with $[\one]_{mm}+Q^{\infty}(V)$
for $m=0, \dots, N'$, $(M^i_{k, j; N'})_{[m]}$ is invariant under
$\vartheta_{M^i_{k, j; N'}}
([\omega]_{mm}+Q^{\infty}(V))$. Since 
$$([\omega]_{mm}+Q^{\infty}(V))\diamond ([\one]_{ll}+Q^{\infty}(V))
=([\omega]_{ll}+Q^{\infty}(V))\diamond 
([\one]_{mm}+Q^{\infty}(V))=0$$ 
for $l\in \N$ not equal to $m$, 
we see that $\vartheta_{M^i_{k, j; N'}}([\omega]_{mm}
+Q^{\infty}(V))$ is $0$ on $(M^i_{k, j; N'})_{[l]}$ . 
Note also that
$\vartheta_{M^i_{k, j; N'}}([\one]_{mm}+Q^{\infty}(V))$ 
on $(M^i_{k, j; N'})_{[m]}$  is the identity on
$(M^i_{k, j; N'})_{[m]}$. From (\ref{generalized-eigenspace}),
we have
\begin{equation}\label{tilde-weight}
\tilde{\vartheta}_{M^i_{k, j; N'}}(([\omega]_{mm}
-(2\pi i)^{2}(r_{j}+m)[\one]_{mm}
+\widetilde{Q}^{\infty}(V))^{\bkdiamd (K-k+1)})(M^i_{k, j; N'})_{[m]}=0.
\end{equation}
By the definition of $\vartheta_{M^i_{k, j; N'}}$ and (\ref{tilde-weight}),
we obtain 
\begin{equation}\label{weight-0}
\vartheta_{M^i_{k, j; N'}}(([\mathcal{U}(1)\omega]_{mm}
-(2\pi i)^{2}(r_{j}+m)[\mathcal{U}(1)\one]_{mm}
+Q^{\infty}(V))^{\diamond \;(K-k+1)})(M^i_{k, j; N'})_{[m]}=0.
\end{equation}
Using $\mathcal{U}(1)\omega=(2\pi i)^{2}(\omega-\frac{c}{24}\one)$ and 
$\mathcal{U}(1)\one=\one$ (see the definition of $\mathcal{U}(1)$ and
Lemma 1.1 in \cite{H-mod-inv-int}), we see that 
(\ref{weight-0}) becomes 
\begin{equation}\label{weight}
\vartheta_{M^i_{k, j; N'}}\left(\left([\omega]_{mm}
-\left(r_{j}+\frac{c}{24}+m\right)[\one]_{mm}
+Q^{\infty}(V)\right)^{\diamond \;(K-k+1)}\right)(M^i_{k, j; N'})_{[m]}=0.
\end{equation}
Thus $(M^i_{k, j; N'})_{[m]}$ is the generalized eigenspace of 
$\vartheta_{M^i_{k, j; N'}}([\omega]_{mm}+Q^{\infty}(V))$ 
with eigenvalue $r_{j}+\frac{c}{24}+m$. 

For $v\in V$, $k, l=0, \dots, N'$ and $w\in (M^i_{k, j; N'})_{[m]}$, 
\begin{align*}
\vartheta_{M^i_{k, j; N'}}([v]_{kl}+Q^{\infty}(V))w
&=\vartheta_{M^i_{k, j; N'}}([v]_{kl}+Q^{\infty}(V))
\vartheta_{M^i_{k, j; N'}}([\one]_{mm}+Q^{\infty}(V))w\nn
&=\vartheta_{M^i_{k, j; N'}}([v]_{kl}\diamond [\one]_{mm}+Q^{\infty}(V))w\nn
&=\delta_{lm}\vartheta_{M^i_{k, j; N'}}([v]_{km}+Q^{\infty}(V))w\nn
&=\delta_{lm}\vartheta_{M^i_{k, j; N'}}([\one]_{kk}
\diamond [v]_{km}+Q^{\infty}(V))w\nn
&=\delta_{lm}\vartheta_{M^i_{k, j; N'}}([\one]_{kk}+Q^{\infty}(V))
\vartheta_{M^i_{k, j; N'}}([v]_{km}+Q^{\infty}(V))w,
\end{align*}
which is $0$ when $l\ne m$ and is in $(M^i_{k, j; N'})_{[k]}$ when $l=m$.
So Condition 1 in Definition 5.1 in \cite{H-aa-va} is satisfied. 

We define  
\begin{align*}
L_{M^i_{k, j; N'}}(0)&=\sum_{m=0}^{N'}\vartheta_{M^i_{k, j; N'}}
([\omega]_{mm}
+Q^{\infty}(V)),\\ 
L_{M^i_{k, j; N'}}(-1)&=\sum_{m=0}^{N'-1}\vartheta_{M^i_{k, j; N'}}
([\omega]_{m+1, m}
+Q^{\infty}(V)).
\end{align*}
Then $(M^i_{k, j; N'})_{[m]}$ for $m=0, \dots, N'$ are generalized eigenspaces
for $L_{M^i_{k, j; N'}}(0)$ and $M^i_{k, j; N'}$ is the direct sum of 
these  generalized eigenspaces. The eigenvalues of $L_{M^i_{k, j; N'}}(0)$
are $r_{j}+\frac{c}{24}+m$ for $m=0, \dots, N'$ 
and the real parts of these eigenvalues has a minimum $\Re(r_{j})$. 
This shows that Condition 2 in Definition 5.1 in \cite{H-aa-va} is satisfied. 
From what we have shown above, we see 
that  $L_{M^i_{k, j; N'}}(-1)$ maps 
$(M^i_{k, j; N'})_{[m]}$ for $m=0, \dots, N'-1$
to $(M^i_{k, j; N'})_{[m+1]}$. So Condition 3 in 
Definition 5.1 in \cite{H-aa-va} is also satisfied. 

From Remark 4.5 in \cite{H-aa-int-op} with $W=V$, we have 
\begin{align}\label{L0-L-1-commutator-1}
[L_{M^i_{k, j; N'}}(0), L_{M^i_{k, j; N'}}(-1)]w
&=\vartheta_{M^i_{k, j; N'}}([\omega]_{m+1, m+1}
[\omega]_{m+1, m}-[\omega]_{m+1, m}[\omega]_{mm}
+Q^{\infty}(V))w\nn
&=\vartheta_{M^i_{k, j; N'}}([(L_{V}(-1)+L_{V}(0))\omega]_{m+1, m}
+Q^{\infty}(V))w
\end{align}
for $w\in (M^i_{k, j; N'})_{[m]}$. By Proposition 2.3 in \cite{H-aa-va},
we know that 
$$[(L_{V}(-1)+L_{V}(0)-1)\omega]_{m+1, m}\in Q^{\infty}(V).$$
Then the right-hand side of (\ref{L0-L-1-commutator-1}) 
is equal to 
\begin{equation}\label{L0-L-1-commutator-2}
\vartheta_{M^i_{k, j; N'}}([\omega]_{m+1, m})w
=L_{M^i_{k, j; N'}}(-1)w.
\end{equation}
From (\ref{L0-L-1-commutator-1})  and (\ref{L0-L-1-commutator-2}),
we obtain the commutator formula
\begin{equation}\label{L0-L-1-commutator-3}
[L_{M^i_{k, j; N'}}(0), L_{M^i_{k, j; N'}}(-1)]
=L_{M^i_{k, j; N'}}(-1).
\end{equation}

For $k, l=0, \dots, N'$ and $v\in \N$, by 
Remark 4.5 in \cite{H-aa-int-op} with $W=V$ and 
Proposition 2.3 in \cite{H-aa-va},
\begin{align*}
&[L_{M^i_{k, j; N'}}(0), \vartheta_{M^i_{k, j; N'}}([v]_{kl}+Q^{\infty}(V))]w\nn
&\quad =\vartheta_{M^i_{k, j; N'}}([\omega]_{kk}+Q^{\infty}(V))
\vartheta_{M^i_{k, j; N'}}([v]_{kl}+Q^{\infty}(V))w\nn
&\quad \quad -\vartheta_{M^i_{k, j; N'}}([v]_{kl}+Q^{\infty}(V))
\vartheta_{M^i_{k, j; N'}}([\omega]_{ll}+Q^{\infty}(V))w\nn
&\quad =\vartheta_{M^i_{k, j; N'}}([\omega]_{kk}
\diamond [v]_{kl}-[v]_{kl}\diamond [\omega]_{ll}+Q^{\infty}(V))w\nn
&\quad =\vartheta_{M^i_{k, j; N'}}([(L_{V}(-1)+L_{V}(0))v]_{kl}
+Q^{\infty}(V))w\nn
&\quad =(k-l)\vartheta_{M^i_{k, j; N'}}([v]_{kl}
+Q^{\infty}(V))w
\end{align*}
for $w\in (M^i_{k, j; N'})_{[l]}$. Then we obtain the commutator formula
\begin{equation}\label{L0-v-kl-commutator}
[L_{M^i_{k, j; N'}}(0), \vartheta_{M^i_{k, j; N'}}([v]_{kl}+Q^{\infty}(V))]
=(k-l)\vartheta_{M^i_{k, j; N'}}([v]_{kl}
+Q^{\infty}(V)).
\end{equation}

For $k=0, \dots, N'-1$, $l=1, \dots, N'$ and $v\in V$, 
also by 
Remark 4.5 in \cite{H-aa-int-op} with $W=V$,
\begin{align*}
&[L_{M^i_{k, j; N'}}(-1), \vartheta_{M^i_{k, j; N'}}([v]_{kl}+Q^{\infty}(V))]w\nn
&\quad =\vartheta_{M^i_{k, j; N'}}([\omega]_{k+1, k}+Q^{\infty}(V))
\vartheta_{M^i_{k, j; N'}}([v]_{kl}+Q^{\infty}(V))w\nn
&\quad \quad -\vartheta_{M^i_{k, j; N'}}([v]_{kl}+Q^{\infty}(V))
\vartheta_{M^i_{k, j; N'}}([\omega]_{l+1, l}+Q^{\infty}(V))w\nn
&\quad =\vartheta_{M^i_{k, j; N'}}([\omega]_{k+1, k}
\diamond [v]_{kl}-[v]_{kl}\diamond [\omega]_{l+1, l}+Q^{\infty}(V))w\nn
&\quad =\vartheta_{M^i_{k, j; N'}}([L_{V}(-1)v]_{kl}
+Q^{\infty}(V))w
\end{align*}
for $w\in (M^i_{k, j; N'})_{[l]}$. 
Then we obtain the commutator formula
\begin{equation}\label{L-1-v-kl-commutator}
[L_{M^i_{k, j; N'}}(-1), \vartheta_{M^i_{k, j; N'}}([v]_{kl}+Q^{\infty}(V))]
=\vartheta_{M^i_{k, j; N'}}([L_{V}(-1)v]_{kl}
+Q^{\infty}(V)).
\end{equation}
 
 From (\ref{L0-L-1-commutator-3}), (\ref{L0-v-kl-commutator})
 and (\ref{L-1-v-kl-commutator}), we see that 
 Condition 4 in Definition 5.1 in \cite{H-aa-va} is  satisfied. 
 Thus we have shown that $M^i_{k, j; N'}$ is indeed a 
 graded $A^{N'}(V)$-module. 

From Section 6 in \cite{H-aa-int-op}, 
we have the lower-bounded generalized $V$-module 
$W^i_{k, j; N'}=S^{N'}_{\rm voa}(M^i_{k, j; N'})$
constructed from $M^i_{k, j; N'}$. By Proposition 6.2 in 
\cite{H-aa-int-op}, we see that $\Omega_{N'}^{0}(W^i_{k, j; N'})=M^i_{k, j; N'}$.
Since $V$ has no nonzero elements of negative weights and 
$C_{2}$-cofinite and $M^i_{k, j; N'}$ is finite dimensional,
by Property 2 in Proposition \ref{properties}, 
$W^i_{k, j; N'}$ as a lower-bounded generalized $V$-module generated by 
$M^i_{k, j; N'}$ is quasi-finite dimensional and is in particular grading restricted. 
By Property 5 in Proposition \ref{properties}, 
$W^i_{k, j; N'}$  is of finite length. 

Now we consider the case $N'=N$. 
Given an element of $P_{k, j; N}^{i}$, its action on the 
right $P_{k, j; N}^{i}$-module
$M^i_{k, j; N}$ is in fact an $A^{N}(V)$-module map from 
$M^i_{k, j; N}$ to itself. By the universal property of 
$W^i_{k, j}=S^{N}_{\rm voa}(M^i_{k, j; N})$,
there is a unique $V$-module map from $W^i_{k, j; N}$ to itself such that 
its restriction to $M^i_{k, j; N}$ is the action of the 
element of $P_{k, j; N}^{i}$ on $M^i_{k, j; N}$. Thus we obtain a 
right action of $P_{k, j; N}^{i}$ on $W^i_{k, j; N}$. 
Since the action of $P_{k, j; N}^{i}$ on $W^i_{k, j; N}$ is 
given by $V$-module maps,  the  
homogeneous subspaces $(W^i_{k, j; N})_{[r_{j}+m]}$ 
of $W^i_{k, j; N}$  for $m\in \N$ also
have right $P_{k, j; N}^{i}$-module structures.
We now show that these right $P_{k, j; N}^{i}$-modules are in fact 
projective. 

Since 
$\Omega_{N}^{0}(W^i_{k, j; N})=M^i_{k, j; N}$ is a right 
projective $P_{k, j; N}^{i}$-module,
$$\Omega_{N}^{0}(W^i_{k, j; N})=\coprod_{m=0}^{N}
(W^i_{k, j; N})_{[r_{j}+m]}$$ 
and 
$(W^i_{k, j; N})_{[r_{j}+m]}=(M^i_{k, j; N})_{[m]}$ for $m=0, \dots,
N$ are right $P_{k, j; N}^{i}$-modules, 
$(W^i_{k, j; N})_{[r_{j}+m]}$ for $m=0, \dots,
N$ as direct summands of $M^i_{k, j; N}$ are also projective as 
right $P_{k, j; N}^{i}$-modules. We still need to prove that 
$(W^i_{k, j; N})_{[r_{j}+N']}$ for $N'\in N+\Z_{+}$ are 
projective as right $P_{k, j; N}^{i}$-modules.

Using the isomorphisms  $\mathcal{U}(1): \widetilde{A}^{N'}(V)
\to A^{N'}(V)$ and the $\mathcal{U}_{W}(1): \widetilde{A}^{N'}(W)
\to A^{N'}(W)$, all the elements and structures we obtained 
from $\widetilde{A}^{N'}(V)$ and $\widetilde{A}^{N'}(W)$
are mapped to the corresponding elements and 
structures that one can obtain from $A^{N'}(V)$
and $A^{N'}(W)$ in the same way. Moreover, these 
elements and structures have completely the same properties. 

Let $e^{N'}_{i'}=\mathcal{U}(1)\tilde{e}^{N'}_{i'}$
for $i'=1, \dots, n'$.
Then $\one^{N'}+Q^{\infty}(V)
=e^{N'}_{1}+\cdots +e^{N'}_{n'}$ and $N\le N'$,
we have 
\begin{align*}
\one^{N}+Q^{\infty}(V)&=(\one^{N'}+Q^{\infty}(V))
\diamond (\one^{N}+Q^{\infty}(V))\nn
&=e^{N'}_{1}\diamond (\one^{N}+\widetilde{Q}^{\infty}(V))
+\cdots +e^{N'}_{n'}\diamond 
(\one^{N}+\widetilde{Q}^{\infty}(V)).
\end{align*}
Using the properties of $\one^{N}$ and 
$e_{i'}^{N'}$ for $i=1, \dots, n'$,
we see that $e^{N'}_{i'}\diamond (\one^{N}+Q^{\infty}(V))$
 for $i'=1, \dots, n'$ are all in $A^{N}(V)$.
Let $e^{N}_{1}, \dots, e_{n}^{N}$ be the nonzero elements 
in the set consisting of the elements 
$e^{N'}_{i'}\diamond (\one^{N}+Q^{\infty}(V))$
 for $i'=1, \dots, n'$. Then we have
$$\one^{N}+Q^{\infty}(V)=e^{N}_{1}+ \cdots + e_{n}^{N}.$$
From the properties of $\one^{N}$ and $e_{i'}^{N'}$ for $i'=1, \dots, n'$
 again, we see that $e^{N}_{1}, \dots, e_{n}^{N}$ are 
 orthogonal primitive central idempotents of $A^{N}(V)$.
 Let $A^{N}_{i}=A^{N}(V)\diamond e^{N}_{i}$. Then we have 
 $$A^{N}(V)=A^{N}_{1}\oplus \cdots \oplus A^{N}_{n}.$$
 Also we have $A^{N}(W)
 =(\one^{N}+Q^{\infty}(V))
\diamond  A^{N'}(W)\diamond
(\one^{N}+Q^{\infty}(V))$. Since 
 $A^{N'}(W)=U^{N'}_{1}\oplus \cdots \oplus U^{N'}_{n'}$,
 where $U^{N'}_{i'}=e^{N'}_{i'}\diamond 
A^{N'}(W_{1})\diamond e^{N'}_{i}$ 
 for $i'=1, \dots, n'$, 
 we have 
 \begin{align*}
A^{N}(W)&=
 (\one^{N}+Q^{\infty}(V))\diamond
A^{N'}(W)\diamond (\one^{N}+Q^{\infty}(V))\nn
 &= (\one^{N}+Q^{\infty}(V))\diamond
 U^{N'}_{1}\diamond (\one^{N}+Q^{\infty}(V))
 \nn
 &\quad \oplus  \cdots \oplus 
  (\one^{N}+Q^{\infty}(V))\diamond  U^{N'}_{n'}
  \diamond (\one^{N}+Q^{\infty}(V))\nn
  &=(e^{N'}_{1}\diamond (\one^{N}+Q^{\infty}(V)))
  \diamond A^{N}(W)
\diamond (e^{N'}_{1}\diamond (\one^{N}+Q^{\infty}(V))) \nn
 &\quad
 \oplus  \cdots \oplus (e^{N'}_{n'}\diamond
 (\one^{N}+Q^{\infty}(V)))  \diamond
A^{N}(W)
\diamond (e^{N'}_{n'}\diamond (\one^{N}+Q^{\infty}(V)))\nn
&=U^{N}_{1}\oplus  \cdots \oplus U^{N}_{n},
  \end{align*}
  where $U^{N}_{i'}=e^{N}_{i'}\diamond
  A^{N}(W_{1})\diamond e^{N}_{i'}$
  is an $A^{N}_{i'}$-bimodule
  for $i'=1, \dots, n'$.

By definition, 
$$B^{N}=A^{N}(V)\oplus 
A^{N}(W_{1})\subset 
A^{N'}(V)\oplus 
A^{N'}(W_{1})=B^{N'}$$ 
as associative algebras. Let $B^{N}_{i}=A_{i}^{N}\oplus U_{i}^{N}$
for $i=1, \dots, n$. Then we have 
$$B^{N}=B^{N}_{1}\oplus\cdots \oplus B^{N}_{n}.$$
For $i$ and $i'$ such that $e^{N}_{i}=
e^{N'}_{i'}\diamond (\one^{N}+Q^{\infty}(V))$
is nonzero, we have $\text{Rad}(\phi^{i'}_{k, j; N'})\cap B^{N}
=\text{Rad}(\phi^i_{k, j; N})$. Then the kernel of
the homomorphism of associative algebras
from $B^{N}$ to $B^{N'}/\text{Rad}(\phi^{i'}_{k, j; N'})$
is $\text{Rad}(\phi^i_{k, j; N'})\cap B^{N}=\text{Rad}(\phi^i_{k, j; N})$. 
In particular, we obtain an injective homomorphism of associative algebras
from $B^{N}/\text{Rad}(\phi^i_{k, j; N})$ to 
$B^{N'}/\text{Rad}(\phi^{i'}_{k, j; N'})$. 
This injective homomorphism maps 
$e^{N}_{i}+\text{Rad}(\phi^i_{k, j; N})$
to $e^{N'}_{i'}+\text{Rad}(\phi^{i'}_{k, j; N'})$. 

 Let $e^{N'}_{i', \mu\nu}=\mathcal{U}(1)
  \tilde{e}^{N'}_{i', \mu\nu}$ for $i'=1, \dots, n'$, 
  $\mu=1, \dots, l'_{i'},
\nu=1, \dots, l'_{i'\mu}$. 
For $i$ and $i'$ such that $e^{N}_{i}=
e^{N'}_{i'}\diamond (\one^{N}+Q^{\infty}(V))$
is nonzero, let $e^{N}_{i, \mu\nu}=(\one^{N}+Q^{\infty}(V)) \diamond
e^{N'}_{i', \mu\nu}
\diamond e^{N}_{i}$ for 
$\mu=1, \dots, l_{i}$,
$\nu=1, \dots, l_{i\mu}$.
Then $\{e^{N}_{i, \mu\nu}\}$ 
is a complete set of orthogonal primitive idempodents of 
$A^{N}_{i}$ satisfying
$$1_{A^{N}_{i}}=e^{N}_{i}
=\sum_{\mu=1}^{l_{i}}\sum_{\nu=1}^{l_{i\mu}}
e^{N}_{i, \mu\nu},$$
$e^{N}_{i, \mu\nu_{1}}A^{N}_{i}\simeq e^{N}_{i, \mu\nu_{2}}
A^{N}_{i}$
and 
$e^{N}_{i, \mu_{1}\nu_{1}}A^{N}_{i}\not\simeq
e^{N}_{i, \mu_{2}\nu_{2}}A^{N}_{i}$ 
for $\mu_{1}\ne \mu_{2}$. 
In fact, it is easy to see that $\{e^{N}_{i, \mu\nu}\}$
 is a complete set of orthogonal idempodents of 
$B^{N}_{i}$ satisfying the conditions above. Here we give only a 
proof that these orthogonal idempodents are prime. To see this, 
we note that the left $A^{N}_{i}$-module 
$A^{N}_{i}e^{N}_{i, \mu\nu}$ is in fact also a
left $A^{N}(V)$-module. Since $A^{N}(V)$ is finite dimensional, 
$A^{N}_{i}e^{N}_{i, \mu\nu}$ is also finite dimensional. 
Then it must be a graded $A^{N}(V)$-module. 
Similarly, $A_{i'}^{N'}e_{i', \mu\nu}^{N'}$ is a graded $A^{N'}(V)$-module.
Conclusion 2 in Proposition \ref{more-properties} in fact says that 
the functor $S_{\rm voa}^{N}$  from the category of 
graded $A^{N}(V)$-modules  to the category of lower-bounded 
generalized $V$-modules 
is an equivalence of categories and its inverse is the functor 
$\Omega_{N}^{0}$. The same is true for the functor $S_{\rm voa}^{N'}$
and $\Omega_{N}^{0}$. Since $e_{i', \mu\nu}^{N'}$ is prime, 
$A_{i'}^{N'}e_{i', \mu\nu}^{N'}$ is an indecomposable 
graded $A^{N'}(V)$-module. Thus the lower-bounded 
generalized $V$-module $W=S_{\rm voa}^{N'}(A_{i'}^{N'}e_{i', \mu\nu}^{N'})$
is also indecomposable. Since $N'\ge N$, we know that 
$A^{N}_{i}e^{N}_{i, \mu\nu}$
is a graded subspace of $A_{i'}^{N'}e_{i', \mu\nu}^{N'}$.
By Conclusion 3 in Proposition \ref{more-properties}, we see that 
$A^{N}_{i}e^{N}_{i, \mu\nu}=\Omega_{N}^{0}(W)$. 
If $A^{N}_{i}e^{N}_{i, \mu\nu}$ can be decomposed into the direct sum 
of two proper graded $A^{N}(V)$-modules, then 
$S_{\rm voa}^{N}(A_{i}^{N}e_{i, \mu\nu}^{N'})$ can  be 
decomposed into the direct sum 
of two proper lower-bounded generalized $V$-modules. 
But by Conclusion 2 in Proposition \ref{more-properties} again,
$$S_{\rm voa}^{N}(A_{i}^{N}e_{i, \mu\nu}^{N'})=
S_{\rm voa}^{N}(\Omega_{N}^{0}(W))\simeq W.$$
Thus we conclude that $W$  can  be 
decomposed into the direct sum 
of two proper lower-bounded generalized $V$-modules. Contradiction. 
So $A^{N}_{i}e^{N}_{i, \mu\nu}$ cannot be decomposed into the direct sum 
of two proper graded $A^{N}(V)$-modules. This is equivalent to the fact 
that $e^{N}_{i, \mu\nu}$ is prime. 

Since  $B_{i}^{N}=A_{i}^{N}+U_{i}^{N}$, 
we see that $\{e^{N}_{i, \mu\nu}\}$ 
is also a complete set of orthogonal primitive idempodents of 
$B^{N}_{i}$ satisfying
$$1_{B^{N}_{i}}=e^{N}_{i}
=\sum_{\mu=1}^{l_{i}}\sum_{\nu=1}^{l_{i\mu}}
e^{N}_{i, \mu\nu},$$
$e^{N}_{i, \mu\nu_{1}}B^{N}_{i}\simeq e^{N}_{i, \mu\nu_{2}}
B^{N}_{i}$
and 
$e^{N}_{i, \mu_{1}\nu_{1}}B^{N}_{i}\not\simeq
e^{N}_{i, \mu_{2}\nu_{2}}B^{N}_{i}$ 
for $\mu_{1}\ne \mu_{2}$. 
Let 
$\epsilon^{N}_{i, \mu\nu}=e^{N}_{i, \mu\nu}
+\text{Rad}(\phi^i_{k, j; N})$ for 
$\mu=1, \dots, l'_{i'},
\nu=1, \dots, l'_{i'\mu}$. Then $\{\epsilon^{N}_{i, \mu\nu}\}$ 
is a set of orthogonal primitive idempodents of 
$B^{N}_{i}/\text{Rad}(\phi^i_{k, j; N})
=B^{N}/\text{Rad}(\phi^i_{k, j; N})$
and  
$$1_{B^{N}/\text{Rad}(\phi^i_{k, j; N})}
=e^{N}_{i}+\text{Rad}(\phi^i_{k, j; N})
=\sum_{\mu=1}^{l_{i}}\sum_{\nu=1}^{l_{i\mu}}
\epsilon^{N}_{i, \mu\nu}.$$
Moreover, we also have 
$$\epsilon^{N}_{i, \mu\nu_{1}}(B^{N}
/\text{Rad}(\phi^i_{k, j; N}))\simeq \epsilon^{N}_{i, \mu\nu_{2}}
(B^{N}/\text{Rad}(\phi^i_{k, j; N}))$$
and 
$$\epsilon^{N}_{i, \mu_{1}\nu_{1}}(B^{N}
/\text{Rad}(\phi^i_{k, j; N}))\not\simeq \epsilon^{N}_{i, \mu_{2}\nu_{2}}
(B^{N}/\text{Rad}(\phi^i_{k, j; N})$$ 
for $\mu_{1}\ne \mu_{2}$.

Let 
$$\epsilon_{i}^{N}=\sum_{\mu=1}^{l_{i}}
\epsilon^{N}_{i, \mu 1}.$$
Then we have 
$$P_{k, j; N}^{i}=\epsilon_{i}^{N}\diamond (B^{N}/
\text{Rad}(\phi^i_{k, j; N}))\diamond \epsilon_{i}^{N}$$
and 
$$M^i_{k, j; N}=(B^{N}/
\text{Rad}(\phi^i_{k, j; N}))\diamond \epsilon_{i}^{N}.$$
Thus we obtain an injective homomorphism of associative algebras
from $P_{k, j; N}^{i}=\epsilon_{i}^{N}\diamond (B^{N}/
\text{Rad}(\phi^i_{k, j; N}))\diamond \epsilon_{i}^{N}$
to $P_{k, j; N'}^{i'}=\epsilon_{i'}^{N'}\diamond (B^{N'}/
\text{Rad}(\phi^{i'}_{k, j; N}))\diamond \epsilon_{i'}^{N}$.
We shall view $P_{k, j; N}^{i}$ as a subalgebra
of $P_{k, j; N'}^{i'}$ from now on. In particular, 
$M^{i'}_{k, j; N'}=(B^{N'}/
\text{Rad}(\phi^{i'}_{k, j; N'}))\diamond \epsilon_{i'}^{N'}$
is also a right $P_{k, j; N}^{i}$-module.
Then we also obtain an injective $P_{k, j; N}^{i}$-module map 
from $M^i_{k, j; N}=(B^{N}/
\text{Rad}(\phi^i_{k, j; N}))\diamond \epsilon_{i}^{N}$
to $M^{i'}_{k, j; N'}=(B^{N'}/
\text{Rad}(\phi^{i'}_{k, j; N'}))\diamond \epsilon_{i'}^{N'}$. 
Moreover, $M^i_{k, j; N}$ and $M^{i'}_{k, j; N'}$
are left $B^{N}$-module and $B^{N'}$-module, respectively, 
and this injective 
$P_{k, j; N}^{i}$-module map induces the left $B^{N}$-module
structure on $M^i_{k, j; N}$ from the left $B^{N'}$-module
structure on $M^{i'}_{k, j; N'}$.  

We know that $M^{i'}_{k, j; N'}$ as a right $P_{k, j; N'}^{i'}$-module
is projective. We now show that 
$M^{i'}_{k, j; N'}$ as a right $P_{k, j; N}^{i}$-module is also projective.
It is enough to show that the right action of $P_{k, j; N'}^{i'}$ 
on $M^{i'}_{k, j; N'}$
is in fact determined by the right action of $P_{k, j; N}^{i}$. 
In fact, we have proved that  $M^i_{k, j; N}$ and $M^{i'}_{k, j; N'}$
are $A^{N}(V)$-$P^{i}_{k, j; N}$-bimodule and 
$A^{N'}(V)$-$P^{i'}_{k, j; N'}$-bimodule, respectively and 
the left $A^{N}(V)$-module structure on $M^i_{k, j; N}$
is induced from the left $A^{N'}(V)$-module structure on 
$M^{i'}_{k, j; N'}$ when we view $M^i_{k, j; N}$ as
a right $P^{i}_{k, j; N}$-submodule of $M^{i'}_{k, j; N'}$.
Then we have
$$M^{i'}_{k, j; N'}=\coprod_{m=0}^{N'}(M^{i'}_{k, j; N'})_{[m]}$$ 
and 
$$M^i_{k, j; N}=
(M^{i'}_{k, j; N'})^{N}=\coprod_{m=0}^{N}(M^{i'}_{k, j; N'})_{[m]}.$$ 
By Proposition \ref{more-properties}, $M^{i'}_{k, j; N'}$ is equivalent 
to the graded $A^{N'}(V)$-module
$\Omega_{N'}^{0}(S_{\rm voa}^{N}(M^i_{k, j; N}))
=\Omega_{N'}^{0}(W^i_{k, j})$. 
We know that for $m=1, \dots, N'$, $(M^{i'}_{k, j; N'})_{[m]}$
are right $P^{i'}_{k, j; N'}$-submodules of $M^{i'}_{k, j; N'}$.
In particular, $M^{i}_{k, j; N}$ is a right 
$P^{i'}_{k, j; N'}$-submodule of $M^{i'}_{k, j; N'}$.
Note that the action of  every element of $P^{i'}_{k, j; N'}$
on $M^{i}_{k, j; N}$ is an $A^{N}(V)$-module map. 
By the universal property of $S_{\rm voa}^{N}(M^i_{k, j; N})$,
such an $A^{N}(V)$-module map gives a unique 
$V$-module map from $S_{\rm voa}^{N}(M^i_{k, j; N})$ to
itself. In particular, such an $A^{N}(V)$-module map gives a unique 
$A^{N'}(V)$-module map from $M^{i'}_{k, j; N'}$. Thus 
we see that the action of the element of $P^{i'}_{k, j; N'}$
on $M^{i'}_{k, j; N'}$ must be the one obtained from its restriction to
$M^{i}_{k, j; N}$. But the restriction to $M^{i}_{k, j; N}$ 
 of the action of $P^{i'}_{k, j; N'}$
on $M^{i'}_{k, j; N'}$  is exactly the action of $P^{i}_{k, j; N}$. 
So we see that the action of $P^{i'}_{k, j; N'}$ 
on $M^{i'}_{k, j; N'}$ is determined by the action 
of $P^{i}_{k, j; N}$
on $M^{i}_{k, j; N}$. Thus $M^{i'}_{k, j; N'}$ as a right 
$P^{i}_{k, j; N}$-module is also projective. In particular,
$(W^i_{k, j})_{[r_{j}+N']}=(M^{i'}_{k, j; N'})_{[N']}$
is projective as a right 
$P^{i}_{k, j; N}$-module. Since $N'$ is arbitrary, we see that 
$W^i_{k, j}$ and its homogeneous subspaces are all 
 projective as a right $P^{i}_{k, j; N}$-modules.

We know that 
$\Omega_{N}^{0}(W^i_{k, j})=
\Omega_{N}^{0}(S_{\rm voa}^{N}(M^i_{k, j; N}))
=M^i_{k, j; N}$. Then by Proposition \ref{more-properties}, 
$(W^i_{k, j})'$ is equivalent to 
$S_{\rm voa}^{N}(\Omega_{N}^{0}((W^i_{k, j})'))$.
In particular, Theorem \ref{main-N} can be applied to the case
that $W_{1}=W$ and 
$W_{2}=W_{3}=W^i_{k, j}=S_{\rm voa}^{N}(M^i_{k, j; N})$. 
By this theorem, we obtain a unique $P^{i}_{k, j; N}$-compatible 
intertwining operator 
$\Y^{i}_{k, j}$ of type $\binom{W^i_{k, j}}
{WW^i_{k, j}}$ such that 
$$\rho^{N}(\Y^{i}_{k, j})=
f^i_{k, j; N}\in \hom_{\widetilde{A}^{N'}(V), P^{i}_{k, j; N}}
(\widetilde{A}^{N}(W)\otimes_{\widetilde{A}^{N}(V)} 
M^i_{k, j; N}, M^i_{k, j; N}).$$
By Proposition \ref{slf-int-op}, we have a symmetric linear function
$\psi_{\Y^{i}_{k, j}, \phi^i_{k, j; N}}$
on $\widetilde{A}^{N}(W)$. 

We actually need only the intertwining operators $\Y^{i}_{0, j}$ 
(that is, the case $k=0$). It is clear that the linear map given by 
\begin{equation}\label{difference-slf}
w\mapsto F(w; \tau)-\sum_{j=1}^{J}\sum_{i=1}^{n}
\overline{F}_{\Y^{i}_{0, j}}^{\phi^i_{0, j; N}}(w; \tau)
\end{equation}
for $w\in W$ is also a genus-one $1$-point conformal block labeled by $W$. Then by Proposition \ref{diff-conf-bk},
$$F(w; \tau)-\sum_{j=1}^{J}\sum_{i=1}^{n}
\overline{F}_{\Y^{i}_{0, j}}^{\phi^i_{0, j; N}}(w; \tau)$$
can be expanded as  
$$\sum_{j=1}^{J^{(1)}}\sum_{m\in \mathbb{N}}
C^{(1)}_{0, j, m}(w)q_{\tau}^{r^{(1)}_{j}+m}
+\sum_{k=1}^{K}\tau^{k} G^{(1)}_{k}(w;  q_{\tau}),$$
where for $j=1, \dots, J^{(1)}$, there exists $j'$ satisfying $1\le j'\le J$
such that $r^{(1)}_{j}-r_{j'}\in \Z_{+}$ and where
$G^{(1)}_{k}(w; q_{\tau})$ for $k=1, \dots, K$
are in $\coprod_{r\in \C}q_{\tau}^{r}\C[[q_{\tau}]]$.
Let $s\in \Z_{+}$ be larger than the maximum of the real parts
of the differences of the lowest weights of 
the (finitely many) irreducible $V$-modules.
Then we repeat the argument above $s$ times to
 find finite-dimensional associative algebras $P^{(s)}_{l}$
with symmetric linear functions $\phi^{(s)}_{l}$, 
grading-restricted generalized $V$-$P^{(s)}$-modules 
$W^{(s)}_{l}$ and $P^{(s)}_{l}$-compatible intertwining operators
$\Y^{(s)}_{l}$ of types $\binom{W^{(s)}_{l}}
{WW^{(s)}_{l}}$  for $l=1, \dots, p^{(s)}$ such that 
the linear map given by
$$w\mapsto F(w; \tau)-\sum_{l=1}^{p^{(s)}}
\overline{F}_{\Y^{(s)}_{l}}^{\phi^{(s)}_{l}}(w; \tau)$$
is a genus-one $1$-point conformal block 
and for $w\in W$, 
\begin{equation}\label{s-difference-slf}
F(w; \tau)-\sum_{l=1}^{p^{(s)}}
\overline{F}_{\Y^{(s)}_{l}}^{\phi^{(s)}_{l}}(w; \tau)
\end{equation}
can be expanded as
\begin{equation}\label{expansion-form}
\sum_{j=1}^{J^{(s)}}\sum_{m\in \mathbb{N}}
C^{(s)}_{0, i, m}(w)q_{\tau}^{r^{(s)}_{j}+m}
+\sum_{k=1}^{K}\tau^{k}G^{(s)}_{k}(w;  q_{\tau}),
\end{equation}
where for $j=1, \dots, J^{(s)}$, there exists $j'$ satisfying $1\le j'\le J$
such that $r^{(s)}_{j}-r_{j'}\in s+\N$ and where
$G^{(s)}_{k}(w;  q_{\tau})$ for $k=1, \dots, K$
is in $\coprod_{r\in \C}q_{\tau}^{r}\C[[q_{\tau}]]$. 
We now show that (\ref{s-difference-slf})
is equal to $0$.

We first prove that the first term in (\ref{expansion-form}) is $0$ 
for all $w\in W$.
Assume that it is not $0$ for some $w\in W$.
Note that we have shown that 
$r_{j}+\frac{c}{24}$ for $j=1, \dots, J$
are the lowest weights of $M^i_{k, j; N}$ and thus
are also the lowest weights of $W^i_{k, j}$. But 
the lowest weight of a lower-bounded 
generalized $V$-module of finite length must also be the 
lowest weight of an irreducible $V$-module. In particular, 
$r_{j}+\frac{c}{24}$ for $j=1, \dots, J$ are 
lowest weights of irreducible $V$-modules. Since 
the first term in (\ref{expansion-form}) is not $0$ for some 
$w\in W$,
the same proof shows that $r_{j}^{(s)}+\frac{c}{24}$ 
for $j=1, \dots, J^{(s)}$ 
are also lowest weights of irreducible $V$-modules. 
But $r^{(s)}_{j}-r_{j'}\in s+\N$, we have $\Re(r^{(s)}_{j}-r_{j'})\ge s$
which is larger than the maximum of the real parts of
the differences of  the lowest weights of 
the (finitely many) irreducible $V$-modules. This is impossible since 
as we have shown above, $r^{(s)}_{j}-r_{j'}$ is the difference 
between the lowest weights $r^{(s)}_{j}$ and $r_{j'}$ of some irreducible $V$-modules. 
Thus the first term in (\ref{expansion-form})  must be $0$. 

In fact, if we write the first term in (\ref{expansion-form}) as 
$G_{0}(w; q_{\tau})$, then (\ref{expansion-form}) 
can be written as 
$$\sum_{k=0}^{K}\tau^{k}G^{(s)}_{k}(w; q_{\tau}).$$
We have proved that $G_{0}(w; q_{\tau})=0$. 
We now use induction to prove that $G_{k}(w; q_{\tau})=0$ for 
$k=0, \dots, K$. Assume that $G_{k}(w; q_{\tau})=0$
for $k=0, \dots, k_{0}$. Then  (\ref{expansion-form}) 
becomes 
$$\sum_{k=k_{0}+1}^{K}\tau^{k}G^{(s)}_{k}(w; q_{\tau}).$$
Note that  (\ref{s-difference-slf}) gives a 
genus-one $1$-point conformal block. In particular, 
the property (\ref{1-pt-cb-2}) for this 
genus-one $1$-point conformal block gives
\begin{equation}\label{mod-inv-8-1-G-k}
2\pi i\frac{\partial}{\partial \tau}
\sum_{k=k_{0}+1}^{K}\tau^{k}G^{(s)}_{k}(w; q_{\tau})
=\sum_{k=k_{0}+1}^{K}\tau^{k}G^{(s)}_{k}
(\res_{x}(\tilde{\wp}_{1}(x; q_{\tau})-\tilde{G}_{2}(q_{\tau})x)
Y_{W}(\omega, x)w; q_{\tau})
\end{equation}
(see (\ref{tilde-wp-g-2})  for the definition of 
$\tilde{\wp}_{1}(x; q_{\tau})-\tilde{G}_{2}(q_{\tau})x$). 

The left-hand side of (\ref{mod-inv-8-1-G-k}) is equal to 
\begin{equation}\label{dev-G-k}
2\pi i
\sum_{k=k_{0}+1}^{K}k\tau^{k-1}G^{(s)}_{k}(w; q_{\tau})
+2\pi i
\sum_{k=k_{0}+1}^{K}\tau^{k}\frac{\partial}{\partial \tau}
G^{(s)}_{k}(w; q_{\tau}).
\end{equation}
Since the coefficients of $\tau^{k_{0}}$ 
in the right-hand side of (\ref{mod-inv-8-1-G-k}) and in  (\ref{dev-G-k}) 
 are $0$ and
$2\pi i(k_{0}+1)
G^{(s)}_{k_{0}+1}(w; q_{\tau})$, respectively, 
we obtain $G^{(s)}_{k_{0}+1}(w; q_{\tau})=0$. 
By induction principle, we obtain $G_{k}(w; q_{\tau})=0$ for 
$k=0, \dots, K$. Thus  (\ref{expansion-form})  is $0$.

This proves (\ref{s-difference-slf}) is $0$, that is, 
\begin{equation}\label{mod-trans=trace}
F(w; \tau)=\sum_{j=l}^{p^{(s)}}
\overline{F}_{\Y^{(s)}_{l}}^{\phi^{(s)}_{l}}(w; q_{\tau})
\end{equation}
for $w\in W$.
Since the right-hand side of (\ref{mod-trans=trace}) is 
in  $\mathcal{F}_{w}$, the left-hand side of 
(\ref{mod-trans=trace}) is also in $\mathcal{F}_{w}$, 
proving our theorem.
\epfv

Let $W_{1}, \dots, W_{n}$ be grading-restricted generalized $V$-modules
and $w_{1}\in W_{1}, \dots, w_{n}\in W_{n}$. Given a finite-dimensional 
associative algebra $P$ with a symmetric linear function $\phi$, 
grading-restricted generalized $V$-modules 
$\widetilde{W}_{1}, \dots, \widetilde{W}_{n-1}$,
a grading-restricted generalized $V$-$P$-bimodule 
$\widetilde{W}_{0}=\widetilde{W}_{n}$ projective as a right $P$-module and
intertwining operators $\Y_{1}, \dots, \Y_{n}$ of types 
$\binom{\widetilde{W}_{0}}{W_{1}\widetilde{W}_{1}}, \dots, 
\binom{\widetilde{W}_{n-1}}{W_{n}\widetilde{W}_{n}}$, respectively, such that the product of 
$\Y_{1}, \dots, \Y_{n}$ is compatible with $P$, we have 
a genus-one correlation $n$-point function
$$\overline{F}^{\phi}_{\mathcal{Y}_{1}, \dots, \mathcal{Y}_{n}}(w_{1}, \dots, w_{n};
z_{1}, \dots, z_{n}; \tau).$$
Let $\mathcal{F}_{w_{1}, \dots, w_{n}}$ be the vector space of such 
genus-one correlation $n$-point functions for all $P$, $\phi$, 
$\widetilde{W}_{1}, \dots, \widetilde{W}_{n-1}$, $\widetilde{W}_{0}=\widetilde{W}_{n}$
and $\Y_{1}, \dots, \Y_{n}$.

We are ready to prove the main theorem of this paper.

\begin{thm}\label{mod-inv}
Let $V$ be a $C_{2}$-cofinite vertex operator algebra 
without nonzero negative weight elements.
Then for a finite-dimensional 
associative algebra $P$ with a symmetric linear function $\phi$, 
grading-restricted generalized $V$-modules 
$\widetilde{W}_{1}, \dots, \widetilde{W}_{n-1}$,
a grading-restricted generalized $V$-$P$-bimodule 
$\widetilde{W}_{0}=\widetilde{W}_{n}$ projective as a right $P$-module,
intertwining operators $\Y_{1}, \dots, \Y_{n}$ of types 
$\binom{\widetilde{W}_{0}}{W_{1}\widetilde{W}_{1}}$, $\dots$, 
$\binom{\widetilde{W}_{n-1}}{W_{n}\widetilde{W}_{n}}$, respectively, such that 
product of 
$\Y_{1}, \dots, \Y_{n}$ is compatible with $P$,
and 
$$\left(\begin{array}{cc}
\alpha&\beta\\
\gamma&\delta
\end{array}\right)\in SL(2, \mathbb{Z}),$$
\begin{align*}
&\overline{F}^{\phi}_{\mathcal{Y}_{1}, \dots, \mathcal{Y}_{n}}
\Biggl(\left(\frac{1}{\gamma\tau+\delta}\right)^{L_{W_{1}}(0)}w_{1}, \dots,
\left(\frac{1}{\gamma\tau+\delta}\right)^{L_{W_{n}}(0)}w_{n};
\frac{z_{1}}{\gamma\tau+\delta}, \dots, \frac{z_{n}}{\gamma\tau+\delta}; 
\frac{\alpha\tau+\beta}{\gamma\tau+\delta}\Biggr)
\end{align*}
is in $\mathcal{F}_{w_{1}, \dots, w_{n}}$. 
\end{thm}
\pf
As in the proof of the modular invariance (Theorem 7.3 in 
\cite{H-mod-inv-int}) in the rational case,
using the genus-one associativity proved 
in \cite{F1} and \cite{F} (see Theorem \ref{F-genus1-assoc}),
the proof in the general $n$ case is reduced to the $n=1$ case. So
we need only prove the $n=1$ case. 

In the $n=1$ case, by Proposition \ref{mod-transf-conf-bk},
the linear map given by 
\begin{align*}
w_{1}\mapsto \overline{F}^{\phi}_{\mathcal{Y}_{1}}
\Biggl(\left(\frac{1}{\gamma\tau+\delta}\right)
^{L_{W_{1}}(0)}w_{1};
\frac{z_{1}}{\gamma\tau+\delta}; 
\frac{\alpha\tau+\beta}{\gamma\tau+\delta}\Biggr)
\end{align*}
is a genus-one $1$-point conformal block. Then by 
Theorem \ref{conf-bk=>pseudo-q-tr}, for $w_{1}\in W_{1}$,
$$\overline{F}^{\phi}_{\mathcal{Y}_{1}}
\Biggl(\left(\frac{1}{\gamma\tau+\delta}\right)
^{L_{W_{1}}(0)}w_{1};
\frac{z_{1}}{\gamma\tau+\delta}; 
\frac{\alpha\tau+\beta}{\gamma\tau+\delta}\Biggr)$$
is in $\mathcal{F}_{w_{1}}$, proving the theorem in this
case.  
\epfv

\appendix

\renewcommand{\theequation}{\thesection.\arabic{equation}}
\renewcommand{\thethm}{\thesection.\arabic{thm}}
\setcounter{equation}{0} \setcounter{thm}{0} 

\section{The Weierstrass $\wp$-function $\wp_{2}(z; \tau)$,
the Weierstrass $\zeta$-function $\wp_{1}(z; \tau)$
and the Eisenstein series $G_{2}(\tau)$}

We recall in this appendix 
some basic facts on the Weierstrass $\wp$-function, denoted by 
$\wp_{2}(z; \tau)$ in this paper,
the Weierstrass $\zeta$-function, denoted by $\wp_{1}(z; \tau)$ in 
this paper,
and the Eisenstein series $G_{2}(\tau)$. For details,
see \cite{L} and \cite{K}. 
The Weierstrass  $\wp$-function $\wp_{2}(z; \tau)$ has the $q$-expansion
\begin{align}\label{mod-inv-9.01}
\wp_{2}(z; \tau)=(2\pi i)^{2}q_{z}(q_{z}-1)^{-2}+(2\pi i)^{2}
\sum_{s\in \Z_{+}}\sum_{l|s}l(q_{z}^{l}+q_{z}^{-l})q_{\tau}^{s}
-\frac{\pi^{2}}{3}-2(2\pi i)^{2}\sum_{l\in \Z_{+}}\sigma(l)q_{\tau}^{l}
\end{align}
in the region given by 
$|q_{\tau}|<|q_{z}|<|q_{\tau}|^{-1}$ and $z\ne 0$,
where $\sigma(l)=\sum_{n|l}n$ for $l\in \Z_{+}$. 
Note that the coefficients of the power series (\ref{mod-inv-9.01})
in $q_{\tau}$ are holomorphic functions of $z$ on the whole complex plane 
except for the the coefficient $(2\pi i)^{2}q_{z}(q_{z}-1)^{-2}$ of $q_{\tau}^{0}$,
which has a pole of order $2$ at $z=0$. 
Let 
\begin{align}\label{tilde-wp-2}
\tilde{\wp}_{2}(x; q)
&=(2\pi i)^{2}e^{2\pi ix}(e^{2\pi ix}-1)^{-2}+(2\pi i)^{2}
\sum_{s\in \Z_{+}}\sum_{l|s}l(e^{2l\pi ix}+e^{-2l\pi ix})q^{s}\nn
&\quad -\frac{\pi^{2}}{3}-2(2\pi i)^{2}\sum_{l\in \Z_{+}}\sigma(l)q^{l},
\end{align}
where $e^{2\pi ix}(e^{2\pi ix}-1)^{-2}$ is understood as the 
formal Laurent series obtained by expanding $e^{2\pi iz}(e^{2\pi iz}-1)^{-2}$ 
as a Laurent series near $z=0$ and then replacing $z$ by $x$,
$e^{2l\pi ix}=\sum_{k\in \N}\frac{(2l\pi ix)^{k}}{k!}$
and $e^{-2l\pi ix}=\sum_{k\in \N}\frac{(-2l\pi ix)^{k}}{k!}$. 
In terms of only formal variables, the formal Laurent series 
$e^{2\pi ix}(e^{2\pi ix}-1)^{-2}$ can also be 
obtained as follows: 
Write $(e^{2\pi ix}-1)^{-2}$ as $(2\pi ix)^{-2}(1+\sum_{k\in \Z_{+}+1}
\frac{(2\pi ix)^{k-1}}{k!})^{-2}$ and then expand $(1+\sum_{k\in \Z_{+}+1}
\frac{(2\pi ix)^{k-1}}{k!})^{-2}$ using binomial expansion 
as a power series in $\sum_{k\in \Z_{+}+1}
\frac{(2\pi ix)^{k-1}}{k!}$. Since the $l$-th power of 
$\sum_{k\in \Z_{+}+1}
\frac{(2\pi ix)^{k-1}}{k!}$ is a power series in $x$ such that the 
coefficients of $x^{j}$ for $j=0, \dots, l-1$ are $0$, this
expansion of $(1+\sum_{k\in \Z_{+}+1}
\frac{(2\pi ix)^{k-1}}{k!})^{-2}$ gives a well-defined formal power series in $x$. 
Multiplying $(2\pi ix)^{-2}$ with this formal power series in $x$,
we obtain a formal Laurent series expansion of $(e^{2\pi ix}-1)^{-2}$.
Multiplying this formal Laurent series expansion of $(e^{2\pi ix}-1)^{-2}$
with  the formal power series 
$e^{2l\pi ix}=\sum_{k\in \N}\frac{(2l\pi ix)^{k}}{k!}$, we obtain 
a formal Laurent series expansion of $e^{2\pi ix}(e^{2\pi ix}-1)^{-2}$
which is the same as the formal Laurent series expansion obtained using 
complex analysis.

Similarly, the Weierstrass $\zeta$-function $\wp_{1}(z; \tau)$ minus 
$G_{2}(\tau)z$ 
has the  $q$-expansion
\begin{align}\label{wp-1-g-2}
\wp_{1}(z; \tau)-G_{2}(\tau)z
&=2\pi iq_{z}(q_{z}-1)^{-1}
-2\pi i\sum_{s\in \Z_{+}}\sum_{l|s}(q_{z}^{l}-q_{z}^{-l})
q_{\tau}^{s}-\pi i
\end{align}
in the region given by 
$|q_{\tau}|<|q_{z}|<|q_{\tau}|^{-1}$ and $z\ne 0$.
Let 
\begin{align}\label{tilde-wp-g-2}
\tilde{\wp}_{1}(x; q)-\widetilde{G}_{2}(q)x
&=2\pi ie^{2\pi ix}(e^{2\pi ix}-1)^{-1}-2\pi i
\sum_{s\in \Z_{+}}
\sum_{l|s}
(e^{2l\pi ix}-e^{-2l\pi ix})q^{s}
-\pi i,
\end{align}
where $e^{2\pi ix}(e^{2\pi ix}-1)^{-1}$ is understood as the 
formal Laurent series obtained by expanding $e^{2\pi iz}(e^{2\pi iz}-1)^{-1}$ 
as a Laurent series near $z=0$ and then replacing $z$ by $x$. 
The formal Laurent series $e^{2\pi ix}(e^{2\pi ix}-1)^{-1}$ can also 
be obtained in terms of only formal variables in a way completely similar to
that for $e^{2\pi ix}(e^{2\pi ix}-1)^{-2}$ above.

We also have the Laurent series expansion in $z$
\begin{align*}
\wp_{2}(z; \tau)&=\frac{1}{z^{2}}+\sum_{k\in \Z_{+}}
(2k+1)G_{2k+2}(\tau)z^{2k},\\ 
\wp_{1}(z; \tau)-G_{2}(\tau)z&=\frac{1}{z}-\sum_{k\in \N}
G_{2k+2}(\tau)z^{2k+1}
\end{align*}
in the region $0<|z|<\min(1, |\tau|)$, where $G_{2k+2}(\tau)$ are the 
Eisenstein series. Let 
\begin{align}
\wp_{2}(x; \tau)&=\frac{1}{x^{2}}+\sum_{k\in \Z_{+}}
(2k+1)G_{2k+2}(\tau)x^{2k},\label{wp-2-x}\\ 
\wp_{1}(x; \tau)-G_{2}(\tau)x&=\frac{1}{x}-\sum_{k\in \N}
G_{2k+2}(\tau)x^{2k+1}\label{wp-1-g-2-x}
\end{align}

\begin{lemma}\label{tilde-wp-2-and-1}
The formal Laurent series 
 $\tilde{\wp}_{2}(x; q)$ and $\tilde{\wp}_{1}(x; q)-\widetilde{G}_{2}(q)x$
in $x$ and $q$ can be obtained by expanding 
the coefficients of $\wp_{2}(x; \tau)$ 
and $\wp_{1}(x; \tau)-G_{2}(\tau)x$, respectively, 
in powers of $x$ as 
power series in $q_{\tau}$ and then replacing $q_{\tau}$ by $q$. 
\end{lemma}
\pf
Since (\ref{mod-inv-9.01}) is absolutely convergent in 
the region given by $|q_{\tau}|<|q_{z}|<|q_{\tau}|^{-1}$ and $z\ne 0$
and it has a pole of order $2$ at $z=0$, the 
double series (\ref{tilde-wp-2})
with $x$ and $q$ replaced by $z$ and $q_{\tau}$  is 
absolutely convergent 
in the region given by $|q_{\tau}|<\epsilon$ and $0<|z|<\frac{-\log \epsilon}{2\pi}$
for any $\epsilon \in (0, 1)$. 
In particular, after we substitute $z$ and $q_{\tau}$ for $x$ and $q$,
we can first sum over the powers of $q_{\tau}$ in the region 
$|q_{\tau}|<\epsilon$ to obtain a Laurent series in $z$ which is absolutely 
convergent in the region $0<|z|<\frac{-\log \epsilon}{2\pi}$ to 
$\wp_{2}(z; \tau)$.
Substituting $x$ for $z$ in this Laurent series, we obtain a 
formal Laurent series in $x$ which by definition is 
equal to $\wp_{2}(x; \tau)$. This is equivalent to the statement for 
$\tilde{\wp}_{2}(x; q)$. 

The proof for the statement for 
$\tilde{\wp}_{1}(x; q)-\widetilde{G}_{2}(q)x$ is the same. 
\epfv

\renewcommand{\theequation}{\thesection.\arabic{equation}}
\renewcommand{\thethm}{\thesection.\arabic{thm}}
\setcounter{equation}{0} \setcounter{thm}{0} 

\section{Identities involving binomial coefficients}

We recall and prove some identities 
involving binomial coefficients in this appendix.

\begin{lemma}
\begin{equation}\label{main-lemma-21}
\sum_{k=0}^{n}\binom{-n-1}{k}
(e^{2\pi i(n+1)x}+(-1)^{-n-k-1}e^{2\pi ikx})
(e^{2\pi ix}-1)^{-n-k-1}=1.
\end{equation}
\end{lemma}
\pf
Using the identity 
$$\sum_{k=0}^{n}\binom{k+n}{n}
\frac{(-1)^{k}(1+x)^{n+1}-(-1)^{n}(1+x)^{k}}{x^{n+k+1}}=1$$
given by Proposition 5.2 in \cite{DLM}, we have
\begin{align*}
&\sum_{k=0}^{n}\binom{-n-1}{k}
(e^{2\pi i(n+1)x}+(-1)^{-n-k-1}e^{2\pi ikx})
(e^{2\pi ix}-1)^{-n-k-1}\nn
&\quad =\sum_{k=0}^{n}\binom{k+n}{n}
\frac{(-1)^{k}e^{2\pi i(n+1)x}-(-1)^{n}e^{2\pi ikx}}
{(e^{2\pi ix}-1)^{n+k+1}}\nn
&\quad =\sum_{k=0}^{n}\binom{k+n}{n}
\frac{(-1)^{k}(1+(e^{2\pi ix}-1))^{n+1}-(-1)^{n}(1+(e^{2\pi ix}-1))^{k}}
{(e^{2\pi ix}-1)^{n+k+1}}\nn
&\quad =1.
\end{align*}
\epfv

\begin{lemma}
For $m, n, k, l\in \Z_{+}$ satisfying $n\ge m\ge l$, we have
\begin{equation}\label{lemma-3-1}
\sum_{j=0}^{m}
\frac{l}{(n-j+k)}\binom{m}{j}\binom{-l-1}{n-j+k-1}
=0.
\end{equation}
\end{lemma}
\pf
For $j=0, \dots, m$, by the definition of
the binomial coefficients, we have
\begin{align*}
\frac{l}{(n-j+k)}\binom{-l-1}{n-j+k-1}
&=\frac{l}{(n-j+k)}\frac{(-l-1)\cdots (-l-1-n+j-k+1+1)}{(n-j+k-1)!}\nn
&=(-1)^{n-j+k-1}\frac{(n-j+l+k-1)\cdots (l+1)l(l-1)!}{(l-1)!(n-j+k)!}\nn
&=(-1)^{n+k-1}(-1)^{j}\frac{(n-j+l+k-1)\cdots (n-j+k+1)}{(l-1)!}.
\end{align*}
Then 
\begin{align}\label{main-lemma-8}
&\sum_{j=0}^{m}
\frac{l}{(n-j+k)}\binom{m}{j}\binom{-l-1}{n-j+k-1}\nn
&\quad=\sum_{j=0}^{m}(-1)^{n+k-1}(-1)^{j}\binom{m}{j}
\frac{(n-j+l+k-1)\cdots (n-j+k+1)}{(l-1)!}\nn
&\quad =\sum_{j=0}^{m}(-1)^{n+k-1}(-1)^{j}\binom{m}{j}
\frac{(n-j+l+k-1)\cdots (n-j+k+1)}{(l-1)!}x^{n-j+k}\lbar_{x=1}\nn
&\quad =\frac{(-1)^{n+k-1}}{(l-1)!}\frac{d^{l-1}}{dx^{l-1}}
\sum_{j=0}^{m}(-1)^{j}\binom{m}{j}
x^{n-j+l+k-1}\lbar_{x=1}\nn
&\quad =\frac{(-1)^{n+k-1}}{(l-1)!}\frac{d^{l-1}}{dx^{l-1}}
x^{n+l+k-1}(1-x^{-1})^{m}
\lbar_{x=1}\nn
&\quad =\frac{(-1)^{n+k-1}}{(l-1)!}
\frac{d^{l-1}}{dx^{l-1}}x^{n-m+l+k-1}(x-1)^{m}
\lbar_{x=1}\nn
&\quad =0,
\end{align}
where in the last step, we have used $l-1< m$. 
\epfv

\begin{lemma}
For $\alpha, \beta\in \C$ and $m, n\in \N$, 
\begin{align}
\sum_{j=0}^{m}\binom{\alpha}{j}\binom{\beta}{m-j}
&=\binom{\alpha+\beta}{m}, \label{s-p-binomial-1}\\
\sum_{j=0}^{m}\binom{m}{j}\binom{\alpha}{j+n}
&=\binom{m+\alpha}{m+n}.\label{s-p-binomial-2}
\end{align}
\end{lemma}
\pf
The identity (\ref{s-p-binomial-1}) is well known and is obtained by
taking the coefficient of $x^{m}$ from both sides of
\begin{align*}
\sum_{m\in \N}\left(\sum_{j=0}^{m}\binom{\alpha}{j}
\binom{\beta}{m-j}\right)x^{m} =(1+x)^{\alpha}(1+x)^{\beta}
=(1+x)^{\alpha+\beta}= \sum_{m\in\N}\binom{\alpha+\beta}{m}x^{m}.
\end{align*}

We have
\begin{align}\label{s-p-binomial-2.1}
\sum_{j=0}^{m}\binom{m}{j}\binom{\alpha}{j+n}&=
\sum_{j=0}^{m}\binom{m}{m-j}\binom{\alpha}{j+n}\nn
&=\sum_{j=0}^{m}\binom{m}{(m+n)-(j+n)}\binom{\alpha}{j+n}\nn
&=\sum_{k=n}^{m+n}\binom{m}{(m+n)-k}\binom{\alpha}{k}.
\end{align}
Using $\binom{m}{(m+n)-k}=0$ for $k=0, \dots, n-1$ and 
(\ref{s-p-binomial-1}), we see that
the right-hand side of (\ref{s-p-binomial-2.1}) is equal to 
\begin{equation}\label{s-p-binomial-2.2}
\sum_{k=0}^{m+n}\binom{m}{(m+n)-k}\binom{\alpha}{k}
=\binom{m+\alpha}{m+n}.
\end{equation}
Combning  (\ref{s-p-binomial-2.1}) and (\ref{s-p-binomial-2.2}),
we obtain (\ref{s-p-binomial-2}).
\epfv

\begin{lemma}
For $m, n\in \N$ satisfying $m>n$ and $l=0, \dots, m$, 
\begin{equation}\label{main-lemma-13}
\sum_{k=0}^{n}\sum_{j=0}^{m}
\binom{-2m+n-1}{k}
\frac{l}{2m-j-n+k}\binom{m}{j}\binom{l-1}{2m-j-n+k-1}=\delta_{l, m-n}.
\end{equation}
\end{lemma}
\pf
In the case $\alpha\in \N$ and $\alpha\le m$, we have $\binom{\alpha}{j}=0$ 
for $\alpha< j\le m$
and (\ref{s-p-binomial-1}) becomes
\begin{equation}\label{main-lemma-12.5}
\sum_{j=0}^{\alpha}\binom{\alpha}{j}\binom{\beta}{m-j}=
\binom{\alpha+\beta}{m}.
\end{equation}
Note that $2m-n+k-1\ge m\ge s\ge l$.
For $k=0, \dots, n$ and $j, l=0, \dots, m$, we have
\begin{align}\label{main-lemma-12}
&\frac{l}{(2m-j-n+k)}\binom{l-1}{2m-j-n+k-1}\nn
&\quad=\frac{l}{(2m-j-n+k)}\frac{(l-1)\cdots
 (l-1-2m+j+n-k+2)}{(2m-j-n+k-1)!}\nn
 &\quad =\frac{l\cdots
 (l-2m+j+n-k+1)}{(2m-j-n+k)!}\nn
 &\quad=\binom{l}{2m-j-n+k}.
 \end{align}
 Using (\ref{main-lemma-12}),  (\ref{main-lemma-12.5})
with $\alpha=m$, $q=l$ and $m$ replaced by $2m-n+k$,
the fact that $\binom{m+l}{2m-n+k}=0$ when $l+n-m<0$,
and (\ref{s-p-binomial-1}) with $\alpha=-2m+n-1$, $\beta=m+l$ and
$m$ replaced by $l+n-m$, 
we have
\begin{align*}
&\sum_{k=0}^{n}\sum_{j=0}^{m}
\binom{-2m+n-1}{k}
\frac{l}{2m-j-n+k}\binom{m}{j}\binom{l-1}{2m-j-n+k-1}\nn
&\quad =\sum_{k=0}^{n}\sum_{j=0}^{m}
\binom{-2m+n-1}{k}
\binom{m}{j}\binom{l}{2m-j-n+k}\nn
&\quad=\sum_{k=0}^{n}
\binom{-2m+n-1}{k}
\binom{m+l}{2m-n+k}\nn
&\quad=\left\{\begin{array}{ll}0&l+n-m<0\\
\displaystyle \sum_{k=0}^{l+n-m}
\binom{-2m+n-1}{k}
\binom{m+l}{l+n-m-k}&l+n-m\ge 0\end{array}\right.\nn
&\quad=\left\{\begin{array}{ll}0&l+n-m<0\\
\displaystyle \binom{l+n-m-1}{l+n-m}&l+n-m\ge 0\end{array}\right.\nn
&\quad=\left\{\begin{array}{ll}0&l+n-m\ne 0\\
1&l+n-m=0.\end{array}\right.\nn
&\quad =\delta_{l, m-n}.
\end{align*}
\epfv

\begin{lemma}
An identity of Andersen (\cite{A}):
For $\alpha\in \C$, $n\in \Z_{+}$ and $k=0, \dots, n$,
\begin{equation}\label{andersen}
\sum_{j=0}^{k}\binom{\alpha}{j}\binom{-\alpha}{m-j}
=\frac{m-k}{m}\binom{\alpha-1}{k}\binom{-\alpha}{m-k}.
\end{equation}
\end{lemma}
\pf
For the proof of (\ref{andersen}), see \cite{A}. 
\epfv

\begin{lemma}
For $m, n, l\in \N$ satisfying $2n>m>n\le l$, 
\begin{align}\label{main-comb-formu}
\sum_{p=0}^{2n-m}(-1)^{-p-m}\sum_{k=0}^{m}
\binom{-2n+m-1}{k}\binom{2n-m+1}{p+m+1-k}
\binom{n+p+l}{p+m}=\delta_{l, m-n}.
\end{align}
\end{lemma}
\pf
Using the identity (\ref{andersen})
with $\alpha, m, k$ replaced by $-2n+m-1, 
p+m+1, m$
and 
$$(-1)^{-p-m}\binom{n+p+l}{p+m}=\binom{m-n-l-1}{p+m},$$
we see that the left-hand side of (\ref{main-comb-formu}) is equal to 
\begin{align}\label{main-comb-formu-1}
&\sum_{p=0}^{2n-m}(-1)^{-p-m}\frac{p+1}{p+m+1}\binom{-2n+m-2}{m}
\binom{2n-m+1}{p+1}
\binom{n+p+l}{p+m}\nn
&\quad =\binom{-2n+m-2}{m}\sum_{p=0}^{2n-m}\frac{2n-m+1}{p+m+1}
\binom{2n-m}{p}\binom{m-n-l-1}{p+m}.
\end{align}
When $l\ne m-n$, we have 
\begin{equation}\label{main-comb-formu-1.5}
\frac{1}{p+m+1}\binom{m-n-l-1}{p+m}
=\frac{1}{m-n-l}\binom{m-n-l}{p+m+1}.
\end{equation}
Using (\ref{main-comb-formu-1.5}) and (\ref{s-p-binomial-2})
with $m, n, \alpha$ replaced by $2n-m, m+1, m-n-l$, respectively, 
the right-hand side of (\ref{main-comb-formu-1}) is equal to 
\begin{align}\label{main-comb-formu-2}
&\frac{2n-m+1}{m-n-l}\binom{-2n+m-2}{m}\sum_{p=0}^{2n-m}
\binom{2n-m}{p}\binom{m-n-l}{p+m+1}\nn
&\quad =\frac{2n-m+1}{m-n-l}\binom{-2n+m-2}{m}\binom{n-l}{2n+1}.
\end{align}
Since $l\le n<2n+1$, $\binom{n-l}{2n+1}=0$. So the right-hand side of 
(\ref{main-comb-formu-2}) and also the right-hand side of
(\ref{main-comb-formu-1})
is $0$ in this case. Thus (\ref{main-comb-formu}) holds
in the case $l\ne m-n$. 
In the case $l=m-n$, using 
$$\binom{-2n+m-2}{m}\binom{-1}{p+m}
=(2n+1) \binom{2n}{m}(-1)^{p}$$
and (\ref{main-comb-formu-5}) with $n$ replaced by $2n-m$, 
we see  that
the right-hand side of  (\ref{main-comb-formu-1})
is equal to 
\begin{align*}
&\binom{-2n+m-2}{m}\sum_{p=0}^{2n-m}\frac{2n-m+1}{p+m+1}
\binom{2n-m}{p}\binom{-1}{p+m}\nn
&\quad=(2n+1)
\binom{2n}{m}\sum_{p=0}^{2n-m}\frac{(-1)^{p}}{p+m+1}
\binom{2n-m}{p}\nn
&\quad =1.
\end{align*}
Thus  (\ref{main-comb-formu}) also holds in the case $l=m-n$.
\epfv

\begin{lemma}
For $m, n\in \N$, 
\begin{equation}\label{main-comb-formu-5}
(n+1+m)\binom{n+m}{m}
\sum_{p=0}^{n}\frac{(-1)^{p}}{p+m+1}\binom{n}{p}=1.
\end{equation}
\end{lemma}
\pf
Multiplying $x^{m}$ to both sides of the binomial expansion
$$\sum_{p=0}^{n}\binom{n}{p}x^{p}=(1+x)^{n},$$
we obtain 
$$\sum_{p=0}^{n}\binom{n}{p}x^{p+m}=(1+x)^{n}x^{m}.$$
Integrating both sides from $0$ to $x$, we obtain 
\begin{align}\label{main-comb-formu-4}
&\sum_{p=0}^{n}\binom{n}{p}\frac{x^{p+m+1}}{p+m+1}\nn
&\quad =\int_{0}^{x}(1+t)^{n}t^{m}dt\nn
&\quad =\sum_{i=0}^{m}\frac{m\cdots (m-i+1)}{(n+1)\cdots (n+1+i)}
(1+x)^{n+1+i}x^{m-i}-(-1)^{m}\frac{m\cdots 1}{(n+1)\cdots (n+1+m)}.
\end{align}
Substituting $-1$ for $x$ in both sides of (\ref{main-comb-formu-4}),
we obtain 
\begin{equation}\label{main-comb-formu-6}
(-1)^{m+1}\sum_{p=0}^{n}\frac{(-1)^{p}}{p+m+1}\binom{n}{p}
=-(-1)^{m}\frac{m\cdots 1}{(n+1)\cdots (n+1+m)}.
\end{equation}
The identity (\ref{main-comb-formu-6}) is equivalent to 
(\ref{main-comb-formu-5}).
\epfv

\begin{lemma}
For $1\le l\le n$, 
\begin{align}
\sum_{j=0}^{n-1}\binom{n}{j}\frac{l}{n-j}\binom{-l-1}{n-j-1}&=1.
\label{2nd-main-lemma-6}\\
\sum_{j=0}^{n-1}\binom{n}{j}\frac{l}{n-j}\binom{l-1}{n-j-1}
&=\binom{n+l}{n}-1.\label{2nd-main-lemma-7}
\end{align}
\end{lemma}
\pf
For $1\le l\le n$, 
\begin{align}\label{2nd-main-lemma-4}
&\sum_{j=0}^{n-1}\binom{n}{j}\frac{l}{n-j}\binom{-l-1}{n-j-1}\nn
&\quad =\sum_{j=0}^{n-1}\binom{n}{j}\frac{l}{n-j}
\frac{(-l-1)\cdots (-l-1-n+j+1+1)}{(n-j-1)!}\nn
&\quad =\sum_{j=0}^{n-1}(-1)^{n-j-1}\binom{n}{j}
\frac{(l+n-j-1)\cdots (l+1)l(l-1)!}{(l-1)!(n-j)!}\nn
&\quad =\sum_{j=0}^{n-1}(-1)^{n-j-1}\binom{n}{j}
\frac{(l+n-j-1)\cdots (n-j+1)}{(l-1)!}\nn
&\quad =\sum_{j=0}^{n-1}(-1)^{n-j-1}\binom{n}{j}
\frac{1}{(l-1)!}\frac{d^{l-1}}
{dx^{l-1}}x^{l+n-j-1}\lbar_{x=1}.
\end{align}
But 
\begin{align*}
&\sum_{j=0}^{n}(-1)^{n-j-1}\binom{n}{j}
\frac{1}{(l-1)!}\frac{d^{l-1}}
{dx^{l-1}}x^{l+n-j-1}\lbar_{x=1}\nn
&\quad =(-1)^{n-1}
\frac{1}{(l-1)!}\frac{d^{l-1}}
\sum_{j=0}^{n}\binom{n}{j}(-1)^{j}
{dx^{l-1}}x^{l+n-j-1}\lbar_{x=1}\nn
&\quad =(-1)^{n-1}
\frac{1}{(l-1)!}\frac{d^{l-1}}
{dx^{l-1}}x^{l+n-j-1}(1-x^{-1})^{n}\lbar_{x=1}\nn
&\quad =(-1)^{n-1}
\frac{1}{(l-1)!}\frac{d^{l-1}}
{dx^{l-1}}x^{l-j-1}(x-1)^{n}\lbar_{x=1}\nn
&\quad =0.
\end{align*}
Then
\begin{align}\label{2nd-main-lemma-5}
&\sum_{j=0}^{n-1}(-1)^{n-j-1}\binom{n}{j}
\frac{1}{(l-1)!}\frac{d^{l-1}}
{dx^{l-1}}x^{l+n-j-1}\lbar_{x=1}\nn
&\quad =\sum_{j=0}^{n}(-1)^{n-j-1}\binom{n}{j}
\frac{1}{(l-1)!}\frac{d^{l-1}}
{dx^{l-1}}x^{l+n-j-1}\lbar_{x=1}+1\nn
&\quad= 1.
\end{align}
From (\ref{2nd-main-lemma-4}) and (\ref{2nd-main-lemma-5}),
we obtain (\ref{2nd-main-lemma-6}).

For $1\le l\le n$, 
\begin{align*}
\sum_{j=0}^{n-1}\binom{n}{j}\frac{l}{n-j}\binom{l-1}{n-j-1}
&=\sum_{j=0}^{n-1}\binom{n}{j}\binom{l}{n-j}\nn
&=\sum_{j=0}^{n}\binom{n}{j}\binom{l}{n-j}-1\nn
&=\binom{n+l}{n}-1,
\end{align*}
where in the last step, we have used (\ref{s-p-binomial-1}). 
\epfv

\begin{lemma}
For $1\le l\le n$, 
\begin{align}
\sum_{m=1}^{n}\sum_{j=0}^{n}
\binom{-n-1}{m}\binom{n}{j}\frac{l}{n-j+m}
\binom{-l-1}{n-j+m-1} &=0,\label{2nd-main-lemma-8}\\
\sum_{m=1}^{n}\sum_{j=0}^{n}
\binom{-n-1}{m}\binom{n}{j}\frac{l}{n-j+m}
\binom{l-1}{n-j+m-1} &=-\binom{n+l}{n}.\label{2nd-main-lemma-9}
\end{align}
\end{lemma}
\pf
For $1\le l\le n$, by (\ref{lemma-3-1}), we have 
\begin{align*}
&\sum_{m=1}^{n}\sum_{j=0}^{n}
\binom{-n-1}{m}\binom{n}{j}\frac{l}{n-j+m}
\binom{-l-1}{n-j+m-1}\nn
&\quad=\sum_{m=1}^{n}\binom{-n-1}{m}\sum_{j=0}^{n}
\binom{n}{j}\frac{l}{n-j+m}
\binom{-l-1}{n-j+m-1}\nn
&\quad =0.
\end{align*}

For $1\le l\le n$, using  the properties of binomial coefficients and  
(\ref{s-p-binomial-1}), we have
\begin{align*}
&\sum_{m=0}^{n}\sum_{j=0}^{n}
\binom{-n-1}{m}\binom{n}{j}\frac{l}{n-j+m}
\binom{l-1}{n-j+m-1}\nn
&\quad =\sum_{m=0}^{n}\sum_{j=0}^{n}
\binom{-n-1}{m}\binom{n}{j}
\binom{l}{n-j+m}\nn
&\quad =\sum_{m=0}^{n}
\binom{-n-1}{m}\sum_{j=0}^{n+m}\binom{n}{j}
\binom{l}{n-j+m}\nn
&\quad =\sum_{m=0}^{n}
\binom{-n-1}{m}\binom{n+l}{n+m}\nn
&\quad =\sum_{m=0}^{l}
\binom{-n-1}{m}\binom{n+l}{l-m}\nn
&\quad =\sum_{m=0}^{l}
\binom{-n-1}{m}\binom{n+l}{l-m}\nn
&\quad =\binom{l-1}{l}\nn
&\quad=0.
\end{align*}
Then 
\begin{align*}
&\sum_{m=1}^{n}\sum_{j=0}^{n}
\binom{-n-1}{m}\binom{n}{j}\frac{l}{n-j+m}
\binom{l-1}{n-j+m-1}\nn
&\quad =\sum_{m=0}^{n}\sum_{j=0}^{n}
\binom{-n-1}{m}\binom{n}{j}\frac{l}{n-j+m}
\binom{l-1}{n-j+m-1} -\sum_{j=0}^{n}
\binom{n}{j}\frac{l}{n-j}
\binom{l-1}{n-j-1}\nn
&\quad =-\sum_{j=0}^{n}
\binom{n}{j}
\binom{l}{n-j}\nn
&\quad =-\binom{n+l}{n},
\end{align*}
where in the last step, we have used (\ref{s-p-binomial-1}). 
\epfv

\noindent {\small \sc Department of Mathematics, Rutgers University,
110 Frelinghuysen Rd., Piscataway, NJ 08854-8019}

\noindent {\em E-mail address}: {\tt yzhuang@math.rutgers.edu}


\begin{thebibliography}{KWAK2}

\bibitem[ABD]{ABD}
T. Abe, G. Buhl, C. Dong, Rationality, Regularity,
 and
$C_2$-cofiniteness, {\it Trans. Amer. Math. Soc.} {\bf 356}
(2004), 3391--3402. 

\bibitem[AbN]{AbN}
T. Abe and K. Nagatomo, Finiteness of conformal blocks over the
projective line, in: {\it Vertex operator algebras in mathematics and physics 
(Toronto, 2000)}, ed. S. Berman, Y. Billig, Y.-Z. Huang and 
J. Lepowsky, Fields Inst. Commun., Vol. 39, 
Amer. Math. Soc., Providence, 2003, 1--12.

\bibitem[AM]{AM}
D. Adamovi\`{c} and A. Milas, On $\mathcal{W}$-algebras 
associated to $(2, p)$ minimal
models and their representations, {\it Int. Math. Res. Not.}
{\bf 20} (2010), 3896--3934.


\bibitem[An]{A}
E. S. Andersen,  Two summation formulae for product sums of 
binomial coefficients,  {\it Math. Scand.} {\bf 1} (1953), 261--262. 

%\bibitem[AF]{AF} 
%F.~ Anderson and K.~Fuller, {\it Rings and Categories
%of Modules}, Graduate Texts in Mathematics 13, Springer-Verlag,
%New York, Heidelberg, Berlin (1991).

\bibitem[Ar]{Ar} 
Y.~Arike, Some remarks on symmetric linear functions
and pseudotrace maps, {\em Proc. Japan Acad. Ser. A Math. Sci.} {\bf 86}
(2010), 119-124.


\bibitem[ArN]{ArN} 
Y.~Arike and K.~Nagatomo, Some remarks on pseudo-trace 
functions for orbifold models associated with symplectic fermions, 
{\em Int. J. Math.} {\bf 24}, (2013) 1350008 

%\bibitem[Br]{Br} M.~Brou\'e, Higman criterion revisited, {\em Mich. J.
%Math.} {\bf 58}, (2009) 125-179.

\bibitem[B]{B} 
G. Buhl, A spanning set for VOA modules, {\it J. Alg.} {\bf 254} (2002), 125--151.

\bibitem[DLM1]{DLM}
C. Dong, H. Li and G. Mason, Vertex operator algebras and
associative algebras, {\em J. Algebra} {\bf 206} (1998), 67-96.

\bibitem[DLM2]{DLM2} 
C. Dong, H. Li and G. Mason, Modular
invariance of trace functions in orbifold theory 
 and generalized moonshine, {\it
Comm. Math. Phys.} {\bf 214} (2000), 1--56.

\bibitem[F1]{F1} 
F.  Fiordalisi, Logarithmic intertwining operator
and genus-one correlation functions, Ph.D. Thesis,
Rutgers University, 2015. 
(\href{https://rucore.libraries.rutgers.edu/rutgers-lib/47367/PDF/1/play/}{Online PDF file} of 
the thesis.)

\bibitem[F2]{F}
F. Fiordalisi, Logarithmic Intertwining Operators and Genus-One Correlation Functions,
{\it Comm. Contemp. Math.} {\bf 18} (2016), 1650026. 

\bibitem[GN]{GN}
M. R. Gaberdiel and A. Neitzke, Rationality, quasirationality
and finite W-algebras, {\it Comm. Math.Phys.} {\bf 238} (2003) 305--331.

\bibitem[H1]{H-diff-eqn} 
Y.-Z. Huang, Differential equations and
intertwining operators, {\em Comm. Contemp. Math.} {\bf 7} (2005),
375--400.

\bibitem[H2]{H-mod-inv-int}
Y.-Z. Huang, Differential equations, duality and modular invariance, 
{\it Comm. Contemp. Math.} {\bf 7} (2005), 649--706. 

\bibitem[H3]{H-verlinde}
Y.-Z. Huang,
Vertex operator algebras and the Verlinde conjecture, 
{\it Comm. Contemp. Math.} {\bf 10} (2008), 103--154. 

\bibitem[H4]{H-rigidity}
Y.-Z. Huang,
Rigidity and modularity of vertex tensor categories, {\it Comm. Contemp. Math.}
{\bf 10} (2008), 871--911. 

\bibitem[H5]{H-cofiniteness} 
Y.-Z. Huang, Cofiniteness conditions, projective
covers and the logarithmic tensor product theory, {\em J. Pure
Appl. Alg.} {\bf 213} (2009), 458--475.

\bibitem[H6]{H-const-twisted-mod}
Y.-Z. Huang, A construction of lower-bounded generalized 
twisted modules for a grading-restricted vertex (super)algebra,
{\it Comm. Math. Phys.} {\bf 377} (2020), 909--945. 

\bibitem[H7]{H-affine-twisted-mod}
Y.-Z. Huang, Lower-bounded and grading-restricted twisted modules for affine vertex (operator) algebras, {\it J. Pure Appl. Alg.} {\bf 225} (2021), Paper no. 106618. 

\bibitem[H8]{H-aa-va}
Y.-Z. Huang, Associative algebras and the representation theory of 
grading-restricted vertex algebras, 
to appear; arXiv:2009.00262. 

\bibitem[H9]{H-aa-int-op}
Y.-Z. Huang, Associative algebras and intertwining operators, 
{\it Comm. Math. Phys.} {\bf 396} (2022), 1--44. 

\bibitem[H10]{H-conv-cft}
Y.-Z. Huang, 
Convergence in conformal field theory, {\it Chin. Ann. Math.} {\bf B43} 
(2022), 1101--1124. 

\bibitem[HLZ1]{HLZ6}
Y.-Z. Huang, J. Lepowsky and L. Zhang, Logarithmic tensor category theory, VI:
Expansion condition, associativity of logarithmic intertwining operators, and the
associativity isomorphisms; arXiv:1012.4202.

\bibitem[HLZ2]{HLZ7}
Y.-Z. Huang, J. Lepowsky and L. Zhang, Logarithmi tensor category theory, VII:
Convergence and extension properties and applications to expansion for intertwining
maps; arXiv:1110.1929.

\bibitem[HY]{HY}
Y.-Z. Huang and J. Yang,  Logarithmic intertwining operators and associative algebras, 
{\it J. Pure Appl. Alg.} {\bf 216} (2011), 1467--1492. Corrigendum, {\it J. Pure Appl.
Alg.} {\bf 226} (2021), 107020.


\bibitem[K]{K}
N. Koblitz, {\it Introduction to elliptic
curves and modular forms}, Second Edition, Graduate Texts in Mathematics, 
Vol. 97, Springer-Verlag, New York, 1993. 


\bibitem[L]{L}
S. Lang, {\it Elliptic functions}, Graduate Texts in Mathematics, 
Vol. 112, Springer-Verlag, New York, 1987.

\bibitem[Mc]{Mc} 
R. McRae, Proof of Miyamoto's Proposition 4.4, private communication. 

\bibitem[MS1]{MS1}
G.~Moore and N.~Seiberg,
Polynomial equations for rational conformal field theories,
{\em Phys. Lett.} {\bf B212} (1988), 451--460.

\bibitem[MS2]{MS2}
G.~Moore and N.~Seiberg,
Classical and quantum conformal field theory,
{\em Comm. Math. Phys.} {\bf 123} (1989), 177--254.

\bibitem[Mi1]{M1} M.~Myiamoto, 
Intertwining operators and modular invariance, 
arXiv:math/0010180.

\bibitem[Mi2]{M} M.~Myiamoto, Modular invariance of vertex operator
algebras satisfying $C_2$-cofiniteness, {\em Duke Math. J.} {\bf 122}
(2004), 51-91.

\bibitem[T]{T}
V. G. Turaev, {\em Quantum invariants of knots and $3$-manifolds},
de Gruyter Studies in Math., Vol. 18, 
Walter de Gruyter, Berlin, 1994.

\bibitem[V]{V}
E. Verlinde, 
Fusion rules and modular transformations in 2D conformal field theory, 
{\em Nucl. Phys.} {\bf B300} (1988), 360--376.

\bibitem[Z]{Z}
Y. Zhu, Modular invariance of characters of vertex operator algebras,
{\em J.
Amer. Math. Soc.} {\bf 9} (1996), 237--307.

\end{thebibliography}
\end{document}